\documentclass[11 pt, reqno]{article}

\usepackage{amsmath,amsfonts,amssymb,amsthm}
\usepackage[colorlinks=true,hyperindex=true]{hyperref}
\usepackage{cancel,bbm}
\usepackage{mathabx} 
\usepackage{comment}
\usepackage{enumitem}
\usepackage{cite}

\usepackage{chngcntr}
\counterwithin*{equation}{section}

\usepackage[left=3cm, right=3cm, top=3 cm, bottom= 3 cm]{geometry}
\usepackage{color}
\usepackage{fancyhdr}
\usepackage{latexsym}

\usepackage{bm}

\newtheorem{ccounter}{ccounter}[section]
\newtheorem{thm}[ccounter]{Theorem}
\newtheorem{lem}[ccounter]{Lemma}
\newtheorem{cor}[ccounter]{Corollary}
\newtheorem{defn}[ccounter]{Definition}
\newtheorem{prop}[ccounter]{Proposition}
\newtheorem{ass}[ccounter]{Assumption}
\newtheorem{ex}[ccounter]{Example}

\def\bet{\begin{thm}}
\def\eet{\end{thm}}
\def\bel{\begin{lem}}
\def\eel{\end{lem}}
\def\bas{\begin{ass}}
\def\eas{\end{ass}}
\def\bec{\begin{cor}}
\def\eec{\end{cor}}
\def\bed{\begin{defn}}
\def\eed{\end{defn}}
\def\bep{\begin{prop}}
\def\eep{\end{prop}}
\def\beq{\begin{equation}}
\def\eeq{\end{equation}}
\def\proof{\noindent {\bf Proof.}\ \ }
\def\bea{\begin{equation*}}
\def\eea{\end{equation*}}
\def\tr{\mathrm{tr}}
\def\bex{\begin{ex}}
\def\eex{\end{ex}}
\def\remark{\noindent{\bf Remark. }}
\usepackage{tikz}
\usepackage{caption}

\def\rr{\mathbb{R}}

\def\cc{\mathbb{C}}
\def\1{\boldsymbol{1}}
\def\Im{\mathrm{Im}}
\def\Re{\mathrm{Re}}
\def\e{\mathrm{e}}
\def\i{\mathrm{i}}
\def\del{\partial}
\def\d{\mathrm{d}}
\def\eps{\varepsilon}
\renewcommand\leq\varleq
\renewcommand\geq\vargeq

\def\O{\mathcal{O}}

\def\pp{\mathbb{P}}

\def\D{\mathcal{D}}
\def\I{\mathcal{I}}
\def\mfa{\mathfrak{m}}
\def\A{\mathcal{A}}

\def\dto{\downarrow}
\def\uto{\uparrow}
\def\D{\mathcal{D}}
\def\C{\mathcal{C}}

\def\tilphi{\tilde{\varphi}}
\def\mfa{\mathfrak{a}}

\def\tilf{\tilde{f}}

\def\mfb{\mathfrak{b}}

\def\B{\mathcal{B}}
\def\C{\mathcal{C}}
\def\E{\mathcal{E}}
\def\tilphi{\tilde{\varphi}}
\def\tileta{\tilde{\eta}}
\def\mfc{\mathfrak{c}}
\def\tiltheta{\tilde{\theta}}
\def\fB{\mathfrak{B}}
\def\tilgamma{\tilde{\gamma}}
\def\uu{\mathbb{U}}
\def\tilm{\tilde{m}}
\def\mfd{\mathfrak{d}}
\def\Ed{\mathcal{G}}

\begin{document}
\title{Unitary BM}

\begin{table}
\centering

\begin{tabular}{c}
\multicolumn{1}{c}{\Large{\bf Local law and rigidity for unitary Brownian motion}}\\
\\
\\
\end{tabular}
\begin{tabular}{ c c c  }
Arka Adhikari$^1$
& \phantom{blah} & 
Benjamin Landon$^2$ 
 \\
 & & \\  
 \small{Stanford University} & & \small{University of Toronto } \\
 \small{Department of Mathematics} & & \small{Department of Mathematics} \\
 \small{\texttt{arkaa@stanford.edu}} & & \small{\texttt{blandon@math.toronto.edu}} \\
  & & \\
\end{tabular}
\\
\begin{tabular}{c}
\multicolumn{1}{c}{\today}\\
\\
\end{tabular}

\begin{tabular}{p{15 cm}}
\small{{\bf Abstract:}  We establish high probability estimates on the eigenvalue locations of  Brownian motion on the $N$-dimensional unitary group, 
 as well as estimates on the number of eigenvalues lying in any interval on the unit circle. These estimates are optimal up to arbitrarily small polynomial factors in $N$. Our results hold at the spectral edges (showing that the extremal eigenvalues are within $\O (N^{-2/3+})$ of the edges of the limiting spectral measure), in the spectral bulk, as well as for times near $4$ at which point the limiting spectral measure forms a cusp. Our methods are dynamical and are based on analyzing the evolution of the Cauchy transform of the empirical spectral measure along the characteristics of the PDE satisfied by the limiting spectral measure, that of the free unitary Brownian motion.}
\end{tabular}
\end{table}

\section{Introduction and main results}

{\let\thefootnote\relax\footnotetext{$1$. The research of A.A. is supported by NSF grant DMS-2102842. $2$. The work of B.L. is partially supported by NSERC.}}In this work we study Brownian motion on the unitary group $\uu(N)$ of dimension $N$.  One can define Brownian motion on $\uu(N)$ by considering the left-invariant Riemannian metric induced by the inner product $\langle A, B \rangle = N \tr (A B^*)$ on the Lie algebra of skew-Hermitian matrices. Unitary Brownian motion  is then the Markov diffusion process starting from the identity matrix with generator given by the Laplacian on $\uu (N)$ associated to this metric.

It will be more convenient for us to consider the following equivalent definition of $U_t$ as the solution of the  It{\^{o}} stochastic differential equation,
\beq \label{eqn:unit-def}
\d U_t = \i U_t \d W_t - \frac{1}{2} U_t \d t, \qquad U_0 = \1
\eeq
where $W_t$ is a standard complex Hermitian Brownian motion. That is, if $X_t$ and $X'_t$ are $N \times N$ matrices of independent standard Brownian motions, then,
\beq \label{eqn:wt-def}
W_t = \frac{1}{ \sqrt{4N}} \left( X_t +X_t^T  + \i (X'_t - (X'_t)^T ) \right).
\eeq
The process \eqref{eqn:unit-def} admits strong solutions by standard results (see, e.g., Theorem 8.3 of \cite{legall}).

 Unitary Brownian motion is well-studied in random matrix theory as well as in the context of free probability  due to its connection with an object called free unitary Brownian motion.  Of particular interest is the  empirical spectral measure  of $U_t$, which is a random, time-dependent measure on the unit circle defined by,
 \beq
 \d \nu_{N, t} (x) := \frac{1}{N} \sum_{i=1}^N \delta_{\lambda_i (t) } (x) \d x.
 \eeq
 In the work \cite{biane1997free}, Biane showed that for fixed $t$, the measure $\nu_{N, t}$ converges almost surely to a measure on the unit circle which we will denote by $\nu_t$. Identifying the unit circle with the angular coordinates $\theta \in (-\pi, \pi]$, the measure $\nu_t$ has a density $\rho_t ( \theta ) $ for any $ t>0$. For $t <4$ the support of $\rho_t$ is given by,
 \beq \label{eqn:int-thetat}
I_t := [ - \Theta_t, \Theta_t] , \qquad \Theta_t :=  \frac{1}{2} \sqrt{ (4-t)t} +2 \arcsin \left( \sqrt{\frac{t}{4}}\right) ,
 \eeq
 whereas for $ t \geq 4$, the support is the entire unit circle. Moreover, for $t \geq 4$ the density is everywhere non-zero unless $t=4$ in which case $\rho_t$ vanishes only at $\pi$.  In fact, $\rho_t$ can be described as the spectral measure of  free unitary Brownian motion, an object appearing in free probability.  The limit of the $\nu_{N, t}$ was also derived independently by Rains in \cite{rains1997combinatorial}.
 
 Since Biane's paper, there have been many works studying the convergence of $\nu_{N, t}$ to $\rho_t$.  Concentration estimates and convergence for the empirical averages $\int f  \d \nu_{N, t}$ were established by Kemp \cite{kemp2017heat} for various classes of $f$ of low regularity. Meckes and Melcher \cite{meckes2018convergence} established explicit convergence rates in terms of the $L^1$-Wasserstein metric. For $0 < t < 4$, the convergence of the spectral edge of $U_t$ to $\pm \Theta_t$ was established by Collins, Dahlqvist and Kemp \cite{collins2018spectral}.  This work also established a multi-time, multi-matrix version of this result. The asymptotic Gaussian fluctuations of the empirical averages of $\int f \d \nu_{N, t}$ were established by L{\'e}vy and Maida \cite{levy2010central}. Multivariate fluctuations for trace polynomials of a two parameter family of diffusion processes (including unitary Brownian motion as a special case) were studied by C{\'e}bron and Kemp \cite{cebron2022fluctuations}.

The main contribution of the present work is to establish almost-optimal rates (i.e., up to polynomial $N^{\eps}$ factors) of convergence of $\nu_{N, t}$ to the limiting distribution $\rho_t$ on the almost-shortest possible scales, as well as almost-optimal estimates on the eigenvalue locations. In the random matrix literature, these estimates are known as \emph{local laws} and \emph{rigidity} estimates, respectively.  Our local laws are stated in Theorem \ref{thm:bulk} and in Corollary \ref{cor:interval} below, showing that the number of eigenvalues in any sub-interval $I$ of the unit circle is given by $N \rho_t (I) + \O (N^{\eps} )$ for any $\eps >0$.  

For times $ t< 4$ we establish almost-optimal rates of convergence of the spectral edge of $U_t$ to $\pm \Theta_t$. That is, the extremal eigenvalues are within distance $\O (N^{-2/3+\eps} )$ of $\pm \Theta_t$ (what is usually termed \emph{edge rigidity} in the literature). Given the square-root behavior of the spectral measure at the edges, this is expected to be optimal up to the $N^{\eps}$ factor.  We also derive almost-optimal edge rigidity results up to $t = 4 - N^{-1/2+\eps}$ at which point the measure $\rho_t$ forms a cusp (i.e., vanishes like a cube root) near $\theta = \pm \pi$. This is expected to be optimal as at later times, the natural inter-particle distance at the spectral edges exceeds the distance between $+ \Theta_t$ and $-\Theta_t$.  Our rigidity estimates  are formulated in Corollaries \ref{cor:edge-rig} and \ref{cor:cusp-rig} below. The scaling behaviors in the various parameter regimes will be given as the results are introduced. 

To our knowledge, our results are the strongest available estimates on eigenvalue locations for unitary Brownian motion. Our methods are completely different from prior works on rates of convergence to $\rho_t$, relying on the method of characteristics from PDEs. Previous works were based on moment calculations and/or concentration estimates for heat kernels on Lie groups.

\subsection{Discussion of methodology}

There are many approaches in the literature for proving local laws and rigidity in random matrix theory. For Wigner matrices and related mean-field random matrices $W$, a multi-scale approach to analyzing the resolvent $(W-z)^{-1}$ giving estimates down to the optimal scale was developed by Erd\H{o}s, Schlein and Yau \cite{erdHos2009local,erdHos2009semicircle,erdHos2010wegner}. For  pedagogical overviews of this strategy see \cite{semi-lec,erdos2017dynamical}. Earlier local laws for random matrices on short scales appeared in \cite{guionnet2000concentration,bai2002convergence}.   

While powerful, the resolvent method makes heavy use of the matrix structure and independence between different entries; while $U_t$ is associated with a matrix process, the correlation structure between entries is complicated and renders this approach intractable.  Instead, the process $U_t$ somewhat resembles what is known as (Hermitian) Dyson Brownian motion, whose definition we now recall. In the work \cite{dyson1962brownian}, Dyson showed that the eigenvalues of $V+W_t$ for any Hermitian $V$ (and $W_t$ as in \eqref{eqn:wt-def}) obey the system of SDEs,
\beq \label{eqn:dbm}
\d \mu_i  = \sqrt{ \frac{2}{N \beta} } \d B_i + \frac{1}{N} \sum_{j \neq i } \frac{1}{ \mu_i - \mu_j} \d t - \frac{1}{2} \mu_i \d t ,
\eeq
where $\beta=2$. Moreover, at the level of formal calculation, Dyson showed that the eigenvalues $ \lambda_i (t)$ of the matrix $U_t$ satisfy the closed system of SDEs,
\beq \label{eqn:unit-dbm}
\d \lambda_i  = \frac{1}{ \sqrt{N}} \i \lambda_i \d B_i - \frac{1}{N} \sum_{j \neq i } \frac{ \lambda_j  \lambda_i}{ \lambda_i  - \lambda_j } \d t - \frac{1}{2} \lambda_i \d t,
\eeq
where the $\{ B_i (t) \}_{i=1}^N$ are a family of independent standard Brownian motions. We will not make use of this SDE, but mention that well-definedness of this process was studied in \cite{cepa2001brownian}.

  Given the similar form of these two processes, it is therefore useful to recall how rigidity has been established in the study of  DBM \eqref{eqn:dbm} for general $\beta$.  Given Dyson's original derivation of \eqref{eqn:dbm} we see that for the special values $\beta=1, 2, 4$,  this eigenvalue process in fact comes from a matrix-valued Brownian motion that can be thought of as an additively deformed Gaussian Orthogonal/Unitary/Symplectic ensemble. In the case that the initial data is the $0$ matrix, this is just a scaled Gaussian matrix, and so the local laws and rigidity follow from those for Wigner matrices. For general initial data, local laws and rigidity were developed by Lee and Schnelli \cite{lee2016extremal,lee2013local} as well as the second author with Yau \cite{convergence,landon2017edge} still using resolvent methods, as the additive structure inherent in DBM allows this approach to work.  However,  for general $\beta$, the process \eqref{eqn:dbm} is no longer naturally associated with any matrix process and so no resolvent methods can work.

Nonetheless, the second author with Huang \cite{huang2019rigidity} showed that local laws and rigidity in fact hold for DBM with general $\beta$ as well as a general potential (the equation \eqref{eqn:dbm} being associated with quadratic $V' ( \mu) = \frac{1}{2} \mu$).  This was further developed by the first author with Huang to study the spectral edges \cite{adhikari2020dyson}.  The method in these works is our starting point and so we review it here.

The  works \cite{adhikari2020dyson,huang2019rigidity} were based on a PDE-style approach to studying the Stieltjes transform,
\beq
m (z, t) := \frac{1}{N} \sum_{i=1}^N \frac{1}{ \mu_i (t) -z } ,
\eeq
of the empirical eigenvalue measure associated to the Hermitian DBM \eqref{eqn:dbm}, on short scales $\Im[z] \sim N^{-1}$. 
For $\beta =2$, this quantity satisfies,
\beq
 \d m (z, t) = \left( m (z, t) + \frac{z}{2} \right) \del_z m(z, t) \d t+ \frac{1}{2} m (z, t) \d t  + \d N_t ,
\eeq
for some Martingale $N_t$ which turns out to be lower-order.  This leads to the limiting complex Burger's equation,
\beq \label{eqn:complex-burger}
\del_t \tilm (z, t) = \left( \tilm (z, t) + \frac{z}{2} \right) \del_z \tilm (z, t) + \frac{1}{2} \tilm (z, t) .
\eeq
These equations were considered by Rogers and Shi \cite{rogers1993interacting} and general potential analogues by Li, Li and Xie \cite{li2020law,li2013generalized}.  The complex Burger's equation may be solved by the elementary  method of characteristics. The works \cite{adhikari2020dyson,huang2019rigidity} were based on tracking the difference between $\tilm(z_t, t)$ and $m(z_t, t)$ along characteristics $z_t$. Note that for any fixed $z$, it is relatively straightforward to see that the limiting $m(z, t)$ must satisfy \eqref{eqn:complex-burger}. The main challenge in \cite{adhikari2020dyson,huang2019rigidity} is to extend this to short scales and obtain optimal error estimates.

The natural analogue of the Stieltjes transform $m (z, t)$ above for measures supported on the unit circle is the following transform,
\beq
f(z, t) := \frac{1}{N} \tr \left( \frac{ U_t + z}{ U_t-z} \right) = \int \frac{ \e^{ \i \theta} + z}{ \e^{ \i \theta } - z} \d \nu_{N, t}  ( \theta ) .
\eeq
We will refer to this as the ``Cauchy transform'' although the Cauchy transform is usually defined slightly differently \cite{cima2006cauchy} (the usual definition differs from ours only by an affine transformation). Our choice of $f$ is to match the work of Biane \cite{biane1997free,biane1997segal} discussed in more detail below.

An application of It{\^{o}}'s formula starting from \eqref{eqn:unit-def} shows that $f(z, t)$ obeys the SDE,
\beq \label{eqn:df}
\d f(z, t) = - \frac{z f(z, t)}{2} \del_z f(z, t) \d t + \d M_t (z)
\eeq
where $M_t(z)$ is a complex-valued Martingale defined via,
\beq \label{eqn:M-def}
\d M_t(z) = - \frac{2 \i z}{N} \tr \left( \frac{ U_t}{ (U_t-z)^2} \d W_t \right) = - \frac{2 \i z}{N} \sum_{i, j=1}^N \left[ \frac{ U_t}{ (U_t-z)^2} \right]_{ij} \d (W_t)_{ij} ,
\eeq
where the second equality explicitly writes out the trace in terms of the matrix elements of $U_t (U_t -z)^{-2}$ and the matrix of stochastic differentials of the Hermitian Brownian motion $W_t$. 
 The covariation process of $M_t$ is easily calculated,
 \begin{align} \label{eqn:M-cov}
 \langle \d M , \d M \rangle &= - \frac{4 z^2}{N^3} \tr \left( \frac{ U_t^2}{ ( U_t -z )^4} \right) \d t \notag \\
 \langle \d M , \d \bar{M} \rangle &= \frac{4 |z|^2}{N^3} \tr \left( \frac{ |U_t|^2}{ |U_t - z|^4 } \right) \d t.
 \end{align}
Based on this, one sees that the limiting Cauchy transform should be the solution to,
\beq \label{eqn:unit-pde}
\del_t \tilf (z, t) = - \frac{z \tilf (z, t)}{2} \del_z \tilf (z, t), \qquad \tilf (z, 0) = \frac{ 1+z}{1-z} ,
\eeq
the analogue of the complex Burger's equation \eqref{eqn:complex-burger}. 
Indeed, these equations were found by Biane \cite{biane1997segal,biane1997free}, and the Cauchy transform  $\tilf(z, t)$ characterizes the limiting spectral measure $\nu_t$.  
 We collect properties of this measure and its density $\rho_t ( \theta)$ in Appendix \ref{a:dos}.

In the present work, we will analyze the equation \eqref{eqn:df} through the characteristics associated to \eqref{eqn:unit-pde}. That is, if $z (t)$ is a time-dependent curve in $\cc$ satisfying,
\beq \label{eqn:dzdt}
\frac{ \d}{\d t} z(t) = \frac{ z(t) \tilf ( z(t), t) }{2},
\eeq
then,
\beq
\frac{ \d }{ \d t} \tilf ( z(t), t) = 0.
\eeq
In fact, Biane's work \cite{biane1997segal} shows that for any $t$ the map $z \to z(t)$ is a conformal map of some domain onto the complement of the unit circle in $\cc$, which allows one to construct $\tilf (z, t)$ from the characteristics. 

At a superficial level, we are translating the methods of \cite{huang2019rigidity,adhikari2020dyson} from the real line to the unit circle. However, the translation is not at all straightforward and there are serious obstacles to be overcome once one moves to our new setting. We briefly mention a few now; they will be presented in more detail below as we discuss our results. First, it is not a-priori clear that this method could even work in the first place. In particular, it is crucial that the Martingale term $M_t (z)$ above can be controlled by the empirical Cauchy transform $f (z, t)$ itself. The precise form of the quadratic variation of $M_t$ is therefore important. Similar considerations hold for controlling the term $\del_z f (z, t)$ by $f (z, t)$. In fact, getting precise constants here is crucial due to the use of Gronwall's inequality in our proof.

One of the new novelties of our work is to deal with times close to $t \approx 4$ where the spectral measure $\rho_t( \theta)$ forms a cusp singularity near $\theta = \pm \pi$.  Local laws near a cusp are in general delicate (see e.g., \cite{alt2018dyson,erdHos2020cusp,cipolloni2019cusp}) and in the random matrix setting have not been dealt with via the characteristics approach before our work (but see the work \cite{huang2021edge} which studies non-intersecting random walks via characteristics in a different context). Secondly, the works  \cite{adhikari2020dyson,huang2019rigidity} dealt only with short times $t = o(1)$ whereas we are interested in times of order $1$. Coupled with the curvature of the unit circle, this requires a more detailed understanding of the behavior of the characteristics, especially near the spectral edges and cusps than was required before. Here, we partially rely on the semi-explicit form of the spectral measure and Cauchy transform $\tilf (z, t)$.  Finally, the results of \cite{adhikari2020dyson} make somewhat strong assumptions on the initial data. This was primarily due to the treatment of general potentials in that work, which allowed for an analysis of the movement of the spectral edge. Here, our initial data is a delta function, falling outside the assumptions of \cite{adhikari2020dyson}. In particular, our short-time analysis is more involved.

The characteristic approach has appeared in a few other works on the short-scale behavior of eigenvalues. Bourgade used characteristics to analyze a ``stochastic advection equation'' derived from a certain coupling between DBMs and obtained fine estimates on the local eigenvalue behavior of general Wigner matrices, including universality of the extreme gaps \cite{bourgade2018extreme}.  Von Soosten and Warzel used random characteristics  to prove local laws for Wigner matrices \cite{von2019random} and to study delocalization in the {R}osenzweig--{P}orter model \cite{von2019non}.

\subsection{Statement of main results}

Cauchy transforms of measures $\mu$ on the unit circle $f_\mu (z)$ are traditionally studied for $z$ in the open unit disc, $\{ |z| < 1 \}$. Due to the identity,
\beq
f_\mu (r \e^{ \i \theta } ) = - \widebar{  f_\mu }( r^{-1} \e^{ \i \theta } ) ,
\eeq
this is equivalent to studying the behavior for $|z| >1$, and it is somewhat conceptually simpler for our techniques to treat $|z| >1$. 
Based on this, our first main result is the following, identifying the rate of convergence of the empirical spectral measure of unitary Brownian motion  down to the optimal scale (up to polynomial factors).  In order to state our results we introduce the notion of \emph{overwhelming probability.} 

\bed \label{def:op} If $\A_i$ are events indexed by some set $i \in \I$ (and may depend on $N$) then we say that the family of events $\A_i$ hold with overwhelming probability if for all $D>0$ there is a $C>0$ so that,
\beq
\sup_{i \in \I } \pp [ \A_i^c] \leq  N^{-D} ,
\eeq
for all $N \geq C$.
\eed
\bet \label{thm:bulk}
Let $T>0$ and $\eps, \delta, \mfc >0$. For any $0 < t< T$, we define the domain,
\beq \label{eqn:Bt-def-u1}
\B_t := \left\{ z \in \cc : 5 > \log |z| > \frac{ N^{\delta}}{N | \Re[ \tilf (z, t) ] |} \vee N^{-\mfc} \right\}.
\eeq
Then, with overwhelming probability we have uniformly for all $0 < t <T$ and $z \in \B_t$ that,
\beq
| f(z, t) - \tilf (z, t)| \leq \frac{N^{\eps}}{N \log |z|}.
\eeq
\eet

Let us now explain the connection between the Cauchy transform and the spectral measure. 
For a general probability measure $\mu$ on the unit circle, one may recover $\mu$ via the weak limits,
\beq
(2 \pi )\d \mu (\theta)  =  \lim_{ r \uto 1} \left( \Re\left[ \int \frac{ \e^{ \i \theta'}  + r \e^{ \i \theta}}{ \e^{ \i \theta'}  - r \e^{ \i \theta}} \d \theta' \right] \right) \d \theta  = -\lim_{ r \dto 1} \left( \Re\left[ \int \frac{ \e^{ \i \theta'}  + r \e^{ \i \theta}}{ \e^{ \i \theta'}  - r \e^{ \i \theta}} \d \theta' \right] \right) \d \theta .
\eeq
However, this is not useful in order to obtain effective estimates on, e.g., the number of eigenvalues in an interval. The Helffer-Sj{\"o}strand formula \cite{davies1995functional} for measures on $\rr$ allows one to relate empirical averages of test functions to integrals of the Stieltjes transform over $\cc$, turning estimates on the Stieltjes transform into effective estimates on the eigenvalues \cite{semi-lec}. In Section \ref{sec:HS} we quickly develop a version of the Helffer-Sj{\"o}strand formula for measures on the unit circle (like the usual HS formula, it is a consequence of Green's theorem). Similar formulas have appeared before in the literature \cite{mbarek2015helffer}, but the form given here is well-adapted to our purposes.   Using this and the above theorem as input, we obtain  the following.
\bec \label{cor:interval}
Let $T>0$ and $\eps >0$. With overwhelming probability, the following holds uniformly over  all intervals $I \subset (-\pi, \pi]$, and all $0 < t < T$,
\beq
\left| \left| \{ \lambda_i (t) = \e^{\i \theta_i (t) } : \theta_i \in I \} \right| - N \int_I \rho_t (\theta) \d \theta \right| \leq N^{\eps} + t^{-1/2} N^{-1/\eps} ,
\eeq
for $N$ large enough.
\eec
The above corollary shows that as long as $t \geq N^{-C}$ for some $C>0$ then the number of eigenvalues in any interval is given by the limiting spectral measure up to an arbitrarily small polynomial error, with very high probability. This is the optimal scaling up to perhaps replacing the $N^{\eps}$ error by some sort of logarithmic error.

 Theorem \ref{thm:bulk} is proven in Section \ref{sec:bulk}, and Corollary \ref{cor:interval} is derived in Section \ref{sec:interval}. This is the most straightforward of our results as it is the most literal translation of the methods of \cite{huang2019rigidity,adhikari2020dyson} to the unitary setting. Nonetheless, estimates as sharp as those of Corollary \ref{cor:interval} were not known before this work.

 Our method relies on an application of Gronwall's inequality to the the function $ t\to f (z (t), t) - \tilf (z (t), t)$ where $z(t)$ is a characteristic as above. The most delicate estimates are associated with estimating the martingale term $M_t (z(t))$ defined above, as well as the argument around \eqref{eqn:bulk-e1} and \eqref{eqn:bulk-e2}. The integral form of Gronwall's inequality we apply involves an exponential term which this latter argument estimates. Here, we cannot lose any constants or else the error term would be far too large to close our argument.

Recall that for $ t<4$ the support of $\rho_t$ is not yet the entire unit circle and is instead the interval $I_t := [ - \Theta_t, \Theta_t]$ where $\Theta_t$ is as in \eqref{eqn:int-thetat}.  The estimates of Theorem \ref{thm:bulk} are insufficient to address the natural question of whether or not there are eigenvalues outside of $I_t$, or their typical distance from the edges of $I_t$. 

As mentioned above, Collins, Dahlqvist and Kemp \cite{collins2018spectral} showed that  with high probability there are no eigenvalues outside any open set containing $I_t$.  However, such a statement does not yield the correct order of fluctuations of the extremal eigenvalues of $U_t$.

Before stating our results regarding the extremal eigenvalues, let us first ascertain what we expect for the order of magnitude of the fluctuations.  For $\delta <t < 4 -\delta$, we show in Appendix \ref{a:dos} that  for $E>0$ sufficiently small,
\beq
\rho_t (\Theta_t - E) = c_t E^{1/2} (1 + \O (E) ).
\eeq
That is, the spectral measure vanishes like a square-root at the edges of its spectrum. 
The square-root behavior near spectral edges is generic in random matrix theory and is associated with limiting Tracy-Widom fluctuations on the order of $\O (N^{-2/3})$.  While not explicitly formulated, we expect that the estimates of \cite{collins2018spectral} in fact show that the extremal eigenvalues of $U_t$ are no more than $\O (N^{-c})$ from the edges of $I_t$ for some small, explicit $c>0$. However, the interparticle distance associated with the square-root behavior is $\O (N^{-2/3} )$ and so such an estimate nonetheless falls short.

Our next result shows that with overwhelming probability, the extremal eigenvalues are in fact within $\O (N^{-2/3+\eps} )$ of the edges of the support $I_t$. These are analogues of the well-known rigidity results in random matrix theory at the edge (see, e.g., \cite{erdos2012rigidity} for the first such estimates for generalized Wigner matrices). This is also the analogue of the results of \cite{adhikari2020dyson} on the Hermitian DBM in the unitary setting.
\bet \label{thm:edge}
Let $\delta >0$ and  $\eps >0$ be sufficiently small. With overwhelming probability, the following holds uniformly for all $t$ satisfying $\delta < t < 4- \delta$,
\beq
\left| \left\{ i : \lambda_i (t) = \e^{ \i \theta}, \theta \in [-\pi, \pi] \backslash [ - \Theta_t - N^{-2/3+\eps}, \Theta_t + N^{-2/3+\eps} ] \right\} \right| = 0 ,
\eeq 
and for all $0 \leq t \leq \delta$,
\beq
\left| \left\{ i : \lambda_i (t) = \e^{ \i \theta}, \theta \in [-\pi, \pi] \backslash [ - \Theta_t - N^{-\eps/6}, \Theta_t + N^{-\eps/6} ] \right\} \right| = 0 .
\eeq 
\eet
For short times $t \ll 1$, the measure $\rho_t$ looks like an approximate semicircle centered at $\theta =0$ of width $\sqrt{t}$ and $\rho_t(0) \asymp t^{-1/2}$. One therefore expects a different scaling of the interparticle distance for short times $t$. The error we obtain is not optimal for short times $t$, but this regime is not the main focus of our work and so we do not try to optimize our approach here.

 Together with Corollary \ref{cor:interval}, we can then deduce the following rigidity estimates. To introduce them, we require some further notation. Note that for all times $t>0$ the joint law of the eigenvalues of $U_t$ has a density with respect to Haar measure, and so for each fixed time $t_0 > 0 $, the eigenvalues are almost surely distinct. On the other hand, it follows from \cite{cepa2001brownian} that for distinct initial data, the solution to \eqref{eqn:unit-dbm} exists as a strong solution for all times $t> t_0$ and moreover the eigenvalues do not intersect. Taking $t_0 \to 0$ it follows that with probability $1$, the eigenvalues are distinct for all times $t>0$. Since the eigenvalue locations are continuous functions of time and all start at the location $z=1$ at $t=0$ it follows that there is a labelling $\{ \lambda_i (t) \}_{i=1}^N$  so that we can write $\lambda_i = \e^{ \i \theta_i (t)}$ for continuous $\theta_i (t)$ starting at $0$ and for all $t>0$ satisfying,
\beq \label{eqn:int-ordering}
\theta_1 (t) < \theta_2 (t) < \dots < \theta_N (t) < \theta_1 (t) + 2 \pi.
\eeq
 Note that it is possible for $\lambda_i (t) = \e^{ \i \theta_i (t)}$ to wrap many times around the unit circle. E.g., each $\theta_i (t) $ can take any value in $\rr$ such that the above ordering is respected.  However,  since the statement of Theorem \ref{thm:edge} holds on an event of overwhelming probability for all $t$ simultaneously, we can rule out any eigenvalues passing through the point $\theta = \pi$ on this event (this could also be concluded by passing from estimates holding on a sufficiently finely-spaced grid of times to all $t$ using the Hoffman-Wielandt inequality).  This observation allows us to formulate the following rigidity estimates.

Denote by $\gamma_i(t)$ the quantiles of $\rho_t$,
\beq
\frac{i}{N} = \int_{-\pi}^{\gamma_i(t)} \rho_t (\theta ) \d \theta ,
\eeq
with the convention that $\gamma_N$ is either the right spectral edge or $\pi$ for $t \geq 4$. We have the following.
\bec \label{cor:edge-rig}
Let $\delta >0$ and $\eps >0$. The following holds with overwhelming probability uniformly for all $t$ satisfying $\delta < t < 4- \delta$ and all $ 1 \leq i \leq N$. We have,
\beq
| \theta_i (t) - \gamma_i (t) | \leq N^{\eps} \frac{1}{N^{2/3} \min \{ i^{1/3} , (N+1 - i )^{1/3} \} }.
\eeq
\eec
Theorem \ref{thm:edge} is proven in Section \ref{sec:edge}.  Corollary \ref{cor:edge-rig} follows in a straightforward manner from Corollary \ref{cor:interval} and Theorem \ref{thm:edge} and so we omit the proof (see, e.g., Section 3.3 of \cite{huang2019rigidity}).

Compared to the work \cite{adhikari2020dyson}, we encounter several new difficulties in our unitary setting in establishing these edge rigidity results. These mostly have to do with the fact that establishing the above results heavily depends on detailed analysis of the behavior of characteristics near the spectral edge. This behavior depends especially on the distance along the unit circle of characteristic from the location $\pm \Theta_t$ (the angular coordinate) as well as the distance from the unit circle. However, these coordinates introduce curvature whereas in the real line setting these are flat Cartesian coordinates of the real and imaginary part of the characteristic. This is further complicated by the fact that for short times $t$, the spectral measure is very peaked and that for long times $t$, the characteristics will leave the small neighbourhoods of the spectral edges for which we can develop expansions of $\tilf (z, t)$.  

We overcome the short-time difficulties mainly by sacrificing obtaining optimal estimates for short times; i.e., for short times we consider only characteristics that are somewhat far from the spectral edge. The second fact we use to overcome the difficulties associated with curvature and long times is monotonicity of the radial coordinate of the characteristics in time. This second fact is lacking for the general potential processes considered in \cite{huang2019rigidity,adhikari2020dyson}, and is one of the sources of the short-time restrictions in those works (``shocks'' can develop for these processes).  Essentially, monotonicity of the characteristics allows us to split the paths into ``short'' and ``long'' time regimes. In the short-time regimes, we can use a combination of analytic approaches to square-root measures and the semi-explicit formulas for $\tilf (z, t)$ to control the behavior of characteristics close to spectral edges. In the long time regime, things are far away from the spectrum and so relevant quantities can usually be bounded by constants.

We now turn our attention to times $t \sim 4$. This is a critical time for unitary Brownian motion, as at $t=4$, the two spectral edges $\pm \Theta_t$ merge, and the support of the density of states becomes the entire unit circle for later times $t$. In fact, as we show in Appendix \ref{a:dos}, the spectral measure at $t=4$ has a cusp singularity,
\beq
\rho_4 (\pi + E) = c |E|^{1/3} (1 + \O ( |E|^{1/3} ) ),
\eeq
for some constant $c>0$.  Moreover, for times $t$ near $4$, the spectral measure undergoes a transition where the two  edges gradually form a ``near-cusp'', become a cusp, and then form a small local minimum as $t$ ranges from slightly less than $4$ to slightly larger than $4$.

This behavior is identical to that found in the theory of the so-called quadratic vector equation. The quadratic vector equation is a generic equation characterizing the density of states of certain classes of mean field self-adjoint random matrix models. In a series of works \cite{qve,sing,univ-wig-type}, Ajanki, Erd{\H{o}}s and Kr{\"u}ger carried out a systematic study of the solutions to the quadratic vector equation. In particular, they found that the only possible singularities that may occur are cusps, where the density of states of vanishes like a cube root, and square roots, occurring at either external or internal edges. Moreover, they characterized transitional regimes where intervals of the density of states merge or split. In such regimes, they showed that the leading order of the density of states is always given by universal shape functions arising from Cardano's formula for the roots of third-degree polynomials. There are two such functions; the first, $\Psi_e$ corresponds to the case when two separate intervals merge and describe a transition from square-root to cubic behavior. The second, $\Psi_m$ describes what occurs after the cusp forms, when the density of states has a small minimum.

In fact, in Appendix \ref{a:dos}  we show that for times $t <4$ that the density of states $\rho_t$ of unitary Brownian motion is described by $\Psi_e$ and for times $t > 4$ by $\Psi_m$, the universal shape functions of \cite{qve}.  On the one hand, this is remarkable as there is no quadratic equation describing the density of states of unitary Brownian motion and moreover that the eigenvalues are on the unit circle instead of the real line. On the other hand, the universal shape functions arise from expansions of the Cauchy or Stieltjes transform near critical points (i.e., the spectral edges or minima) as soon as one is guaranteed that either the coefficients of the quadratic or cubic terms is non-degenerate and so this behavior is somewhat expected.

Theorems \ref{thm:bulk} and \ref{thm:edge} above do not capture the behavior of the extremal eigenvalues in the case of the near-cusp, when times $t$ are very close to $4$. Using the formula for $\Theta_t$ above, we have for $t < 4$ that the gap between the two spectral edges scales like,
\beq
\Delta_t := 2 ( \pi - \Theta_t ) = \frac{1}{3} (4 -t )^{3/2} (1 + \O (4-t) ) .
\eeq
By our calculations of $\rho_t$ and the asymptotics of the shape function $\Psi_e$ we have that in a vicinity of the edge,
 the density of states behaves like a re-scaled square-root,$^1${\let\thefootnote\relax\footnotetext{1. The notation $\asymp$ means that for two positive possibly $N$-dependent quantities, $a_N \asymp b_N$ implies that there is a constant $C>0$ so that $C^{-1} a_N \leq b_N \leq C a_N$.}}
\beq
\rho_t ( \Theta_t - E) \asymp \frac{ E^{1/2}}{ (4-t)^{1/4} } \asymp \frac{E^{1/2}}{\Delta_t^{1/6}}, \qquad 0 \leq E \leq \Delta_t .
\eeq
From the behavior of $\rho_t$ it follows that the natural fluctuation scale of the extremal eigenvalues is $\Delta_t^{1/9} N^{-2/3}$. This is of the same order of magnitude as $\Delta_t$ when $4-t = N^{-1/2}$. It follows that for $t \ll 4 - N^{-1/2}$ one expects that the extremal eigenvalues are still located near their respective edges. For larger $t$, the fluctuations of the extremal eigenvalues is larger than $\Delta_t$ and so no such rigidity estimate is expected.  The first statement is the content of the following theorem.
\bet \label{thm:cusp}
Let $\delta >0$ and let $\eps >0$. With overwhelming probability the following holds. Uniformly for all $t$ satisfying $2 < t < 4 - N^{-1/2+\delta}$ we have that,
\beq
\left| \left\{ i : \lambda_i (t) = \e^{ \i \theta}, \theta \in [-\pi, \pi] \backslash [ - \Theta_t - N^{-2/3+\eps}\Delta_t^{1/9}, \Theta_t + \Delta_t^{1/9} N^{-2/3+\eps} ] \right\} \right| = 0
\eeq
\eet 
The above theorem is proven in Section \ref{sec:cusp}.  In principle, its proof could be absorbed into the proof of Theorem \ref{thm:edge}. However, handling the three separate scaling regimes, when $t$ is small, of intermediate size, and close to $4$, would require significant additional notation and burden the reader by overly complicating the proofs. We have chosen instead to treat the ``short'' and ``long'' time regimes separately; in fact we will use the result of Theorem \ref{thm:edge} in our proof of Theorem \ref{thm:cusp}, initializing the dynamics at an intermediate time $0 \ll t \ll 4$, conditional on the results of Theorem \ref{thm:edge} holding.

 Local laws and rigidity for random matrices exhibiting cusps were established in \cite{erdHos2020cusp} using resolvent methods. The work \cite{cipolloni2019cusp} also establishes rigidity results for certain interpolating ensembles using a dynamical, PDE-based approach not related to our approach, although both works study the formation of cusps under eigenvalue dynamics.

The main obstacle in proving the above theorem is to understand how the cusp scaling affects the behavior of the characteristics. In particular, we must understand how the angular and radial coordinates of the characteristics are affected by this new scaling. Luckily, in the regime where we expect to prove edge rigidity, there is still a small interval where the density of states behaves like a square root, albeit rescaled by a factor involving $\Delta_t$.   The characteristics relevant to edge rigidity start close to the spectral edge, and some of our calculations of square-root behavior in the earlier short-time regime of Theorem \ref{thm:edge} are applicable here, after finding appropriate re-scalings by $\Delta_t$ of the angular and radial characteristic coordinates.  Nonetheless, we still need to handle the behavior of the characteristics for times of order $1$, and so more analysis of the characteristics is required. This is further complicated by the curvature of our coordinate system as well as the fact that the scaling factor $\Delta_t$ is itself time-dependent and will in general differ by several orders of magnitude over the time intervals we consider.

We can use Theorem \ref{thm:cusp} together with Corollary \ref{cor:interval} to deduce the following rigidity estimates, in a similar manner to Corollary \ref{cor:edge-rig}. We omit the proof.
\bec \label{cor:cusp-rig}
Let $\delta >0$ and $\eps >0$.  For $N^{-1/2+\delta} \leq 4-t \leq 10^{-1}$ we have that the following estimates hold with overwhelming probability. Uniformly for all $i$ satisfying $1 \leq i \leq N (4-t)^2$ we have,
\beq \label{eqn:cusp-rig-u1}
| \theta_i (t)  - \gamma_i (t) | \leq N^{\eps} (4-t)^{1/6} \frac{1}{N^{2/3} i^{1/3} }
\eeq
and for $N (4-t)^2 \leq i \leq N/2$ we have,
\beq \label{eqn:cusp-rig-u2}
| \theta_i (t) - \gamma_i (t) | \leq N^{\eps} \frac{1}{N^{3/4} i^{1/4}}.
\eeq
Analogous estimates hold for indices $i$ near $N$.
\eec

\noindent{\bf Remark}. The reason for the two regimes of indices less than or greater than $N(4-t)^2$ is due to the behavior of the limiting spectral measure near $-1 = \e^{ \i \pi}$ for times $t$ close to $4$. With $s = 4-t$, we have that $\rho_t (\Theta_t -E) \asymp E^{1/2} /s^{1/4}$ for $E \leq s^{3/2}$ and $\rho_t(\Theta_t - E) \asymp E^{1/3}$ for $0.1 \geq E \geq s^{3/2}$ (see Proposition \ref{prop:cusp-shape-u1} and \eqref{eqn:cusp-shape-u1}). The form of the RHS of the estimates \eqref{eqn:cusp-rig-u1} and \eqref{eqn:cusp-rig-u2} reflect the different interparticle distances in each of these regimes. \qed

For larger times $t \gtrsim  4 - N^{-1/2}$, the optimal estimates for the Cauchy transform are in fact included in Theorem \ref{thm:bulk}; compare with, e.g., the local laws of \cite{erdHos2020cusp}.  Note that Theorem \ref{thm:bulk} alone is insufficient to conclude rigidity estimates similar to Corollaries \ref{cor:edge-rig} or \ref{cor:cusp-rig}. In particular, the above results cannot rule out the case that after time $t \sim 4$, that all of the eigenvalues wrap around the unit circle many times.

However, due to the fact that $\theta_N (t)$ and $\theta_1 (t)$ cannot cross, the ``winding number'' (we use this term  loosely) can be determined from the behavior of the center of mass,
\beq
\bar{\theta} (t) := \frac{1}{N} \sum_{i=1}^N \theta_i (t).
\eeq
From either \eqref{eqn:unit-def} or \eqref{eqn:unit-dbm} one can check that \emph{formally} $ \d \bar{ \theta} = N^{-1} \d B$ for a Brownian motion $B$. In Section \ref{sec:com} we justify this using a careful application of the analytic functional calculus.
\bep \label{prop:com-evolution}
For any $t >0$ we have almost surely that,
\beq
\bar{\theta} (t) = \frac{1}{N} \tr (W_t)
\eeq
where $W_t$ is the standard complex Hermitian Brownian motion in \eqref{eqn:unit-def}.  
\eep

We expect that a version of the above statement could also be deduced from Lemma 5 of \cite{meckes2018convergence} but we provide a direct proof, aspects of which may be useful in other settings.

This allows us to deduce the following. Let us denote the extended quantiles $\tilde{\gamma}_i (t)$ of $\rho_i (t)$ as follows. For $1 \leq i \leq N$ we let $\tilgamma_i (t) = \gamma_i (t)$. For $i > N$ we 
\beq
\tilgamma_i (t) = 2 \pi + \gamma_{i-N} (t)
\eeq
and for $i < 1$ we let
\beq
\tilgamma_i (t) = - 2 \pi + \gamma_{i+N} (t).
\eeq
\bec \label{cor:bulk-rig}
Let $\eps >0$ and $\delta >0$. With overwhelming probability we have uniformly for all $t$ satisfying $\delta < t < \delta^{-1}$ that,
\beq
 \tilgamma_{i - N^{\eps} } (t) \leq \theta_i (t) \leq \tilgamma_{i + N^{\eps}} (t).
\eeq
\eec
Note that this is weaker for the edge eigenvalues for times $t < 4 - N^{-1/2+\eps}$ than Corollaries \ref{cor:edge-rig} and \ref{cor:cusp-rig}, and is  only useful in the regime where we can no longer rule out the existence of eigenvalues in the gap between the spectral edges, or there is no longer a gap. Corollary \ref{cor:bulk-rig} is proven in Section \ref{sec:bulk-rig}. 

\subsection{Further discussion and motivation}

In addition to being an intrinsic question about the nature of the process $U_t$, these local law and rigidity estimates have been well studied in the context of Hermitian random matrix theory. There, a primary motivation is the study of the universality of the local eigenvalue statistics: that is, whether or not the limiting local eigenvalue statistics coincide with those of the Gaussian ensembles, which admit exact formulas. 


The short scale behaviors of the repulsive interaction terms of the eigenvalue process of the Hermitian and Unitary Brownian motions \eqref{eqn:dbm} and \eqref{eqn:unit-dbm}   are similar. Based on this and the general belief in the universality of large correlated systems, it is natural to conjecture that the local eigenvalue statistics of \eqref{eqn:unit-dbm} should be given by the same statistics as the GUE in the limit $N \to \infty$.

Indeed, there have been many developments in the universality theory of the local eigenvalue statistics both within the larger context of random matrix theory, as well as that of the Hermitian Dyson Brownian motion started from general initial data. Note that if the initial data of \eqref{eqn:dbm} is not the $0$ matrix, then the joint eigenvalue distribution of $X_t$ is no longer that of a re-scaled GUE. Nonetheless, local scaling limits of DBM with general initial data has been obtained in \cite{convergence,fixed,es,landon2017edge}.  

Universality has been established for wide classes of Hermitian random matrices.  We refer the interested reader to, e.g., the book \cite{erdos2017dynamical} for an overview of these developments as well as to the seminal papers of Tao and Vu \cite{taovu1,tao2010random}, and Erd\H{o}s, Schlein and Yau \cite{local-relax}.

Given these advances in the theory of Hermitian random matrices, it is natural to turn to the question of universality of unitary Brownian motion. An important tool in many of the proofs of Hermitian universality are the aforementioned rigidity and local law estimates. 

 The main contribution of the present work is to establish these results.  It is then a subject of current investigation to use these estimates to prove the local universality of unitary Brownian motion: that the local eigenvalues statistics of $U_t$ are given by the Tracy-Widom, and Pearcey and Sine kernels in the various scaling regimes of interest.

\subsection{Notation} 

The notion of overwhelming probability was defined above in Definition \ref{def:op}.   We let $c >0$ and $C>0$ denote small and large constants respectively. In general, we allow them to increase or decrease from line to line. For two positive $N$-dependent quantities (or quantities depending on some auxiliary parameters, usually time $t$) the notation $a_N \asymp b_N$ means that there is a constant $C>0$ so that $C^{-1} a_N \leq b_N \leq C a_N$.   The notation $a_N \ll b_N$ means $a_N / b_N \to 0$ as $N \to \infty$. We will use this notation sparingly, but somewhat informally. When used we always have an explicit estimate, e.g., $a_N \leq b_N / \log(N)$.  For complex $c_N$ and $d_N$, the notation $c_N = \O (d_N)$ means $|c_N| \leq C | d_N|$ for some $C>0$.

\subsection{Organization of paper}

Sections \ref{sec:bulk}, \ref{sec:edge} and \ref{sec:cusp} are meant to be read in a relatively linear fashion. These sections prove Theorems \ref{thm:bulk}, \ref{thm:edge} and \ref{thm:edge}, respectively. The analysis in each section directly builds off that of the previous section. Section \ref{sec:bulk} treats the ``bulk'' local law (i.e., a local law with an error that is optimal only in the bulk) and the treatment of the characteristics is relatively straightforward. Section \ref{sec:edge} treats the cases $t < 4 - \delta$ where $\rho_t$ has a regular square-root edge. Here, the treatment of characteristics (and the resulting estimates of the quantities such as the Martingale term in evolution equation of $f(z, t)$) is more complicated. Finally, Section \ref{sec:cusp} deals with the formation of the cusp. 

In the short Section \ref{sec:com} we prove Proposition \ref{prop:com-evolution}, that the centre of mass, or averaged winding number, is described by a Brownian motion. In Section \ref{sec:rigidity} we establish the analog of the Helffer-Sj{\"o}strand formula and use it to deduce Corollary \ref{cor:interval}. The latter is very similar to what has appeared in \cite{huang2019rigidity} and so not all details are provided. 

In Appendix \ref{a:dos} we establish various properties of the limiting spectral measure $\rho_t$. We use as input its characterization in terms of conformal maps of Biane \cite{biane1997segal}, as well as arguments of Ajanki, Erd{\H{o}}s and Kr{\"u}ger \cite{qve} involving solutions of Cauchy/Stieltjes transforms of approximate cubic equations, Cardano's formula and the universal shape functions $\Psi_m$ and $\Psi_e$. Finally, Appendix \ref{a:calc} collects various calculus-type inequalities used in the proof.

A shorter version of this paper has been submitted to PTRF. This is the longer version that contains complete proofs of all results in the paper, appearing as v3 of arXiv:2202.06714 on the arXiv.

\section{Bulk estimates} \label{sec:bulk}

In this section we prove Theorem \ref{thm:bulk}. Fix a final time $T>0$. This time may depend on $N$, but stays bounded above. We introduce the characteristic maps via
\beq
\C_t (z) = z \exp\left[  - \frac{ (T-t) \tilf (z, T) }{2} \right].
\eeq
That is, the function $t \to \C_t (z)$ satisfies the characteristic equation \eqref{eqn:dzdt} and has final condition $\C_T (z) = z$.

From the above, it is clear that the real parts of $f (z, t)$ and $\tilf (z, t)$ will play  important roles. We record here the identity,
\beq
\Re[ f(z, t) ] = \frac{1}{N} \sum_{i=1}^N \frac{1- |z|^2}{ | \lambda_i - z|^2} \label{eqn:realpart} .
\eeq
We will also have use for,
\beq \label{eqn:delzf}
\del_z f (z, t) = \frac{2}{N} \sum_{i=1}^N \frac{ \lambda_i (t) }{ ( \lambda_i (t) - z)^2} ,
\eeq
and the inequality
\beq \label{eqn:delzf-bd}
| \del_z f(z, t) | \leq \frac{2}{ 1- |z|^2} \Re[ f (z, t) ]  .
\eeq
In the remainder of the section we will also make use of the spectral domains $\B_t$ that were defined above in \eqref{eqn:Bt-def-u1}. 
We collect some elementary properties of the characteristics. 
\bel \label{lem:bulk-char}
For any $|z| >1$, the map,
\beq
t \to \log | \C_t (z) |
\eeq
is decreasing in time. Secondly, there is a constant $C>0$ so that for any $z \in \B_T$ we have,
\beq
| \C_t(z) | \leq C\e^{CT}.
\eeq
\eel
\proof The first claim follows from the fact that $\Re[ \tilf (z, t) ] <0$ for $|z| >1$.  For the second, let $z_t = \C_t (z)$. Note that $\tilf (z_t, t) = \tilf (z, T)$ for all $t$. If at any time $t$ we have $|z_t| > 2$, then $| \tilf (z_t, t) | \leq 4$ and so $| \C_t (z) | \leq |z|\e^{2 T}$. So either $|z_t| < 2$ for all $t$ or $|z_t| \leq |z|\e^{2 T}$ for all $t$. This yields the claim. \qed

We fix a final point $z \in \B_T$. We will first prove that Theorem \ref{thm:bulk} holds at a single $z$ by tracking the evolution of $f(z, t)$ along the characteristic,
\beq
z_t := \C_t (z).
\eeq
The extension to all $z$ and all $0 < t < T$ will be detailed later.  We now introduce the stopping time $\tau$ via,
\beq
\tau := \inf \left\{ s \in [0, T] : | f (z_s, s) - \tilf (z_s, s) | > \frac{ N^{\eps}}{N \log |z_s | } \right\} \wedge T
\eeq
where we choose $\eps < \delta/10$, with $\delta$ as in the definition of $\B_T$.  
We first note that,
\begin{align}
\d \left( f (z_t, t) - \tilf (z_t, t) \right) &= - \frac{ z_t \del_z f(z_t, t) }{2} \left( f (z_t, t) - \tilf (z_t, t) \right) + \d M_t (z_t) ,
\end{align}
with the Martingale term defined above in \eqref{eqn:M-def}. 
Hence,
\beq
f ( z_\tau , \tau ) - \tilf (z_\tau, \tau ) = \E_1 ( \tau) + \E_2 ( \tau)
\eeq
where
\beq
\E_1(t) := - \int_0^{t} \frac{ z_s \del_z  f (z_s, s) }{2} ( f (z_s , s) - \tilf (z_s, s) ) \d s
\eeq
and
\beq
\E_2 (t) := \int_0^{t} \d M_s (z_s ).
\eeq
We first prove the following estimate on the martingale term.
\bep \label{prop:bulk-bdg}
For all $\eps_1 >0$ we have,
\beq
\pp \left[ \exists t \in [0, \tau] : | \E_2 (t) | > \frac{ N^{\eps_1}}{N \log |z_t | } \right] \leq C \log(N) \e^{ - c N^{\eps_1}}.
\eeq
\eep
\proof We fix a sequence of intermediate times $t_k$ with $k=1, \dots , M$ in $[0, T]$ such that $\log |z_{t_k} | \leq 2 \log |z_{t_{k+1} } |$ and $t_M = T$.  Then $M \leq C \log(N)$ for some $C>0$ since $\log |z_0|$ is bounded by Lemma \ref{lem:bulk-char}. Let $\tau_k = \tau \wedge t_k$. The quadratic variation of $ \E_2 ( \tau_k )$ satisfies,
\begin{align}
\langle \E_2 ( \tau_k), \bar{\E}_2 ( \tau_k) \rangle &= \frac{4}{N^2} \int_0^{\tau_k} |z_s|^2 \frac{1}{N} \sum_{i=1}^N \frac{1}{ | \lambda_i (s) - z_s |^4} \d s \notag\\
&\leq \frac{4}{N^2} \int_0^{ \tau_k } \frac{ |z_s|^2}{ (|z_s| -1)^2 } \frac{1}{N} \sum_{i=1}^N \frac{1}{ | \lambda_i (s) - z_s |^2} \d s \notag\\
&\leq \frac{C}{N^2} \int_0^{ \tau_k } \frac{ | \Re[ f (z_s, s) ] |}{ ( \log |z_s | )^3} \d s.
\end{align}
In the second line we used the trivial estimates $| \lambda_i (s) - z_s | \geq |z_s| -1$. In the last line we used Lemma \ref{lem:bulk-char} to bound $|z_s|$ as well as the representation \eqref{eqn:realpart}. We also used that $\log(r) \leq r-1$ for $r>1$.  For $s < \tau$ we have 
\beq
| \Re[ f(z_s, s) ] - \Re[ \tilf (z_s, s) ] | \leq \frac{N^{\eps}}{N \log |z_s|  } \leq \frac{ N^{\eps}}{N \log |z_T|}.
\eeq
By definition of $\B_T$ we see that
\beq
\frac{ N^{\eps}}{N \log |z_T|} \leq N^{\eps-\delta} | \Re[ \tilf (z_T, T) ] | = N^{\eps-\delta} | \Re[ \tilf (z_s, s) ] |
\eeq
where we used that $\tilf$ is constant along characteristics. Since $\eps < \delta$ we therefore see that,
\beq \label{eqn:ref-char-bd}
| \Re[ f (z_s, s) ] | \leq (1 + N^{\eps-\delta} ) | \Re[ \tilf (z_s, s) ] | \leq \left( 1  + \frac{1}{ \log(N) } \right) | \Re[ \tilf (z_s, s) ] |
\eeq
for $s < \tau$. Therefore,
\begin{align}
\int_0^{\tau_k} \frac{ | \Re[ f (z_s, s) ] |}{ ( \log |z_s | )^3} \d s &\leq  2 \int_0^{\tau_k} \frac{ | \Re[ \tilf (z_s, s) ] |}{ ( \log |z_s | )^3} \d s \notag\\
&\leq 2 \int_0^{t_k} \frac{ | \Re[ \tilf (z_s, s) ] |}{ ( \log |z_s | )^3} \d s \notag\\
&\leq \frac{2}{ (\log |z_{t_k} | )^2}.
\end{align}
In the last line we used that $\del_u \log |z_u | = - \frac{1}{2}| \Re[\tilf (z_u, u) ] |$.  By the Burkholder-Davis-Gundy (BDG) inequality (see, e.g., \cite{revuz2013continuous}) and a union bound, we conclude that,
\beq \label{eqn:bulk-a1}
\pp \left[ \exists k : \sup_{0 \leq t \leq \tau_k } | \E_2 ( t) | > \frac{ N^{\eps_1}}{ 2 N \log |z_{t_k } | } \right] \leq C \log(N) \e^{ - c N^{\eps_1} }.
\eeq
Let $0 < t < \tau$ and let $k$ be such that $\tau_k \leq t < \tau_{k+1}$. On the complement of the event on the LHS of \eqref{eqn:bulk-a1} we have,
\beq
| \E_2 (t) | \leq \frac{N^{\eps_1}}{2  N \log |z_{t_k } | } \leq \frac{ N^{\eps_1}}{ N \log |z_t | }
\eeq
by the choice of the $t_k$'s. This completes the proof. \qed

\noindent{\bf Proof of Theorem \ref{thm:bulk}}. Let $z \in \B_T$ with characteristic $z_t$ as above. Let $\tau$ be the stopping time as defined above.  Let $\Upsilon$ be the event of Proposition \ref{prop:bulk-bdg} and define,
\beq
g(s) := \frac{ |z_s \del_z f (z_s, s) |}{2}. 
\eeq
By Gronwall's inequality we have for all $0<t < \tau$ on the event $\Upsilon$,
\begin{align} \label{eqn:bulk-a2}
| f (z_t, t) - \tilf (z_t, t) | \leq \int_0^t g (s) \exp \left[ \int_s^t g(u) \d u \right]  \frac{ N^{\eps_1}}{N \log |z_s| } \d s + \frac{N^{\eps_1}}{N \log |z_t| }.
\end{align}
We now estimate the integral of the function $g$ that appears above in the argument of the exponential function. Using first \eqref{eqn:delzf-bd} we have,
\begin{align} \label{eqn:bulk-e1}
\int_s^t g(u) \d u &= \int_s^t \frac{ |z_u | | \del_z f(z_u, u) |}{2} \d u \notag\\
&\leq \int_s^t \frac{ |z_u| | \Re[ f(z_u, u) ] | }{ |z_u|^2 -1 } \d u \notag\\
&\leq \int_s^t \frac{ | \Re[ f (z_u, u) ] |}{ 2 \log |z_u | } \d u.
\end{align}
In the last inequality we used the elementary inequality,
\beq \label{eqn:bulk-e2}
\frac{x}{x^2-1} \leq \frac{1}{ 2 \log (x) }
\eeq
which holds for $x >1$ (see Lemma \ref{lem:calc}). Using now \eqref{eqn:ref-char-bd} we have,
\begin{align}
\int_s^t \frac{ | \Re[ f (z_u, u) ] |}{ 2 \log |z_u | } \d u \leq& \left( 1+ \frac{1}{ \log(N) } \right) \int_s^t \frac{ | \Re[ \tilf (z_u, u) ]|}{2  \log |z_u | } \d u  \notag\\
= & \left( 1 + \frac{1}{ \log(N) } \right) \log \left( \frac{ \log |z_s|}{ \log |z_t| } \right)
\end{align}
where we used in the last line that $\del_u \log |z_u | = - \frac{1}{2}| \Re[ \tilf (z_u, u) ] |$. Therefore,
\beq \label{eqn:bulk-a3}
 \exp \left[ \int_s^t g(u) \d u \right] \leq C \frac{ \log |z_s|}{ \log |z_t| }  .
\eeq
We used the fact that the assumption $ C\geq  \log |z_u| \geq N^{-\mfc}$ from the definition of $\B_T$ implies that
\beq
\left(  \frac{ \log |z_s|}{ \log |z_t| }  \right)^{1/\log(N) } \leq C.
\eeq
Substituting \eqref{eqn:bulk-a3} into \eqref{eqn:bulk-a2} yields,
\begin{align}
|f (z_t, t) - \tilf (z_t, t) | \leq \frac{ C N^{\eps_1}}{N \log |z_t|} \int_0^t g(s) \d s + \frac{ N^{\eps_1}}{N \log |z_t| }.
\end{align}
From the definition of $g$ and \eqref{eqn:delzf-bd} and \eqref{eqn:ref-char-bd} we see that,
\beq
g(s) \leq C \frac{ | \Re[ \tilf (z_s, s) ] |}{ \log |z_s |}
\eeq
for $s < \tau$. Therefore,
\beq
\int_0^t g(s) \d s \leq C \log(N)
\eeq
using that $\log |z_t| \geq N^{-\mfc}$. We conclude that for all $ 0 < t < \tau$ that,
\beq
|f (z_t, t) - \tilf (z_t, t) | \leq  C \frac{ \log(N) N^{\eps_1}}{N \log |z_t|}.
\eeq
Taking $\eps_1 < \eps/2$, where $\eps >0$ is in the definition of the stopping time, we see that on the event $\Upsilon$, we must have $\tau = T$ for all $N$ large enough. This proves the estimate of Theorem \ref{thm:bulk} at the final time $T$ and at the point $z$. The extension to all points $z \in \B_T$ may be done by first taking a union bound over $\O (N^C)$ points and then using
\beq
| f_\mu (z) - f_\mu (w) | \leq C|z-w| \left( \frac{1}{ \log |z| } + \frac{1}{ \log |w|} \right)^2
\eeq
valid for any Cauchy transform and $|z|, |w|$ in bounded regions of $\cc$. The above estimate is elementary and can be proven directly from the definition
$$
f_\mu (z) := \int_0^{2 \pi} \frac{ \e^{ \i \theta} + z}{ \e^{ \i \theta }- z  }\d \mu ( \theta)
$$
and the fact that for $|z| >1$ we have for the quantity in the denominator $|\e^{ \i \theta} -z | \geq |z| -1 \geq c \log|z|$.

The extension to all times $t < T$ is done in a similar manner; one first proves the estimate for all $z \in \B_t$ for all times $t$ in a well-spaced grid of $[0, T]$ of at most size $\O (N^A)$, some $A>0$ to be determined. Then one notes that the proof also gives that the estimate holds along the entire characteristic, and along each characteristic, we have, e.g.,  $|z_t - z_s | \leq C N^{C} |t-s|$, for some $C>0$ depending on $\mfc >0$. We need only take $A$ large enough depending on $C >0$.  This completes the proof of Theorem \ref{thm:bulk}. \qed

\section{Edge estimates} \label{sec:edge}

Fix $T < 4$. We will need to consider characteristics ending at many times and so we introduce the characteristic map,
\beq
\C_{s, t} (z) = z \exp \left[ - \frac{ (t-s) \tilf (z, t) }{2} \right]
\eeq
for $0  \leq s \leq t$.  We will need to establish several properties about the behavior of characteristics near the edge $\Theta_t$.  The proof of the following lemma, an exercise in calculus, is postponed to Appendix \ref{a:calc}. 

\bel \label{lem:char-path}
Let $\rho$ be a measure on $[-\pi, \pi]$ such that $\rho ( \theta ) = \rho (-\theta)$ that is supported in $[-E, E]$ for $0 < E < \pi$. Assume either that $\rho ( \theta) \leq M$ or $E < \pi/8$. For any $\eps >0$ there is a $c_\eps >0$, depending only on $E, M$ and $\eps >0$ so that for $0 < r-1 < c_\eps$ and $E <  \theta < \pi - \eps$ we have,
\beq
 \Im [f ( r \e^{ \i \theta} ) ] < \Im [ f ( \e^{ \i \theta } ) ] < \Im [ f ( \e^{ \i E } ) ]
\eeq
where $f(z)$ is the Cauchy transform of $\rho$. We also have for $E < \theta \leq \varphi \leq \pi$,
\beq
0 \leq \Im [ f ( \e^{ \i \varphi } ) ] \leq \Im [ f ( \e^{ \i \theta } ) ] ,
\eeq
with equality occuring only in trivial cases.
\eel

We now prove the following.
\bel \label{lem:char-der}
Let $T<4$ and $ z= r \e^{ \i (\Theta_t + \kappa)}$ for $\kappa >0$. There is a small $c>0$ and $d>0$ so that,
\beq \label{eqn:imf-est-edge}
\Im[  \tilf (z, t) ] - \Im [ \tilf (\e^{ \i \Theta_t } , t ) ] \leq - d \sqrt{ \kappa}
\eeq
for $0 < t < T$ and all $0 < r-1 < c$ and $\kappa < \pi - \Theta_t$.
\eel
\proof By Lemma \ref{lem:calc-3}, the desired estimate hold for all $\kappa < \eps $ and $0 < r-1 < \eps$ for some $\eps >0$. In particular, 
\beq
\Im[  \tilf (\e^{ \i ( \Theta_t + \eps/2)}, t) ] - \Im [ \tilf (\e^{ \i \Theta_t } , t ) ]  \leq - c_1
\eeq
for some $c_1>0$.  Since for $t < T$ the measures are supported away from $\pi$ we see that there is a $\delta >0$ so that 
\beq
| \Im[ \tilf ( z, t) ] - \Im[ \tilf (-1, t ) ] | \leq \frac{c_1}{2},
\eeq
for $|z+1| < \delta$. It follows  that for all $|z+1| < \delta$,
\beq
\Im[ \tilf (z, t) ] - \Im[ \tilf (\e^{ \i \Theta_t} ) ] \leq \frac{c_1}{2} + \Im[  \tilf (\e^{ \i ( \Theta_t + \eps/2}), t) ] - \Im[ \tilf (\e^{ \i \Theta_t} ) ]  \leq - \frac{c_1}{2}
\eeq
where we used the second estimate of Lemma \ref{lem:char-path}. This proves the desired estimate for $|z+1| < \delta$ after possibly decreasing the value of $d>0$. 

On the other hand, by Lemma \ref{lem:char-path} we conclude that there is a $c_2 >0$ so that for $0 < r-1 < c_2$ and $\eps/2 < \kappa < \pi - \Theta_t - \delta/2$ that,
\beq
\Im[ \tilf (r \e^{ \i ( \Theta_t + \kappa)}, t) ] - \Im[ \tilf (\e^{ \i \Theta_t} ) ] \leq \Im[ \tilf ( \e^{ \i ( \Theta_t + \eps/2)}, t) ] - \Im[ \tilf (\e^{ \i \Theta_t} ) ] \leq -c_1.
\eeq
This concludes the proof. \qed

The following contains the properties of the characteristics that we will need.

\bep \label{prop:edge-char}
Let $T <4$.  There is a $\mfa>0$  and $\mfb >0$ depending on $T$ so that the following holds.  Let $z_s = \C_{s, t} ( z)$ for any $t< T$, where $z = r \e^{ \i \theta}$ satisfies $0 < r-1 < \mfa$ and $ \Theta_t < \theta < \pi$. Denote $z_s = r_s \e^{ \i (\Theta_s + \kappa_s )}$. Let $s_*$ be,
\beq
s_* = \inf \{ s < t : (r_s-1) < \mfa \}.
\eeq
Then for $s_* < s < t$ we have,
\beq \label{eqn:edge-c3}
\sqrt{ \kappa_s} \geq \sqrt{ \kappa_t} + \mfb (t-s)
\eeq
and for $s \leq s_*$ we have $r_s \geq 1 + \mfa $.  Furthermore, let  $D(s)$ be a function that obeys,
\beq
\sqrt{ D(s)} \leq \sqrt{D(t)} + \frac{\mfb }{2} (t-s), \qquad s < t.
\eeq
If $\kappa_t \geq D (t)$ then for $s_* < s < t$ we have,
\beq \label{eqn:edge-c1}
\kappa_s \geq D(s) +\frac{\mfb}{2} \sqrt{ \kappa_t} (t-s).
\eeq
Finally, characteristics do not cross the real axis in the complex plane.
\eep
\proof We take $\mfa$ so small so that the conclusion of Lemma \ref{lem:char-der} for $z = r \e^{ \i \theta}$ with $r-1 < 10 \mfa$. For $s < s_*$ we have $r_s > \mfa$ because the radial coordinate is decreasing in time. Let $s_1$ be,
\beq
s_1 = \inf \{ s < t : \kappa_s > 0 \}.
\eeq
Note that $s_1 < t$ as we assume $\kappa_t >0$. We claim that $s_1 \leq s_*$. For $t > s> s_1 \vee s_*$ we have, the following calculation,
\begin{align} \label{eqn:edge-a1}
\del_s \kappa_s & =   \frac{1}{2}\Im[\tilde{f}(z_s,s)] - \frac{1}{2}\Im[\tilde{f}(\e^{ \i \Theta_s},s)] \notag\\
&= \frac{1}{2}\Im[\tilde{f}(z_t,t)] - \frac{1}{2}\Im[\tilde{f}(\e^{ \i \Theta_t},t)] - \frac{1}{2}(\sqrt{\frac{4 -s}{s}} - \sqrt{\frac{4-t}{t}}) \notag \\
   & \leq \frac{1}{2}[\Im[\tilde{f}(z_t,t)] - \frac{1}{2}\Im[\tilde{f}(\e^{ \i \Theta_t},t)]] - c(t-s) \notag\\
   &  \leq -d\sqrt{\kappa_t} -c(t-s) 
\end{align}
The second line is from the fact that $\tilf (\e^{ \i \Theta_t}, t) = \i \sqrt{ 4t^{-1} -1 }$ for all $t$ and that $\tilf$ is constant along characteristics. The third line is straightforward. The last line follows from Lemma \ref{lem:char-der}. In particular, we see that $\kappa_s$ is a decreasing function for $s > s_* \vee s_1$. Since $\kappa_t >0$ it follows that $s_1 \leq s_*$. Therefore,  the final inequality in \eqref{eqn:edge-a1} holds for all $s_* < s < t$. Therefore via integration we obtain,
\beq
\kappa_s \geq \kappa_t + 2c_2 (t-s) \sqrt{ \kappa_t} + c_2 (t-s)^2
\eeq
for some small $c_2 >0$ for all $s_* < s < t$. This is equivalent to the desired inequality.  For the last inequality of the proposition, we have that
    \begin{equation}
        \sqrt{\kappa_s} \geq \sqrt{\kappa_t} + \mfb (t-s) \geq \sqrt{D(t)} + \mfb (t-s) \geq \sqrt{D(s)} + \frac{\mfb}{2}(t-s).
    \end{equation}
We can square this to get
$\kappa_s \geq D(s) + \mfb (t-s) \sqrt{D(s)}$ as well as $\kappa_s \geq \kappa_t + 2 \mfb (t-s) \sqrt{\kappa_t}$. If $D(s) \geq \kappa_t$, we can use the first inequality. Otherwise, we can use the second. 

The final claim that characteristics not cross the real axis is due to the fact that the imaginary part of $\tilf (z, t)$ vanishes for purely real $z$, due to the symmetry of $\rho_t$. \qed 

The real part of the Cauchy transform $f(z, t)$ can be used to detect the presence of outlying eigenvalues. In order to use it, we first need the following estimate on how $\tilf(z, t)$ behaves. The proof of the following appears in Appendix \ref{a:edge-sqr-proof}. 
\bel \label{lem:edge-sqr}
Let $\delta >0$. Let $z = (1+\eta) \e^{ \i (\Theta_t + \kappa) }$. Uniformly in the region,
\beq
0 < \eta < 5, \qquad 0 < \kappa < \pi - \Theta_t 
\eeq
we have for $\delta < t < 4 - \delta$ that,
\beq \label{eqn:sqrroot}
c \frac{ \eta}{ \sqrt{ \kappa+\eta} } \leq | \Re[ \tilf (z, t) ] |  \leq C \frac{ \eta}{ \sqrt{ \kappa+\eta} }
\eeq
for some $c, C>0$.  For $ 0 < t < \frac{1}{2}$ we also have,$^1${\let\thefootnote\relax\footnotetext{$1$. An anonymous referee pointed out that this could be improved to $| \Re[ \tilf (z, t) ] | \asymp t^{-3/4} \frac{ \eta}{ \sqrt{ \kappa + \eta}}$ if $|z-\Theta_t| \leq \sqrt{t}$ and $\asymp \eta  |z|^{-2}$ if $|z- \Theta_t| > \sqrt{t}$ (at least for $z$ in the upper-half plane and $\kappa >0$) using the short-time behavior of the spectral measure proved in Proposition \ref{prop:short-time-shape-u1}.}}
\beq \label{eqn:trivial-edge}
\eta c \leq | \Re[ \tilf (z,t) ] | \leq C \frac{ \eta}{ \kappa^2} .
\eeq
\eel

As advertised, the following lemma allows us to use estimates for the Cauchy transform outside the spectrum to rule out the presence of outlying eigenvalues.

\bel \label{lem:edge-no-eig}
Let $\delta >0$ and assume $\delta < t < 4 - \delta$. 
Let $\eps >0$ and $0 \leq k \leq \frac{2}{3}$.  Suppose that the estimate,
\beq \label{eqn:edgerigidity}
\left|  f(z, t) - \tilf (z, t) \right| \leq \frac{N^{\eps}}{N \sqrt{ \eta + \kappa } \sqrt{ \eta} }
\eeq
holds for all $ z = r \e^{ \i \theta}$ for any $r$ and $\theta$ satisfying,
\beq
\Theta_t + N^{-2/3+k + 5 \eps} \leq | \theta | \leq \pi, \qquad N^{-2/3+k/ 4+ \eps } \leq r -1 \leq c
\eeq
 where  $c >0$ is any positive constant. Then there are no eigenvalues $\lambda = \e^{ \i \theta}$ for $\Theta_t + N^{-2/3+k + 5 \eps } \leq | \theta | \leq \pi$.  The same conclusion holds for $t < \delta$ if we take $k = 2/3 - 5 \eps - \eps/6$. 
\eel
\proof Let $\lambda$ be an eigenvalue in $f(z,t)$. We will show that if $z$ is of the form $r \e^{\i \theta}$ with $r>1$ and $|\theta - \arg \lambda| < r-1$, then $\Re[f(z,t)] \geq \frac{1}{N (r-1)}$. By rotational invariance along the unit circle, it suffices to consider $\lambda= 1$.

We see by direct calculation that
    \begin{equation}
        \Re \left[ \frac{z-1}{z+1} \right] = \frac{r^2-1}{(r-1)^2 + 2 r(1- \cos \theta)} \geq  \frac{r^2-1}{(r-1)^2 + r \theta^2} \geq \frac{r^2-1}{(r-1)^2(1+r)} \geq \frac{1}{r-1}.
    \end{equation}
    To get the second inequality, we used $1- \cos \theta \leq \frac{\theta^2}{2}$. The third inequality comes from $\theta < r-1$.
    
    Thus, we see that if $z$ of the form $r \e^{\i \theta}$ with $r>1$ and $|\theta - \arg \lambda| < r-1$, we know that $\Re[f(z,t)] \geq \frac{1}{N(r-1)}$ provided that $\lambda$ is an eigenvalue appearing in the empirical measure for $f (z, t)$.  
    
    Now, consider the point $z=(1+ N^{-2/3 + k/4+ \eps} ) \e^{\i \theta}$, where $\theta$ is an angle in the region $ |\theta| \in [\Theta_t + N^{-2/3 + k+ 5 \eps},  \pi ]$.  
    
    Applying the triangle inequality and the estimates \eqref{eqn:edgerigidity} and \eqref{eqn:sqrroot}, we see that for large enough $N$,
    \begin{equation}
        |\Re[f(z,t)] |\leq | f(z,t) - \tilde{f}(z,t)| + |\Re[\tilde{f}(z,t)] | < C \frac{\eta}{\sqrt{\kappa + \eta}} + \frac{N^\eps}{N \sqrt{\eta} \sqrt{\kappa + \eta}} < \frac{1}{N \eta}.
    \end{equation}
   Indeed, for the final inequality, we use the hypotheses on $\eta$ and $\theta$ to estimate,
   \begin{equation}
   \begin{aligned}
       &\frac{\eta}{\sqrt{\kappa + \eta}} \leq N^{-1/3- k/4 -3/2 \eps },\\
       &  \frac{N^{\eps}}{N \sqrt{\eta} \sqrt{\kappa+ \eta}} \leq N^{-1/3 -5k/8- 2 \eps},\\
     &  \frac{1}{N \eta} = N^{-1/3 -k/4 -  \eps}.
    \end{aligned}
      \end{equation}
    This shows $|\Re[f(z,t)] | < \frac{1}{N \eta}$ and therefore $\e^{\i \theta}$ cannot be an eigenvalue of the empirical measure associated to $f(z, t)$. 
    
    We now consider the case $0 < t < \delta$ and $k = 2/3 - 5\eps - \eps/6$.  We can instead apply the trivial bound $| \Re[\tilde{f}(z,t)] | \leq C \frac{\eta}{\kappa^2}$ \eqref{eqn:trivial-edge}. We see that at the same choice of $\eta$ and $\kappa$, we have 
 \begin{equation}
   \begin{aligned}
       &\frac{\eta}{\kappa^2} \leq N^{-1/2 + \eps/24 },\\
       &  \frac{N^{\varepsilon}}{N \sqrt{\eta} \sqrt{\kappa+ \eta}} \leq N^{-3/4 +59/48  \eps},\\
     &  \frac{1}{N \eta} = N^{-1/2 +7/24  \eps}.
    \end{aligned}
      \end{equation}
Therefore, for $z = r \e^{ \i \theta}$ as in the statement of the lemma we see that $| \Re[f (z, t)  ] | < \frac{1}{N \eta}$. This completes the proof. \qed

With the above results in hand, we can begin our proof of Theorem \ref{thm:edge}. Before doing so, we need to introduce further notation. Let $T<4$.  We consider times $t < T$. Let $\mfa,\mfb>0$ be the constants of Proposition \ref{prop:edge-char}. Let $\delta >0$ and $\eps >0$. Assume $\eps < 10^{-6}$. These parameters will be fixed until the end of the proof of Theorem \ref{thm:edge}.

 Introduce 
 $D(t)$ to be the function
    \begin{equation}
    D(t) = N^{-\eps/6}, 0 \leq t \leq \delta, \qquad
    D(t) = \max \left\{ (N^{-\eps/12} - \frac{\mfb}{10}(t-\delta))^2, N^{-2/3 + 5 \eps} \right\} .
    \end{equation}
    The function $D(t)$ satisfies the hypotheses of Proposition \ref{prop:edge-char}. In fact, with this choice of $D(t)$ we have for any choices of $ s< t$ that
    \beq \label{eqn:edge-c2}
\sqrt{ D(s) } \leq \sqrt{ D(t) } + \frac{\mfb}{10} (t-s) .
\eeq
    Additionally, define $k(t)$ to be the solution of $N^{-2/3 + k(t) + 5 \eps}= D(t)$.  Define the spectral domains $\Ed_t$ by,
    \beq
    \Ed_t := \{ z = r \e^{ \i \theta} : \Theta_t + N^{-2/3+k(t)+5 \eps } \leq | \theta | \leq \pi , N^{-2/3+k(t)/4+ \eps} \leq r -1 \leq \mfa/2 \}.
    \eeq
    For any characteristic $z_t = r \e^{ \i \theta}$ with $ r>1$ define $\kappa (z_t) = |\theta| - \Theta_t$ and $\eta(z_t) = r-1$. Consider the control parameter,
    \beq \label{eqn:Bz-def}
    B(z_t) := \frac{1}{N \sqrt{ \kappa (z_t) +\eta (z_t) }\sqrt{\eta (z_t)} } \1_{ \kappa (z_t) > 0 } + \frac{1}{N \eta(z_t) } \1_{\kappa (z_t) < 0}.
    \eeq
    \bel \label{lem:edge-ref}
    Let $0 < t < T$ and let $z_s  = C_{s, t} (z)$ where $z \in \Ed_t$. Then, for all $0 < s < t$ we have,
    \beq
    N^{\eps} B(z_s) \leq \frac{1}{ \log(N) } | \Re[ \tilf (z_s, s) ] |.
    \eeq 
    \eel
    \proof Let $s_*$ be as in Proposition \ref{prop:edge-char}. For $s \leq s_*$ we have that
    \beq
    N^{\eps} B (z_s) \leq C N^{\eps-1}.
    \eeq
    On the other hand, by Lemma \ref{lem:edge-sqr}, $| \Re[ \tilf (z_s, s) ] | \geq c$ for some $c>0$ for such $z_s$. We now consider $s_* < s < t$. By Proposition \ref{prop:edge-char} it follows that,
   \beq
   \kappa (z_s ) \geq \kappa (z_t) \geq N^{-2/3+k(t)+5 \eps}
   \eeq
   and $\eta (z_s) \geq N^{-2/3+k(t)/4 + \eps}$.  First consider $t\leq \delta$. Then, $k(t) = 2/3-5\eps - \eps/6$ and so
   \beq
   N^{\eps} B( z_s) \leq N^{-3/4+10 \eps}.
   \eeq
   Moreover, $\Re[ \tilf (z_s, s) ] \geq c \eta (z_s)  \geq N^{-1/2-5 \eps}$.  This proves the estimate for such $t$. Consider now $t > \delta$. In this case, the desired inequality can be rewritten as
   \beq
 \log(N) N^{\eps} B(z_s) \leq |\Re[\tilf (z_s, s) ] | = | \Re[ \tilf (z_t, t) ] |   ,
   \eeq
   using the fact that $\tilf$ is constant along characteristics. 
   Since $B (z_s) \leq B (z_t)$ for $s > s_*$, this 
    reduces to whether,
   \beq
   \eta(z_t)^{3/2} \geq \log(N) N^{\eps-1}.
   \eeq
    But this holds since $\eta(z_s) \geq N^{-2/3+ \eps}$. \qed
    
	We introduce a grid of times $t_i = i TN^{-10}$ for $i = 1, \dots, N^{10}$. For each such $t_i$ we let $\{ u^i_j \}_{j=1}^{N^{10}}$ be a well-spaced grid of $\Ed_{t_i}$. We then define the characteristics,
	\beq
	z^i_j (s) := \C_{s, t_i} (u^i_j ).
	\eeq
	We define each characteristic $z^i_j(s)$  only for $0 < s < t_i$. We introduce the stopping times,
	\beq
	\tau_{ij} := \inf \{ s \in (0, t_i ] : | f( z^i_j (s),s) - \tilf (z^i_j (s), s) | > N^{\eps/2} B(z^i_j (s) ) \}
	\eeq
	where the infimum of the empty set is $+ \infty$. Then, we introduce
	\beq
	\tau := \min_{i,j} \tau_{ij} \wedge T.
	\eeq
	We want to prove that $\tau =T$.  First observe that since $|\tilf (z, t)| \leq N$ for any $z \in \Ed_t$, we have $|z_{s_1} - z_{s_2} | \leq C N |s_1-s_2|$ for any characteristic $z_s$ ending in some $\Ed_t$. It follows that for each $i$ and $s$ satisfying $t_i - TN^{-10} < s \leq  t_i$ that the points $\{ z_j^i (s) \}_{j=1}^{N^{10}}$ are a well-spaced grid of $\Ed_s$ in the sense that every $z \in \Ed_s$ is no further than $N^{-8}$ away from the closest of these points. Since $f(z, t)$ and $\tilf(z, t)$ are Lipschitz functions with Lipschitz constant less than $CN$ for $|z| > 1 + N^{-1}$ we see that for every $ s< \tau$ that,
	\beq
	| f (z, s) - \tilf (z, s) | \leq N^{\eps} B(z)
	\eeq
	for any $z \in \Ed_s$.  From Lemma \ref{lem:edge-no-eig} we conclude that at any time $s < \tau$ there are no eigenvalues of the form $\lambda = \e^{ \i \theta}$ for $\Theta_s + N^{-2/3+k(s)+5 \eps} < | \theta| \leq \pi$.
	
	As in the proof of Theorem \ref{thm:bulk} we write, for any $ t < t_i \wedge \tau$,
	\begin{align}
	f ( z^i_j (t) ) - \tilf ( z^i_j (t) ) = \E_1 ( t )^i_j + \E_2 (t)^i_j
	\end{align}
	where
	\beq
	\E_1 (t)^i_j =  - \frac{1}{2} \int_0^s z^i_j (s) (\del_z f) ( z^i_j (s), s) ( f ( z^i_j (s), s) - \tilf ( z^i_j (s), s) ) \d s
	\eeq
	and
	\beq
	\E_2(t)^i_j = \int_0^t \d M_s ( z^i_j (s) ).
	\eeq
We now prove the following estimate on the stochastic term.
\bep \label{prop:edge-bdg}
Let $\eps_1 >0$. Then for all $i, j$ as above,
\beq
\pp \left[ \exists t \in (0, t_i \wedge \tau ) : \left| \E_2 (t)^i_j \right| > N^{\eps_1+\eps/8} B ( z^i_j (t) ) \right] \leq C \log(N) \e^{ - c N^{\eps_1} }.
\eeq
\eep
\proof For simplicity of notation let us denote $z_t := z^i_j (t)$. We fix a sequence of intermediate times $0 <s_1, s_2, \dots, s_M = t_i$ such that,
\beq \label{eqn:edge-c4}
\frac{1}{2} B ( z_{s_k} ) \leq B ( z_{s_{k+1} } ) \leq 2 B ( z_{s_{k} } ) ,
\eeq 
for every $i$. We can take $M \leq C \log(N)$. Recall also the time $s_* < t_i$ defined in Proposition \ref{prop:edge-char}.  Define also $\kappa_t = \kappa ( z_t)$ and $\eta_t = \eta (z_t)$. 
We calculate the quadratic variation,
\begin{align}
\langle \bar{\E}_2 (s_k\wedge \tau ) \E_2 ( s_k \wedge \tau ) \rangle &\leq \frac{C}{N^2} \int_0^{s_k \wedge \tau} \frac{1}{N} \sum_{n=1}^N \frac{1}{ | \lambda_i (t) - z_t |^4} \d t \notag\\
&\leq \frac{C}{N^2} \int_{s_* \wedge \tau}^{s_k \wedge \tau} \frac{1}{N} \sum_{n=1}^N \frac{1}{ | \lambda_n - z_t|^4} \d t + \frac{C}{N^2}
\end{align}
where the integral in the last line is interpreted as $0$ if $s_k \leq s_*$.  When $s_k> s_*$ we continue to estimate the integral. For $s_* < t < \tau $ we have,
\begin{align}
 | \lambda_n (t) - z_t |^2 \geq c \left( ( \kappa_t -N^{-2/3+k(t) + 5\eps } )^2 + \eta_t^2 \right)
\end{align}
Note that in particular, we also used that $\kappa_t \geq N^{-2/3+k(t) + 5\eps } = D(t)$ for $s_* < t < s_k$ due to \eqref{eqn:edge-c1} and the fact that $\kappa_{s_k} \geq D(s_k)$. 
Therefore,
\begin{align} \label{eqn:general-u2}
\frac{1}{N^2} \int_{s_* \wedge \tau}^{s_k \wedge \tau} \frac{1}{N} \sum_{n=1}^N \frac{1}{ | \lambda_n - z_t|^4} \d t &\leq \frac{C}{N^2}  \int_{s_* \wedge \tau}^{s_k \wedge \tau}  \frac{ | \Re[ f (z_t, t) ] |}{ \eta_t \left( ( \kappa_t -N^{-2/3+k(t) + 5\eps } )^2 + \eta_t^2 \right)} \d  t \notag \\
&\leq \frac{C}{N^2}  \int_{s_* \wedge \tau}^{s_k \wedge \tau}  \frac{ | \Re[ \tilf (z_t, t) ] |}{ \eta_t \left( ( \kappa_t -N^{-2/3+k(t) + 5\eps } )^2 + \eta_t^2 \right)} \d  t ,
\end{align}
where in the second line we used that Lemma \ref{lem:edge-ref} together with the definition of the stopping time $\tau$ implies that for $t < \tau$ we have,
$$
| \Re[  f(z_t, t) ] | \leq | \Re [ \tilf (z_t, t) ] | + N^{\eps/2} B (z_t) \leq 2 | \Re[ \tilf (z_t, t) ] |.
$$
We use now, 
\begin{align}
\kappa_t - N^{-2/3+k(t) + 5\eps } &= \kappa_t - D(t)  \notag\\
&\geq \kappa_{s_k} + 2\mfb  \sqrt{ \kappa_{s_k}} (s_k - t) + \mfb^2 (s_k - t)^2 - D(t) \notag\\
&\geq \kappa_{s_k}  + 2 \mfb \sqrt{ \kappa_{s_k} } (s_k-t)+\mfb^2 (s_k - t)^2 \notag\\
 & - D(s_k) - \mfb \sqrt{ D(s_k ) } (s_k - t) - \frac{\mfb^2}{4} (s_k - t)^2 \notag\\
&\geq \mfb \sqrt{ \kappa_{s_k} } (s_k - t).
\end{align}
In the first inequality we used \eqref{eqn:edge-c3}. In the second inequality we used the square of the inequality $\sqrt{D (t) } \leq \sqrt{ D(s_k ) } + \frac{\mfb}{2} (s_k - t)$. In the last inequality we used $\kappa_{s_k} \geq D(s_k)$.

Therefore,
\begin{align}
\frac{1}{N^2} \int_{s_* \wedge \tau}^{s_k \wedge \tau}  \frac{ | \Re[ \tilf (z_t, t) ] |}{ \eta_t \left( ( \kappa_t -N^{-2/3+k(t) + 5\eps } )^2 + \eta_t^2 \right)}  \d t \leq \frac{C}{N^2}  \int_{s_*}^{s_k }  \frac{ | \Re[ \tilf (z_t, t) ] |}{ \eta_t \left( ( \kappa_{s_k} (t- s_k)^2 + \eta_t^2 \right)} \d t.
\end{align}
We need to consider two cases. First, consider $\eta_{s_k} \leq \kappa_{s_k}$.  From the fact that $\kappa (s_k) \geq D (s_k)$ and that $D(s_k) = N^{-\eps/6}$ if $s_k < \delta$ we see from \eqref{eqn:sqrroot} and \eqref{eqn:trivial-edge} that,
\beq
| \Re[ \tilf (z_{s_k}, s_k ) ] | \leq C \frac{ \eta_{s_k}}{ \sqrt{ \kappa_{s_k} } } N^{\eps/4}.
\eeq
Then using this, as well as that $\eta_t$ is decreasing and $\tilf (z_t, t)$ is constant along characteristics we have,
\begin{align}
\frac{1}{N^2} \int_{s_*}^{s_k }  \frac{ | \Re[ \tilf (z_t, t) ] |}{ \eta_t \left( ( \kappa_{s_k} (t- s_k) + \eta_t^2 \right)} \d t &\leq \frac{C N^{\eps/4}}{N^2} \frac{1}{ \sqrt{ \kappa_{s_k} + \eta_{s_k}}} \int_{s*}^{s_k} \frac{1}{ \kappa_{s_k} (t-s_k)^2 + \eta_{s_k }^2  } \d t \notag\\
& \leq \frac{C N^{\eps/4}}{N^2} \frac{1}{ \sqrt{ \kappa_{s_k} + \eta_{s_k} } } \frac{1}{ \sqrt{ \kappa_{s_k}} \eta_{s_k} } \leq \frac{C N^{\eps/4}}{N^2} \frac{1}{ ( \kappa_{s_k} + \eta_{s_k} )(\eta_{s_k} )}.
\end{align}
The second estimate follows via direct integration and in the last inequality we used the assumption $\kappa_{s_k} \geq \eta_{s_k}$. We now consider the case $\eta_{s_k} \geq \kappa_{s_k}$. In this case we proceed similarly to Proposition \ref{prop:bulk-bdg},
\begin{align} \label{eqn:general-u3}
\frac{1}{N^2} \int_{s_*}^{s_k }  \frac{ | \Re[ \tilf (z_t, t) ] |}{ \eta_t \left( ( \kappa_{s_k} (t- s_k) + \eta_t^2 \right)} \d t &\leq \frac{1}{N^2} \int_{s_*}^{s_k} \frac{| \Re [\tilf (z_t, t) ] |}{ \eta_t^3} \leq \frac{C}{N^2 \eta_{s_k}^2}\leq \frac{C N^{\eps/4}}{N^2 \eta_{s_k} ( \kappa_{s_k} + \eta_{s_k} )}.
\end{align}
By the BDG inequality,
\beq
\pp\left[ \sup_{ s \in (0, s_k )} \left| \E_2 ( s \wedge \tau ) \right| > N^{\eps_1+\eps/8} B (z_{s_k} ) \right] \leq C \e^{ - c N^{\eps_1} }.
\eeq
Taking a union bound over the $O ( \log(N) )$ choices of $s_k$ and using \eqref{eqn:edge-c4} we conclude the proof, similarly as in the proof of Proposition \ref{prop:bulk-bdg}. \qed

\vspace{5 pt}

\noindent{\bf Proof of Theorem \ref{thm:edge}.} With Proposition \ref{prop:edge-bdg} in hand, the proof of Theorem \ref{thm:edge} is similar to the proof of Theorem \ref{thm:bulk}. Let $\Upsilon$ be the event of Proposition \ref{prop:edge-bdg}, after taking an intersection over all choices of $i, j \leq N^{10}$. Let $z_t := z(t)^i_j$ for notational simplicity. On the event of $\Upsilon$ we have the inequality for $0 < t < \tau$,
\beq
| f (z_t, t) - \tilf (z_t, t) | \leq \int_0^t g(s) | f(z_s, s)  - \tilf (z_s, s) | \d s + N^{\eps_1+\eps/8} B ( z_t ),
\eeq
where we denoted,
\beq
g(s) = \frac{ | z_s  (\del_z f ) (z_s, s) |}{2}.
\eeq
Hence, via Gronwall's inequality we obtain for $0 < t < \tau$,
\beq
| f (z_t, t) - \tilf (z_t, t) | \leq \int_0^t g(s) \exp \left[ \int_s^t g(u) \d u \right] N^{\eps_1+\eps/8} B ( z_s ) \d s + N^{\eps_1+\eps/8} B (z_t).
\eeq
Similar to the proof of Theorem \ref{thm:bulk}, now using Lemma \ref{lem:edge-ref}, we have
\beq
\exp \left[ \int_s^t g(u) \d u \right] \leq C \frac{ \log |z_s |}{ \log |z_t|}
\eeq
as well as
\beq
g(s) \leq C \frac{ \Re[ \tilf (z_s, s) ]}{ \eta_s}.
\eeq
Therefore,
\beq
| f (z_t, t) - \tilf (z_t, t) | \leq C \frac{ N^{\eps_1+\eps/8}}{\eta_t} \int_0^t | \Re [ \tilf (z_s, s) ]  | B ( z_s ) \d s + N^{\eps_1+\eps/8} B(z_t).
\eeq
We must consider a few different cases. First, let us consider the case that $t < s_*$. Then, we see that
\beq
 \frac{ N^{\eps_1+\eps/8}}{\eta_t} \int_0^t | \Re [ \tilf (z_s, s) ]  | B ( z_s ) \d s + N^{\eps_1+\eps/8} B(z_t) \leq C \frac{N^{\eps_1+\eps/8}}{N}.
\eeq
So in the remainder of the proof we consider the case $t \geq s_*$. 
We now consider a few different cases depending on whether or not $t_i$ (the end-time of the characteristic $z_t$) is small or large. First, assume that $t_i < \delta$. Then, using \eqref{eqn:trivial-edge},
\beq
\frac{ | \Re[ \tilf (z_s, s) ] | }{ \eta_t }  = \frac{ | \Re[ \tilf (z_t, t) ] | }{ \eta_t } \leq C\left( |\kappa_t|^{-2} \right) \leq C N^{\eps/3}
\eeq
we have,
\beq
 \frac{ N^{\eps_1+\eps/8}}{\eta_t} \int_0^t | \Re [ \tilf (z_s, s) ]  | B ( z_s ) \d s \leq N^{11 \eps/24+\eps_1} B(z_t)
\eeq
because $B(z_t)$ is effectively decreasing along characteristics (it may not be decreasing for $t < s_*$ but for such $t$ it is $\O(N^{-1})$). We take $\eps_1 < \eps/100$. 

Now, we assume that $t_i > \delta$. Then for $ s < \delta/2$ we have that either $\eta_s \geq c$ or $\kappa_s \geq c$ by Proposition \ref{prop:edge-char} depending on whether $s$ is smaller or larger than $s_*$. Then for such $s$ we have $B (z_s) \leq CN^{-1}\eta_t^{-1/2}$ and so,
\beq
\frac{N^{\eps_1+\eps/8}}{\eta_t} \int_0^{t\wedge \delta/2} | \Re[ \tilf (z_s, s) ] | B(z_s) \d s \leq C N^{\eps_1+\eps/8} \frac{1}{N \eta_t^{1/2}} \frac{  | \Re[\tilf (z_t, t) ] |}{\eta_t } \leq C N^{\eps_1+\eps/8} B (z_t).
\eeq
In the final inequality we used $\Re[\tilf (z_t, t) ] \eta_t^{-1} \leq C ( \kappa_t + \eta_t)^{-1/2}$ (recall that $t \geq s_*$) if $ t> \delta/2$ and \eqref{eqn:sqrroot}. If $t < \delta/2$ then since $t \geq s_*$ we have $\kappa_t \geq c$ since $t < t_i - \delta/2$ and so $| \Re[ \tilf (z_t, t) ]| \eta_t^{-1}$ is bounded.

Now we have still to estimate the integral over $[t\wedge \delta/2, t]$ in the case that $t_i > \delta$.  Since this integral is $0$ if $t < \delta/2$ we may assume that $t > \delta/2$. Then, using freely \eqref{eqn:sqrroot}, we have,
\begin{align}
\frac{N^{\eps_1+\eps/8}}{\eta_t} \int_{\delta/2 \vee s_*}^t | \Re[\tilf (z_s, s) ] | B(z_s) \d s & \leq \frac{C}{N} \frac{ N^{\eps_1+\eps/8} }{ \sqrt{ \kappa_t + \eta_t } } \int_{\delta/2 \vee s_*}^t \frac{ \eta_s}{ \sqrt{ \kappa_s + \eta_s } \eta_s^{3/2}}\d s \notag\\
&\leq \frac{C}{N} \frac{N^{\eps_1+\eps/8}}{\sqrt{ \kappa_t + \eta_t } } \int_{\delta/2 \vee s_*}^t \frac{ | \Re[ \tilf (z_s, s) ]|  } {\eta_s^{3/2}}\d s \notag\\
& \leq \frac{C N^{\eps_1+\eps/8}}{N \sqrt{ \kappa_t + \eta_t} \sqrt{\eta_t}} 
\end{align}
The integral over $[\delta/2, \delta/2 \vee s_*]$ contributes $N^{\eps_1+\eps/8-1} ( \kappa_t + \eta_t)^{-1/2}$ because $\eta_s \geq c$ here. Therefore, taking $\eps_1 < \eps/100$ we see that we have proven that,
\beq
| f (z_t, t) - \tilf (z_t, t) | \leq C N^{23 \eps/48} B(z_t) \ll N^{\eps/2} B(z_t).
\eeq
on the event $\Upsilon$. It follows that $\tau = T$. \qed

\section{Cusp estimates} \label{sec:cusp}

In this section we prove Theorem \ref{thm:cusp}. We will use the following reverse time parameterization. We will denote by uppercase $T$, and $S$ the usual forward time parameterization, so that we will consider $S$ and $T$ close to $4$.  We will the introduce,
\beq
T = 4 - t, \qquad S = 4 - s
\eeq
so that $t$ and $s$ will usually obey,
\beq
N^{-1/2+\delta} \leq s, t \leq \frac{1}{10}.
\eeq
We will often substitute in $T$ or $t$ into functions of time. When we write $t$ it is understood that we are evaluating the function at $4-t$. For example, the gap between edges is,
\beq
\Delta_t = \Delta_T = 2( \pi - \Theta_T) = \frac{1}{3} t^{3/2} ( 1 + \O (t) ).
\eeq
The proof of Theorem \ref{thm:cusp} is similar in structure to the proof of Theorem \ref{thm:edge}. We first establish analogues of the estimates proven there before proceeding to the main body of the proof.
\bel \label{lem:ref-cusp}
 There is a constant $C>0$ so that the following holds.  For $0 < t < \frac{1}{10}$ and $z = (1+ \eta) \e^{ \i (\Theta_t + \kappa ) }$ with $\eta$ and $\kappa$ satisfying,
\beq
0 < \eta < \Delta_t, \qquad 0 < \kappa < \pi - \Theta_t
\eeq
we have that
\beq
\frac{1}{C} \frac{ \eta}{ \sqrt{ \eta + \kappa} } \leq \Delta_t^{1/6} | \Re[ \tilf (z, t) ] | \leq C \frac{ \eta}{ \sqrt{ \eta + \kappa} }
\eeq
\eel
\proof Denoting $\theta = \kappa + \Theta_t$ we have,
\begin{align}
- \Re[ \tilf (z,t) ] \asymp \eta \int \frac{ \rho_t (x) }{ \eta^2 + \sin^2 (\frac{ \theta - x}{2}) } \d x.
\end{align}
Due to symmetry of the density,
\beq
 \int \frac{ \rho_t (x) }{ \eta^2 + \sin^2 (\frac{ \theta - x}{2}) } \d x \asymp   \int_0^\pi \frac{ \rho_t (x) }{ \eta^2 + \sin^2 (\frac{ \theta - x}{2}) } \d x.
\eeq
Changing coordinates and using the fact that $\sin^2(x) \asymp x^2$ for $|x| \leq \pi/2$ we have,
\beq
 \int_0^\pi \frac{ \rho_t (x) }{ \eta^2 + \sin^2 (\frac{ \theta - x}{2}) } \d x \asymp \int_{0}^{\Theta_t} \frac{ \rho_t (\Theta_t - x ) }{ \eta^2 + ( \kappa + x)^2} \d x.
\eeq
For an upper bound we split the integral into the regions $[0, \Delta_t]$ and $[\Delta_t, \Theta_t]$. For the latter region we can bound $\rho_t(\Theta_t -x  ) \leq C x^{1/3}$,
\begin{align}
\int_{\Delta_t}^{\Theta_t} \frac{ \rho_t (\Theta_t - x ) }{ \eta^2 + ( \kappa + x)^2} \d x \leq C \int_{\Delta_t}^{\Theta_t} ( \Delta_t + x)^{-5/3} \d x \leq C \Delta_t^{-2/3} \leq C \Delta_t^{-1/6} \frac{1}{ \sqrt{ \kappa + \eta} } ,
\end{align}
where in the last inequality we used the assumption $\kappa + \eta \leq 2 \Delta_t$.  In the region $[0, \Delta_t ]$ we use $\rho_t (\Theta_t - x) \Delta_t^{1/6} \leq C x^{1/2}$,
\begin{align}
\int_0^{\Delta_t} \frac{ \rho_t (\Theta_t - x ) }{ \eta^2 + ( \kappa + x)^2} \d x  \leq & C \Delta_t^{-1/6} \int_0^{\Delta_t} \frac{ x^{1/2}}{ \eta^2 + ( \kappa + x)^2} \d x \notag \\
\leq & C \Delta_t^{-1/6} \int_0^{\Delta_t} ( x + \kappa + \eta )^{-3/2} \d x \notag\\
\leq & C \Delta_t^{-1/6} \frac{1}{ \sqrt{ \kappa + \eta } }.
\end{align}
This completes the proof of the upper bound. For the lower bound, 
\begin{align}
\int_{0}^{\Theta_t} \frac{ \rho_t (\Theta_t - x ) }{ \eta^2 + ( \kappa + x)^2} \d x  & \geq  \int_0^{ 2 \Delta_t }  \frac{ \rho_t (\Theta_t - x ) }{ \eta^2 + ( \kappa + x)^2} \d x \notag\\
& \geq c \Delta_t^{-1/6} \int_0^{2 \Delta_t} \frac{ \sqrt{x}}{ \eta^2 + ( \kappa + x)^2} \d x \notag\\
&\geq c \Delta_t^{-1/6} \int_{ \frac{\kappa+\eta}{2}}^{ \kappa+\eta}  \frac{ \sqrt{x}}{ \eta^2 + ( \kappa + x)^2} \d x \notag\\
& \geq c \Delta_t^{-1/6} \int_{ \frac{\kappa+\eta}{2}}^{ \kappa+\eta} (\kappa+\eta)^{-3/2} \d x \geq c  \Delta_t^{-1/6} \frac{1}{ \sqrt{ \kappa + \eta}}.
\end{align}
We used above the assumption that $\kappa+\eta \leq 2 \Delta_t$. This completes the proof. \qed

\bel \label{lem:cusp-char-1}
Recall the convention $T= 4-t$. There are $\mfa, \mfb >0$ so that the following holds. Let $z_s = \C_{S, T} (z)$ denote a characteristic ending at a point $ z = (1+ \eta) \e^{ \i (\Theta_t + \kappa ) }$ where $0 < \eta < \mfa$ and $0 < \kappa < \Delta_T/2$. Let $S_*$ be defined as,
\beq
S_* = \sup \{ S : |z_S| \geq 1 + \mfa \Delta_S \}.
\eeq
Then for $S_* \vee (4- 10^{-1} ) < S < T$ we have,
\beq \label{eqn:cusp-kappa}
\Delta_S^{1/6} \sqrt{ \kappa_S } > \Delta_T^{1/6} \sqrt{ \kappa_T} + \mfb (T-S) ,
\eeq
and also that $\kappa_S$ is decreasing. 
For $S<S_*$ we have $\eta_S \geq \eta_{S_*} \geq \mfa \Delta_{S_*} \geq \mfa \Delta_T$.
\eel
\proof We let $100 \mfa$ be the constant from Proposition \ref{prop:a-cusp}. For $S \in (S_*, T)$ we can apply the estimate from this proposition to obtain,
\begin{align}
\frac{ \d}{ \d S} \kappa_S &= \frac{1}{2} \left( \Im [ \tilf (z_S, S) ]  - \Im[ \tilf (\e^{ \i \Theta_S}, S) ) ] \right) \notag\\
& \leq - c \Delta_S^{-1/6} \sqrt{ \kappa_S},
\end{align}
as long as $\kappa_S \geq 0$. We see that $\kappa_S$ is decreasing so if $\kappa_T >0$ then $\kappa_S > 0$ for $S > S_* \vee (4- 10^{-1})$. Moreover, since $\Delta_S$ is decreasing, we see that
\beq
\frac{ \d }{ \d S} \left( \Delta^{1/3}_S \kappa_S \right) \leq - c \left( \Delta_S^{1/3} \kappa_S \right)^{1/2}.
\eeq
The differential inequality $\del_t g \leq -c  g^{1/2}$ is solved by considering $\del_t g^{1/2} \leq -c/2$ and so we see that for $S \in (S_* \vee (4- 10^{-1} ), T)$ that
\beq
\Delta_S^{1/6} \sqrt{ \kappa_S} \geq \Delta_T^{1/6} \sqrt{ \kappa_T} + c (T-S)
\eeq
for some $c>0$ 
as desired. This completes the proof. \qed

The following establishes estimates on the behavior of the characteristics near the spectral edge.

\bep \label{prop:cusp-char}
Let $z_s = \C_{s, t} (z)$ be a characteristic ending at a point $z$ at final time $t < 10^{-1}$ (recall the reversed time convention $T = 4-t$). Assume at the final time $t$ the estimates,
\beq
\kappa_t > N^{-2/3+5 \eps} \Delta_t^{1/9}, \qquad \Delta_t \geq N^{-3/4+9\eps}
\eeq
hold. Moreover assume $0 < \eta_t < \mfa \Delta_t$. Let,
\beq
s_* = \inf \{ s : \eta_s > \mfa \Delta_s \}.
\eeq
Note $S_* = 4 - s_*$ where $S_*$ is as in Lemma \ref{lem:cusp-char-1}. 
For $t < s < s_* \wedge 10^{-1}$ we have,
\beq \label{eqn:posder}
\del_s ( \kappa_s - N^{-2/3+5 \eps} \Delta_s ) \geq 0 .
\eeq
For the next two estimates, continue to assume $t < s < s_* \wedge 10^{-1}$. 
There is furthermore a constant $\mfd >0$ so that the following holds.  If $s-t \leq \Delta^{1/6}_t \sqrt{ \kappa_t} $ then,
\beq \label{eqn:first-kappa}
 (\kappa_s - N^{-2/3 + 5 \eps} \Delta^{1/9}_s)^2 \geq \mfd (s-t)^2 \Delta^{-1/3}_t \kappa_t. 
\eeq 
If $s-t \geq \Delta^{1/6}_t \sqrt{ \kappa_t}$ then,
\beq \label{eqn:second-kappa}
(\kappa_s - N^{-2/3 + 5\eps} \Delta_s^{1/9})^2 \geq \mfd \kappa_t^2.
\eeq
\eep
\proof By direct calculation, one can see that $| \del_t \Theta_t | \leq C t^{1/2} \leq C \Delta_t^{1/3}$. Therefore, proceeding as in the proof of Lemma \ref{lem:cusp-char-1} we have for $t < s < s_*$,
\beq \label{eqn:someder}
\frac{ \d}{ \d s} ( \kappa_s - N^{-2/3+5 \eps} \Delta_s^{1/9} ) \geq c \Delta^{-1/6}_s \sqrt{ \kappa_s} - C N^{-2/3+5\eps} \Delta^{-5/9}_s,
\eeq
for some $c, C>0$. Consider the set
\beq
\A := \{ s: t \leq s < s_* : \kappa_s < N^{-2/3+5\eps} \Delta_s^{1/9} \}.
\eeq
Note that by assumption there is a small interval around $t$ not contained in $\A$. Assume that $\A$ is not empty and let $F = \inf \A$, so that $s_* > F>t$. At $F$ we clearly must have that
\beq
\kappa_F = N^{-2/3+5\eps} \Delta_F^{1/9}.
\eeq
We now show that the RHS of \eqref{eqn:someder} is strictly positive at $s= F$. Indeed, we evaluate
\begin{align}
\Delta_F^{10/9-1/3}\kappa_F = N^{-2/3+5\eps} \Delta_F^{8/9} \geq N^{-2/3+5\eps} \Delta_t^{8/9} \geq N^{3 \eps} ( N^{-2/3+5\eps} )^2
\end{align}
which shows that the RHS of \eqref{eqn:someder} is strictly positive for $N$ large enough. In particular, the derivative on the LHS of \eqref{eqn:someder} is strictly positive at $s=F$, showing that $\kappa_s \geq N^{-2/3+5\eps} \Delta_s^{1/9}$ in an open interval containing $F$. This contradicts the assumption that $F$ was an infimum, proving that $\A$ is empty. 

Substituting the lower bound $\kappa_s \geq N^{-2/3+5 \eps} \Delta_s^{1/9}$ into the RHS of \eqref{eqn:someder} we then find that,
\begin{align}
\frac{ \d}{ \d s} ( \kappa_s - N^{-2/3+5 \eps} \Delta_s^{1/9} )  \geq c \Delta_s^{-1/9} N^{-1/3+5\eps/2} - C N^{-2/3+5 \eps} \Delta_s^{-5/9}.
\end{align}
Positivity of the RHS is equivalent to 
\beq
\Delta_s^{4/9} \geq \frac{C}{c} N^{-1/3+5\eps/2}.
\eeq
However, by assumption $\Delta_t^{4/9} \geq N^{-1/3+4 \eps}$. We conclude the proof of \eqref{eqn:posder}.   

We turn now to the proof of the remaining inequalities. From \eqref{eqn:cusp-kappa} we have
\begin{align}
    \Delta^{1/3}(s) \kappa(s) - N^{-2/3 + 5 \varepsilon} \Delta^{4/9}(s) 
    \geq & [\Delta^{1/3}(t) \kappa(t) - N^{-2/3 + 5 \varepsilon} \Delta^{4/9}(t)] \nonumber \\ +&  [N^{-2/3 + 5 \varepsilon} \Delta^{4/9}(t) - N^{-2/3 + 5 \varepsilon } \Delta^{4/9}(s)] \nonumber \\
    +& 2 c(s-t) \Delta^{1/6}(t) \sqrt{\kappa(t)},
\end{align}
for some $c >0$.  Note that the first term on the RHS is positive. First, assume that $s-t \leq \Delta_t^{1/6} \sqrt{ \kappa_t}$. From the fact that $\kappa_t \leq \Delta_t$ and $\Delta_t \leq C t^{3/2}$ we see that $s \leq C t$ for some $C>0$. It follows from the mean value theorem that,
\begin{align}
N^{-2/3+5\eps} | \Delta^{4/9}_t - \Delta_s^{4/9} | & \leq C N^{-2/3+5\eps} (s-t) \Delta^{-2/9}_t \notag\\
& \leq C (s-t) \Delta_t^{1/6} \sqrt{ \kappa_t} N^{-3 \eps/2}.
\end{align}
The estimate \eqref{eqn:first-kappa} easily follows.  Consider now the case $s-t \geq \Delta_t^{1/6} \sqrt{ \kappa_t}$. Let $t_*$ be the time such that
\beq
t_*  = t +  \Delta_t^{1/6} \sqrt{ \kappa_t}.
\eeq
Applying \eqref{eqn:posder} we see that
\beq
\kappa_s - N^{-2/3+5 \eps} \Delta_s \geq \kappa_{t_*} - N^{-2/3+5\eps} \Delta_{t_*}.
\eeq
Since the estimate \eqref{eqn:first-kappa} holds at $s = t_*$ we see that
\beq
( \kappa_{t_*} - N^{-2/3+5\eps} \Delta_{t_*} )^2 \geq  c \kappa_t^2
\eeq
for some $c>0$. This completes the proof of the proposition. \qed

The proof of the following is similar to Lemma \ref{lem:edge-no-eig} and appears in Appendix \cite{a:cusp-no-eig-proof}. 
\bel  \label{lem:cusp-no-eig} Let $\eps >0$. Consider the domain,
\beq
\D := \{ z = (1+\eta ) \e^{ \i \theta} : \Delta_t^{1/9} N^{-2/3+\eps} < \eta < 2 \Delta_t^{1/9} N^{-2/3+\eps} , \Theta_t + N^{-2/3+5\eps}  \Delta_t^{1/9} < | \theta | \leq \pi \}.
\eeq
Assume for all $ z \in \D$ we have the estimate,
\beq
| \tilf (z, t)   - f(z, t) | \leq \frac{ N^{\eps}}{ N \sqrt{ \eta} \sqrt{ \kappa+\eta} }.
\eeq
Then there are no eigenvalues of the form $\lambda = \e^{ \i \theta}$ with $\Theta_t + N^{-2/3+5\eps} \Delta_t < | \theta | \leq \pi $.
\eel

\bel \label{lem:cusp-good-char}
Let $\eps >0$ and let $t < N^{-1/10}$ with $\Delta_t \geq N^{-3/4+9\eps}$. Let $z_s$ be a characteristic that ends at time $t$ in the region,
\beq
\{ z= (1+\eta) \e^{ \i \theta} : N^{-2/3+\eps} \Delta_t^{1/9} \leq \eta \leq 2N^{-2/3+\eps} \Delta_t^{1/9}, \Theta_t + N^{-2/3+5\eps} \Delta_t^{1/9} \leq | \theta| \leq \pi \}.
\eeq
Then for $t < s < 10^{-1}$ we have,
\beq
\eta_s \leq C N^{-\eps/2} \Delta_s .
\eeq
\eel
\proof Define the functions
\beq
h_1 (u) = \log |z_t| + u N^{\eps/2} \frac{ \eta_t}{ \Delta_t^{1/6} \sqrt{ \kappa_t} }
\eeq
and
\beq
h_2 (u) = N^{-\eps/2} \Delta_t + u t^{1/2} N^{-\eps}. 
\eeq
By the definition of the characteristics,
\beq
\log |z_s| = \log |z_t| + (s-t) \frac{ | \Re[ \tilf (z, t) ] |}{2},
\eeq
and so for $0 \leq u < 10^{-1}$ we have,
\beq
\eta_{t+u} \leq C  h_1 (u) 
\eeq
for some large $C$, using Lemma \ref{lem:ref-cusp}. Since $\del_u \Delta_u \geq c u^{1/2}$ we have that
\beq
h_2 (u) \leq N^{-\eps/2} \Delta_{t+u} ,
\eeq
for $N$ large enough. Note that $h_2(0) > h_1 (0)$. We claim for $u  < 10^{-1}$ that
\beq
h_2'(u) \geq h_1' (u)
\eeq
which will yield the claim. This is equivalent to,
\beq
N^{3\eps/2} \eta_t \leq t^{1/2} \Delta_t^{1/6} \sqrt{ \kappa_t}.
\eeq
The LHS is less than $2 N^{5\eps/2-2/3} \Delta_t^{1/9}$. The RHS is larger than $c N^{5\eps/2-1/3} \Delta_t^{1/2} \Delta_t^{1/18}$.  But,
\beq
\Delta_t^{1/2-1/9+1/18} = \Delta_t^{4/9} \geq N^{-1/3+4 \eps}
\eeq
by assumption. \qed

\bel \label{lem:cusp-qv}
Fix a time $t$. Let $z_s = \C_{s, t} (z)$ be a characteristic terminating at a point $z = (1+\eta) \e^{ \i ( \Theta_t + \kappa)}$ where
\beq
N^{-2/3+\eps} \Delta_t^{1/9} < \eta < 2 N^{-2/3+\eps} \Delta_t^{1/9}, \qquad N^{-2/3+5 \eps} \Delta_t^{1/9} \leq \kappa < \pi - \Theta_t.
\eeq
Let $t_1$ and $t_2$ be two times $t \leq t_1 < t_2 < 10^{-1}$ and assume that for $t_1 < s < t_2$ the estimate
\beq \label{eqn:cusp-qv-condition}
| \tilf (z_s, s) - f (z_s, s) | \leq 2 | \Re[ \tilf (z_s, s) ] |
\eeq
holds and that at each time $s$ there are no eigenvalues of the form $\lambda = \e^{ \i \theta}$ for $\Theta_s + N^{-2/3+5\eps} \Delta_s^{1/9} < | \theta | \leq \pi $. Then,
\begin{align}
\int_{t_1}^{t_2} \frac{1}{N} \sum_{i=1}^N \frac{1}{ | \lambda_i (s) - z_s |^4}  \d s \leq C \log(N)  B (z_{t_1} )^2
\end{align}
\eel
\proof  First we consider the case that $ \mfa \kappa_{t_1}  \leq 10 \eta_{t_1}$. In this case, it suffices to estimate the integrand by $C \log(N) \eta_{t_1}^{-2}$. This argument is similar to that appearing in the proof of Proposition \ref{prop:bulk-bdg} and so we omit it. 

For the remainder of the proof we therefore assume that $\kappa_{t_1} \mfa \geq 10 \eta_{t_1}$.  In particular, this implies that $\eta_{t_1} \leq \mfa \Delta_{t_1} /10$. 

 Let $s_*$ be as in Lemma \ref{lem:cusp-char-1}, that is
 \beq
s_* :=  \inf \{ s > t : \eta_s > \mfa \Delta_s \}.
 \eeq
 Note that by Lemma \ref{lem:cusp-good-char} we have $s_* > t_2$. 
   For $t_1 < s  <t_2$ we know from \eqref{eqn:posder} that $\kappa_s \geq N^{-2/3+5\eps} \Delta_s^{1/9}$ because this is satisfied at $s=t$. In particular, define
\beq
\tilde{t} = t_1 + \Delta_{t_1}^{1/6} \sqrt{ \kappa_{t_1} }.
\eeq  
We assume $ t_1 < \tilde{t} <  t_2$. The other cases are easy to deal with, as one has to treat only one of the regions of integration described below.   
Then for $t_1 < s < \tilde{t} $ we have,
\beq \label{eqn:a-first-kappa}
( \kappa_s - N^{-2/3+5\eps} \Delta^{1/9}_s )^2 \geq \mfd (s-t_1)^2 \Delta^{-1/3}_{t_1} \kappa_{t_1}
\eeq
and for $\tilde{t}  < s <  t_2$,
\beq \label{eqn:a-second-kappa}
( \kappa_s - N^{-2/3+5 \eps} \Delta_s^{1/9} )^2 \geq \mfd \kappa_{t_1}^2,
\eeq
by applying Proposition \ref{prop:cusp-char}.  We will use the estimate,
\begin{align} \label{eqn:general-u4}
\int_{t_1}^{ t_2} \frac{1}{N} \sum_{i=1}^N \frac{1}{ | \lambda_i (s) - z_s|^4} \d s  & \leq \int_{t_1}^{ t_2} \frac{1}{ ( \kappa_s - N^{-2/3+5\eps} \Delta_s^{1/9} )^2 + \eta_s^2} \frac{ | \Re[ \tilf (z_s, s) ] | }{ \eta_s} \d s .
\end{align}
  We start with the region $ s > \tilde{t}$. For such $s$ we can apply \eqref{eqn:a-second-kappa} and estimate,
\begin{align}
\int_{\tilde{t}}^{ t_2} \frac{1}{ ( \kappa_s - N^{-2/3+5\eps} \Delta_s^{1/9} )^2 + \eta_s^2} \frac{ | \Re[ \tilf (z_s, s) ] | }{ \eta_s} \d s & \leq \frac{C}{ \kappa_{t_1}^2} \int_{\tilde{t}}^{ t_2} \frac{ | \Re[ \tilf (z_s, s) ] | }{ \eta_s} \d s \notag\\
& \leq \frac{C}{ \kappa_{t_1}^2} \log(N) \leq \frac{C}{ \eta_{t_1} ( \kappa_{t_1} + \eta_{t_1})} \log(N)
\end{align}
 We consider now the region $s < \tilde{t}$. By applying \eqref{eqn:a-first-kappa} and Lemma \ref{lem:ref-cusp} (as well as the constancy of $\tilf$ along characteristics, and the fact that $\eta_s$ is increasing in $s$) we obtain,
\begin{align} \label{eqn:general-u5}
\int^{\tilde{t}}_{ t_1} \frac{1}{ ( \kappa_s - N^{-2/3+5\eps} \Delta_s^{1/9} )^2 + \eta_s^2} \frac{ | \Re[ \tilf (z_s, s) ] | }{ \eta_s} \d s & \leq \frac{C}{\Delta_{t_1}^{1/6} \sqrt{ \kappa_{t_1}+\eta_{t_1}} } \int_{t_1}^{\tilde{t}} \frac{1}{(s-t)^2  \Delta_{t_1}^{-1/3} \kappa_{t_1} + \eta_{t_1}^2} \d s \notag\\
&\leq \frac{C}{ \eta_{t_1} \sqrt{ \kappa_{t_1}} \sqrt{ \kappa_{t_1} + \eta_{t_1}}}.
\end{align}
This completes the proof. \qed

Define the domain,
\beq
\fB_s :=\{ z = (1+\eta) \e^{ \i \theta} : N^{-2/3+\eps} \Delta^{1/9}_s < \eta < 2  N^{-2/3+\eps} \Delta^{1/9}_s , N^{-2/3+5\eps} \Delta_s^{1/9} + \Theta_s < | \theta | \leq \pi \} .
\eeq
We will assume $s > N^{-1/2+9\eps}$ so that $\Delta_s > c N^{-3/4+10\eps}$.  We fix $S_0 = 4 - s_0$ with $s_0 = 10^{-1}$, and the final time,
\beq
S_f := 4 - s_f, \qquad s_f = N^{-1/2+9\eps}
\eeq
Similar to the proof of Theorem \ref{thm:edge}, we introduce a polynomial number of characteristics as follows. We introduce times $T_i$ by
\beq
T_0 = S_0, \qquad T_i = T_{i-1} + \frac{S_f-S_0}{N^{10}},
\eeq
for $i=1, \dots, N^{10}$. At each $T_i$ we introduce a well-spaced mesh of $\B_{T_i}$ of size $N^{10}$ denoted by $\{ u_j^i \}_{j=1}^{10}$.  We introduce the characteristics,
\beq
z_j^i (S) := \C_{S, T_i} ( u_j^i ).
\eeq
Note that each characteristic is defined only for $ 0 \leq s \leq T_i$. For $i \geq 1$ we introduce the stopping times,
\beq
\tau_{ij} = \inf \{ S \in (S_0, T_i ) : | f (z_j^i (S), S) - \tilf (z_j^i (S), S) | > N^{\eps/2} B (z_j^i (S) ) \},
\eeq
with the infimum of the empty set being $+ \infty$. Let $\A_0$ be the event,
\beq
\A_0 = \{ \exists (i, j) : | f ( z_i^j (S_0), S_0 ) - \tilf ( z_i^j (S_0), S_0 ) | > N^{\eps_1} \inf_{ S_0 < S < T_i } B ( z_i^j (S) ) \}^c
\eeq
\bel
The event $\A_0$ holds with overwhelming probability.
\eel
\proof We will show that the desired estimates hold on the event of Proposition \ref{prop:a-extended}, with $\eps, \delta$ in that proposition statement chosen sufficiently small.

 Fix a single characteristic $z_S := z_i^j (S)$ and let $S_0 < S < T_i$ be fixed. By \eqref{eqn:posder} and Lemma \ref{lem:cusp-good-char} it follows that $\kappa_S >0$ for all $S$ and that $\kappa_{S_0} \geq c N^{-2/3+5 \eps}$. From \eqref{eqn:ext-ll-2} we see that 
 \beq
 | f ( z_{S_0}, S_0 ) - \tilf ( z_{S_0} , S_0 ) | \leq N^{\eps_1} B( z_{S_0} ).
 \eeq
 Since $\kappa_S >0$ and $\kappa_S$ and $\eta_S$ are both decreasing in $S$ that,
 \beq
 B(z_S) \geq B( z_{S_0} )
 \eeq
 for $S  > S_0$. This yields the claim. \qed

We also let $\tau_0$ be the stopping time that equals $+\infty$ on $\A_0$ and is $S_0$ on $\A_0^c$. We now introduce the stopping time,
\beq
\tau = \left( \min_{i, j} \tau_{ij} \right) \wedge \tau_0 \wedge S_f.
\eeq
Using the Lipschitz continuity of $\tilf (z, S)$ and $f(z, S)$ on the domains $\B_S$ we have from Lemma \ref{lem:cusp-no-eig} that for any $S_0 < S < \tau$ that there are no eigenvalues of the form $\lambda = \e^{ \i \theta}$ for $\Theta_S + N^{-2/3+5 \eps} \Delta_S^{1/9} < | \theta | \leq \pi$. 

\bel \label{lem:cusp-small-error}
Let $z_S$ be a characteristic terminating at time $T_i$ in the domain,
\beq
\{ z = (1+ \eta) \e^{ \i \theta} : N^{-2/3+\eps} \Delta_{T_i}^{1/9} \leq \eta \leq \mfa \Delta_{T_i} / 10, \Theta_{T_i} + N^{-2/3+5\eps} \Delta_{T_i}^{1/9} \leq |\theta | \leq \pi \}.
\eeq
 Then for $S_0 < S < T_i$ we have,
\beq
N^{\eps} B(z_S) \leq \frac{1}{ \log(N)} | \Re[ \tilf (z_S, S) ]|.
\eeq
\eel
\proof Note that $\tilf (z_S, S)$ is constant. At $S = T_i$ we have,
\beq
N^{\eps} B(z_{T_i} )  = N^{\eps} \frac{1}{ N \sqrt{ \eta_{T_i} (\eta_{T_i} + \kappa_{T_i})  }} \leq C \frac{N^{\eps} \Delta^{1/6}_{T_i} }{ N \eta_{T_i}^{3/2}} | \Re[ \tilf (z_{T_i} , T_i ) ] | \leq C N^{-\eps/2} | \Re[ \tilf (z_{T_i}, T_i ) ] |
\eeq
where we used $\eta_{T_i} \geq \Delta_{T_i}^{1/9} N^{-2/3+\eps}$ which holds by assumption.  Let $S_*$ be as in Lemma \ref{lem:cusp-char-1}. By Lemma \ref{lem:cusp-char-1} we have
\beq
B( z_S) \leq B (z_{T_i} )
\eeq
for all $S > S_*$. If $S_* > S_0$ then for $S < S_*$ we have $\eta_S \geq \eta_{S_*} \geq c \Delta_{S_*} \geq c \Delta_{T_i} \geq c \kappa_{T_i}$ and so 
\beq
B(z_S) \leq (N \eta_S)^{-1} \leq C B (z_{T_i} ).
\eeq
This completes the proof. \qed

With similar notation as in the other sections we have for any $S_0 \leq S \leq T_i \wedge \tau$,
\begin{align} \label{eqn:cusp-ito}
f ( z_j^i (S) , S) - \tilf (z_j^i (S), S) &= f ( z_j^i (S_0) , S_0) - \tilf (z_j^i (S_0), S_0)  \notag\\
= & \E_1 (S)_j^i + \E_2 (S)_j^i ,
\end{align}
where
\beq
\E_1 (S)_j^i = - \frac{1}{2} \int_{S_0}^S z_j^i (U) ( \del_z f) ( z_j^i (U), U) ( f ( z_j^i (U), U) - \tilf ( z_j^i (U), U) ) \d U
\eeq
and
\beq
\E_2(S)_j^i = \int_{S_0}^S \d M_U ( z_j^i (U) ).
\eeq
For the martingale term we have the following.
\bep \label{prop:cusp-bdg}
For any $i, j$ and $\eps_1 >0$ we have,
\beq
\pp \left[ \exists S \in (S_0, T_i \wedge \tau ) : | \E_2(S)_j^i | > N^{\eps_1} B(z_j^i (S) ) \right] \leq C \log(N) \e^{ - c N^{\eps_1} }.
\eeq
\eep
\proof This is proven in an almost identical manner to Proposition \ref{prop:edge-bdg}. The quadratic variation is bounded using Lemma \ref{lem:cusp-qv}.  Note that the condition \eqref{eqn:cusp-qv-condition} is a consequence of Lemma \ref{lem:cusp-small-error}.  \qed

\noindent{\bf Proof of Theorem \ref{thm:cusp}.}  Let $\Upsilon$ be the intersection of the event of Proposition \ref{prop:cusp-bdg} and $\A_0$ so that $\Upsilon$ holds with overwhelming probability.  Let $z_S = z_j^i(S)$ be a characteristic. On the event $\Upsilon$ we have for $S_0 < S < T_i \wedge \tau$ that,
\beq
| f ( z_S, S) - \tilf (z_S, S) | \leq \int_0^S g (U) | f(z_U, U) - \tilf (z_U, U) | \d U + 2 N^{\eps_1} B (z_S)
\eeq
by the definition of $\A_0$ and \eqref{eqn:cusp-ito}, where
\beq
g(U) = \frac{ \left| z_U ( \del_z f ) ( z_U, U ) \right|}{2} .
\eeq
Hence, via Gronwall's inequality we obtain,
\beq
| f (z_S,S) - \tilf (z_S, s) | \leq 2 \int_0^S g(u) \exp \left[ \int_U^S g(w) \d w \right] N^{\eps_1} B(z_U) \d u + 2 N^{\eps_1} B ( z_S).
\eeq
As in the proofs of Theorems \ref{thm:bulk} and \ref{thm:edge} we have, using Lemma \ref{lem:cusp-small-error},
\beq
 \exp \left[ \int_U^S g(w) \d w \right] \leq C \frac{  \eta_U}{\eta_S},
\eeq
as well as 
\beq
g(U) \leq C \frac{ | \Re[ \tilf (z_U, U) ] | }{\eta_U}.
\eeq
Hence,
\beq
| f (z_S, S) - \tilf (z_S, S) | \leq C \frac{ N^{\eps_1}}{\eta_S} \int_0^S | \Re[ \tilf (z_U, U) ] | B (z_U) \d U + 2 N^{\eps_1} B(z_S).
\eeq
We now split into a few different cases. 
 Suppose first that $\eta_S \geq \kappa_S$. Then,
\beq
\int_0^S | \Re[ \tilf (z_U, U) ] | B (z_U) \d U \leq \frac{1}{N} \int_0^S | \Re[ \tilf (z_U, U) ] | \eta_U^{-1} \d U \leq C N^{-1} \log(N)
\eeq
and so,
\beq
| f ( z_S, S) - \tilf (z_S, S) | \leq C N^{\eps_1} B(z_S) + C \log(N) N^{\eps_1} \eta_S \leq C \log(N) N^{\eps_1} B ( z_S)
\eeq
where we used the assumption that $\eta_S \geq  \kappa_S$ in the last inequality. Now assume that $\eta_S \leq \kappa_S$.  Similar to the convention of Lemma \ref{lem:cusp-char-1} we define,
\beq
S_* = \sup \{ S_0 < U < S : |z_U| \geq 1 + \mfa \Delta_U \} .
\eeq
with the convention that the supremum of the empty set is $-\infty$.  By Lemma \ref{lem:cusp-good-char} it follows that $S_* = - \infty$.  Therefore, for all $S_0 < U < S$ we have,
\beq
| \Re[ \tilf (z_U, U) ] | \asymp \frac{1}{ \Delta_U^{1/6}} \frac{ \eta_U}{ \sqrt{ \kappa_U + \eta_U}},
\eeq
as $\kappa_U >0$ for all $ S_0 < U < S$.

Let us define $S_1 = 4 - 2s$.  Let us assume that $S_0 < S_1$. The other case is easier, requiring consideration of only one of the regions of integration below. We consider first the region $U \in (S_1, S)$. Then, for such $U$ we have $\Delta_U \leq C \Delta_S$, and so
\begin{align}
\frac{1}{ \eta_S} \int_{S_1}^{S} | \Re[ \tilf (z_U, U) ] | N B (z_U) \d U &\leq \frac{1}{\sqrt{ \kappa_S + \eta_S} \Delta_S^{1/6} } \int_{S_1}^S \frac{1}{ \sqrt{ \kappa_U + \eta_U} \sqrt{ \eta_U} } \d U \notag\\
&\leq C \frac{1}{\sqrt{ \kappa_S + \eta_S} \Delta_S^{1/6} } \int_{S_1}^S \Delta_U^{1/6} | \Re[ \tilf (z_U, U) ] | \eta_U^{-3/2} \d U \notag\\
&\leq C \frac{1}{\sqrt{ \kappa_S + \eta_S}  } \int_{S_1}^S  | \Re[ \tilf (z_U, U) ] | \eta_U^{-3/2} \d U \notag\\
&\leq C \frac{1}{ \sqrt{ \kappa_S + \eta_S} \eta_S^{1/2} }.
\end{align}
Now we consider $U \in (S_0, S_1)$. For such $U$, let $U = 4 - u$. Note $u \geq 2 s$. In particular, from \eqref{eqn:cusp-kappa} we conclude that
\beq
\Delta_U^{1/6} \sqrt{ \kappa_U} \geq c u
\eeq
which implies that $C \Delta_U \geq \kappa_U \geq c \Delta_U$.  From Lemma \ref{lem:ref-cusp} we see that for $U \in (S_0, S_1)$,
\begin{align}
c \frac{ \eta_U}{ \Delta_U^{2/3} } \leq | \Re[ \tilf (z_U, U) ] | \leq C \frac{ \eta_U}{ \Delta_U^{2/3}}.
\end{align}
Since $\tilf (z_U, U)$ is constant along characteristics we deduce from this that
\beq \label{eqn:cusp-b1}
\Delta_U^{1/6} \leq C \frac{ \eta_U^{1/4} \Delta_{S_1}^{1/6}}{\eta_{S_1}^{1/4}} \leq C \frac{ \eta_U^{1/4} \Delta_{S}^{1/6}}{\eta_{S}^{1/4}} 
\eeq
The second inequality uses that $\eta_{S_1} \geq \eta_S$ and that $\Delta_{S_1} \asymp \Delta_S$ by the choice of $S_1$. 
Hence,
\begin{align}
\frac{1}{ \eta_S} \int_{S_0}^{S_1} | \Re[ \tilf (z_U, U) ] | N B (z_U ) \d U &\leq \frac{C}{\sqrt{ \kappa_S + \eta_S} \Delta_S^{1/6} } \int_{S_0}^{S_1} \frac{1}{ \sqrt{ \kappa_U + \eta_U} \sqrt{ \eta_U} } \d U \notag\\
&\leq \frac{C}{\sqrt{ \kappa_S + \eta_S} \Delta_S^{1/6} } \int_{S_0}^{S_1} \Delta_U^{1/6} | \Re[\tilf (z_U, U) ] | \eta_U^{-3/2} \d U \notag\\
&\leq \frac{C}{ \sqrt{ \kappa_S + \eta_S} \eta_S^{1/4} }\int_{S_0}^{S} | \Re[ \tilf (z_U, U) ] | \eta_U^{-5/4} \d U \notag\\
& \leq C \frac{1}{ \sqrt{ \kappa_S + \eta_S} \eta_S^{1/2} }.
\end{align}
In the first inequality we used constancy of $\tilf$ along characteristics and Lemma \ref{lem:ref-cusp}. In the second inequality we used Lemma \ref{lem:ref-cusp} again. In the third inequality we used \eqref{eqn:cusp-b1}. 

Summarizing, we see that on the event $\Upsilon$ we have for any $S_0 < S < \tau \wedge T_i$ that
\beq
| f (z_S, S) - \tilf (z_S, S) | \leq C \log(N) N^{\eps_1} B( z_S).
\eeq
Taking $\eps_1 < \eps /100$ we see that we must have $\tau > T_i$ for every $i$. Therefore, $\tau = T$ and we conclude the theorem. \qed

\section{Center of mass evolution} \label{sec:com}

In this section we prove Proposition \ref{prop:com-evolution}.  We will not track dependence of constants on the parameter $N$, as Proposition \ref{prop:com-evolution} is a statement in soft analysis. The reader should note that we will therefore change notation and consider unitary matrices of size $n$ instead of $N$ evolving according to \eqref{eqn:unit-def}. We do this in order to emphasize the ineffectiveness of the estimates in the dimension of the system.

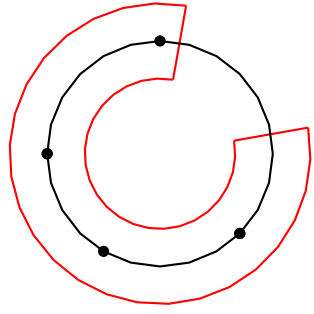
\begin{figure}

\centering
\begin{tikzpicture}
   \draw [red,thick,domain=320:1480] plot ({cos(\x/4)}, {sin(\x/4)});
   \draw [red,thick,domain=320:1480] plot ({2*cos(\x/4)}, {2*sin(\x/4)});
   \draw [red, thick, domain=100:200] plot ({\x*cos(80)/100}, {\x*sin(80)/100} );
   \draw [red, thick, domain=100:200] plot ({\x*cos(370)/100}, {\x*sin(370)/100} );
   
   \draw [black, thick, domain = 0:360] plot ( {1.5*cos(\x) }, {1.5*sin(\x)} ) ;
   
   \node at (0,1.5) [circle,fill,inner sep=1.5pt]{};
   \node at (-1.5,0) [circle,fill,inner sep=1.5pt]{};
   
   \node at (-1.5,0) [circle,fill,inner sep=1.5pt]{};
  
  \node at (1.0607,-1.0607) [circle,fill,inner sep=1.5pt]{};
  \node at (1.0607,-1.0607) [circle,fill,inner sep=1.5pt]{};
  \node at (-0.75,-1.299) [circle,fill,inner sep=1.5pt]{};
  
\end{tikzpicture}
  \captionsetup{width=.8\linewidth}
  
\caption{Schematic diagram of contour $\Gamma$ defined in \eqref{eqn:contour-def-u1}. Contour in red; unit circle in black. Eigenvalues are the thick black dots enclosed by the red contour. }
 \label{fig:contour}
\end{figure}

Fix some time $t_0 > 0$ and assume that there are no eigenvalues in the set $\{ z = \e^{ \i \theta} : \theta \in I_0\}$ for $I_0 = [ \theta_0 -L/2 , \theta_0 + L/2]$ at time $t_0$ and $\theta_0 \in [0, 2 \pi )$. Let $\tau > t_0$ be the first time an eigenvalue enters this set. Let $\Gamma$ be the contour,
\begin{align} \label{eqn:contour-def-u1}
\Gamma :=& \{ z = r \e^{ \i \theta} : 1/2 < r < 3/2, \theta = \theta_0 \pm L/4 \} \notag\\
\cup &\{ z = r \e^{ \i \theta} : r = 1/2, 3/2, \theta \notin [\theta_0 - L/4, \theta_0+L/4] \}.
\end{align}
That is, it encloses all of the eigenvalues on the unit circle but avoids the ray $\{ r \e^{ \i \theta_0 } : r \geq 0 \}$. A schematic diagram of the contour $\Gamma$ is given in Figure \ref{fig:contour}.

We use this contour as we will take a logarithm with branch cut being this ray. 
Then for any time $t_0 < t_1 < \tau$ we have,
\beq
\i \theta_i (t_1) - \i \theta_i (t_0) = \frac{1}{ 2\pi \i }  \int_{\Gamma} g_{\theta_0} ( z) \left( \frac{1}{ z - \lambda_i (t_1)} - \frac{1}{ z - \lambda_i (t_0) } \right) \d z
\eeq
where $g_{\theta_0}(z)$ is the branch of the logarithm holomorphic in $\cc \backslash \{ z = r \e^{ \i \theta_0}, r \geq 0 \}$. This formula holds due to the fact that no eigenvalue crosses the angle $\theta_0$ in this time interval, and so for this entire time interval, $\theta_i (t) = 2 \pi k + \tilde{\theta}_i (t)$ where $k$ is a constant integer and $\theta_0 < \tilde{\theta}_i (t) < \theta_0 + 2 \pi $.

Introducing,
\beq
m(z, t) = \frac{1}{n} \sum_{i=1}^n \frac{1}{ \lambda_i (t) -z } = \frac{1}{n} \tr \frac{1}{ U_t - z}
\eeq 
we therefore have,
\beq
\bar{\theta} (t_1 \wedge \tau ) - \bar{\theta} (t_0) = \frac{1}{ 2 \pi \i} \int_\Gamma \frac{ - g_{\theta_0} ( z) }{ \i } ( m (z, t_1 \wedge \tau ) - m (z, t_0) ) \d z.
\eeq
For $t_0 < t < \tau$, $m(z, t)$ is well-defined for $z \in \Gamma$ and so we may apply the It{\^{o}} formula.  Since,
\beq
f(z, t) = 1 + 2 z m(z, t)
\eeq
we have from \eqref{eqn:df} that (abbreviating $g = g_{\theta_0}$ and $m = m(z, t)$),
\begin{align} \label{eqn:com-a1}
\d \frac{1}{ 2 \pi \i } \int_{\Gamma} \i g(z) m(z, t) \d z &= \frac{1}{ 2 \pi \i } \int_{\Gamma} \frac{- \i g}{2} (1+ 2 z m)(m + zm') \d z \d t \notag\\
&+ \frac{1}{ 2 \pi \i } \int_{\Gamma} g(z) \frac{1}{n} \tr \left( \frac{ U_t}{(U_t -z )^2} \d W \right) \d z.
\end{align}

Observe that by the holomorphic functional calculus,
\beq
\frac{1}{ 2 \pi \i } \int_{\Gamma} g(z) \frac{1}{ (U_t - z)^2} \d z = U_t^{-1},
\eeq
and so
\beq
\frac{1}{ 2 \pi \i } \int_{\Gamma} g(z) \frac{1}{n} \tr \left( \frac{ U_t}{(U_t -z )^2} \d W \right) \d z = \frac{1}{n} \tr \left( \d W \right).
\eeq
The first line of \eqref{eqn:com-a1} turns out to vanish identically. To see this, we will evaluate all of the integrals explicitly using the Cauchy integral formula as well as the variants
\beq
\frac{1}{ 2 \pi \i } \int_{\Gamma} \frac{F(z)}{ (z-a)(z-b) } \d z = \frac{ F(b) - F(a)}{b-a}
\eeq
and 
\beq
\frac{1}{ 2 \pi \i } \int_{\Gamma} \frac{ F(z) }{ (z-a)^2 (z-b)} \d z = - \frac{F'(a)}{b-a} + \frac{ F(b) - F(a)}{(b-a)^2}
\eeq
valid for $F$ analytic in the appropriate domain and $\Gamma$ encircling $a, b$. We have,
\begin{align}
& \frac{1}{ 2 \pi \i} \int g(z) (1 + 2 z m )( m + z m' ) \d z = \frac{1}{ 2 \pi \i} \int g(z) (m + 2 zm^2 + z m'  + 2 z^2 m m' ) \d z \notag\\
= & \frac{1}{ 2 \pi \i} \int \d z g(z) \bigg\{  \frac{1}{n} \sum_i \frac{1}{ \lambda_i - z} + 2 \frac{z}{n^2} \sum_{i \neq j } \frac{1}{( \lambda_i - z) ( \lambda_j - z ) } + z  ( 2n^{-1} +1) \frac{1}{n} \sum_i \frac{1}{ ( \lambda_i - z )^2}  \notag\\
&+ 2 \frac{z^2}{n^2} \sum_{i \neq j } \frac{1}{ ( \lambda_i -z)( \lambda_j - z)^2} + 2\frac{z^2}{n^2} \sum_i \frac{1}{ ( \lambda_i  -z)^3} \bigg\} \notag\\
= & - \frac{1}{n} \sum_i g ( \lambda_i ) + \frac{2}{n^2} \sum_{i \neq j } \frac{ \lambda_i g ( \lambda_i ) - \lambda_j g ( \lambda_j ) }{ \lambda_i - \lambda_j } +\frac{ 2n^{-1}+1}{n} \sum_i ( g ( \lambda_i ) + 1 ) \notag\\
+ & \frac{1}{n^2} \sum_{i \neq j } \frac{ 4 \lambda_j g ( \lambda_j ) +2 \lambda_j }{ \lambda_i - \lambda_j } - \frac{2}{n^2} \sum_{i \neq j } \frac{ \lambda_i^2 g ( \lambda_i ) - \lambda_j^2 g ( \lambda_j ) }{ ( \lambda_i - \lambda_j )^2} - \frac{1}{n^2} \sum_i ( 2 g ( \lambda_i ) + 3).
\end{align}
Now, note that 
\beq
 \frac{2}{n^2} \sum_{i \neq j } \frac{ \lambda_i^2 g ( \lambda_i ) - \lambda_j^2 g ( \lambda_j ) }{ ( \lambda_i - \lambda_j )^2} = 0
\eeq
by symmetry as well as,
\beq
- \frac{1}{n} \sum_i g ( \lambda_i )+\frac{ 2n^{-1}+1}{n} \sum_i ( g ( \lambda_i ) + 1 )- \frac{1}{n^2} \sum_i ( 2 g ( \lambda_i ) + 3) =1 - \frac{1}{n}
\eeq
Finally,
\begin{align}
\frac{2}{n^2} \sum_{i \neq j } \frac{ \lambda_i g ( \lambda_i ) - \lambda_j g ( \lambda_j ) }{ \lambda_i - \lambda_j } +  \frac{1}{n^2} \sum_{i \neq j } \frac{ 4 \lambda_j g ( \lambda_j ) +2 \lambda_j }{ \lambda_i - \lambda_j } = \frac{2}{n^2} \sum_{i \neq j } \frac{ \lambda_j }{ \lambda_i - \lambda_j } = - \frac{n(n-1)}{n^2} 
\end{align}
and so indeed the term on the first line of \eqref{eqn:com-a1} vanishes. It follows that for any $t_0 < t_1 \leq \tau$ we have
\beq
\bar{\theta} (t_1) - \bar{ \theta} (t_0) = \frac{1}{n} \left( \tr W_{t_1} - \tr W_{t_0}  \right).
\eeq
Fix a final time $T >1$ and a large integer $m>0$, and intermediate times $t_i = i T / m$.  Let $\theta_k$, $k=1, 2, \dots,  4n$, be $4n$ equally spaced points along the unit circle. For any $\ell >0$ define the arcs,
\beq
Z_{k, \ell} = \{ z= \e^{ \i \theta} : | \theta - \theta_k | < \ell / ( 100 n) \}.
\eeq
Throughout the remainder of the proof we will make use of the arcs $\{ Z_{k, 1}\}_k$ and $\{ Z_{k, 1/2} \}_k$, i.e., the choices $\ell=1, 1/2$.

For each $0 \leq i \leq m$ and $1 \leq k \leq 4 n$ define the stopping time $\tau_{i, k}$ as follows. If at time $t_i$ there is an eigenvalue in $Z_{k, 1}$ let $\tau_{i, k} = t_i$. Otherwise, $\tau_{i, k}$ be the first time $t_i < t\leq t_{i+1}$ that an eigenvalue hits the set $Z_{k, 1/2}$. If this does not occur, let $\tau_{i, k} = \infty$. Then, let $\tau_i = \min_{k} \tau_{i, k}$. By the pigeonhole principle we have $\tau_i > t_i$ for all $i$. Finally, let $\tau^{(m)} = \min_i \tau_i$. 

From the above discussion we have on the event $\tau^{(m)} = \infty$ and by telescoping that for any $0 < t < T$ we have
\beq
\bar{ \theta} (t)  = \frac{1}{n} \tr \left( W_t \right).
\eeq
The proof of Proposition \ref{prop:com-evolution} will be complete once we prove that
\beq
\lim_{m \to \infty} \pp\left[ \tau^{(m)} = \infty \right] = 1.
\eeq
The remainder of this section is devoted to this.  Consider the points,
\beq
z_{\pm, k } = (1+\eta) \exp \left[ \i \left( \theta_k \pm \frac{3}{4} \frac{1}{ 100 n } \right) \right].
\eeq
Note that for any $k$ that if $Z_{k, 1}$ contains no eigenvalues at time $t_i$, then
\beq
| \Re[ f (z_{\pm, k} , t_i ) ] | \leq  C_1 n^2 \eta
\eeq
for some $C_1 >0$. On the other hand, if for this $k$ we have that $ \tau_{i, k} < \infty$ then there is some $t \in [t_i, t_{i+1}]$ such that,
\beq
| \Re[ f (z_{+, k} , t) ] |  + | \Re[ f (z_{-, k} , t) ] | > \frac{c_1}{n \eta}.
\eeq
Choosing $\eta$ small enough so that
\beq
\frac{c_1}{ n \eta} > 10 C_1 n^2 \eta
\eeq
we have that
\beq
\{ \tau_{i} < \infty \} \subseteq \bigcup_{k=1}^{4n} \left\{ \exists t \in [t_i, t_{i+1} ] : | f ( z_{+, k}, t) - f(z_{+, k}, t_i ) | + | f ( z_{-, k}, t) - f(z_{-, k}, t_i ) |> \frac{c_2}{ 10 n \eta } \right\} .
\eeq
The parameter $\eta>0$ is fixed for the remainder of the proof. 
We see by \eqref{eqn:df}  that (abbreviating $ z= z_{\pm, k}$)
\beq
| f ( z, t) - f(z, t_i ) | \leq CT \frac{1}{m \eta^3} + (M_t - M_{t_i} ).
\eeq
By the BDG inequality,
\beq
\pp\left[ \sup_{t_i < t < t_{i+1} } | M_{t_i} - M_t | > s \right] \leq C \exp\left[ - c s m^{1/2} n^{-2} T^{-1} \right].
\eeq
So taking $s = c_2 / (100 n \eta )$ we see that for all $m$ large enough, there are constants depending on $n$ and $\eta$ such that
\beq
\pp\left[ \tau_i < \infty \right] \leq C \e^{- m^{1/2} c}.
\eeq
This completes the proof of Proposition \ref{prop:com-evolution}. \qed

\section{Rigidity} \label{sec:rigidity}

\subsection{Helffer-Sj{\"o}strand formula} \label{sec:HS}

This section establishes an analog of the Helffer-Sj{\"o}strand formula for measures on the unit circle.  Recall,
\beq
\del_{\bar{z}} = \frac{1}{2} ( \del_x + \i \del_y)
\eeq
as well as Green's theorem,
\beq
F( \lambda ) = \frac{1}{ \pi} \int_{\rr^2}\frac{ (\del_{\bar{z}} F) (x, y) }{ \lambda - (x+ \i y) } \d x \d y
\eeq
for any $F \in C^2$ of compact support.  In polar coordinates we recall,
\beq
\del_{\bar{z}} F (r, \theta)  = \frac{ \e^{ \i \theta}}{2} \left( \del_r F + \frac{ \i}{r} \del_\theta F \right).
\eeq
For any $F$ supported in an annulus,
\beq
F( \lambda ) = \frac{1}{ \pi} \int_{\rr^2} \frac{ \del_{\bar{z}} F }{ 2z} \frac{ \lambda+z}{ \lambda- z} \d x \d y.
\eeq
Let $\varphi : [0, 2\pi ] \to \rr$  be a function on the unit circle that extends to a smooth $2\pi$-periodic function on $\rr$. We define the quasi-analytic extension of $\varphi$ by,
\beq
\tilphi (r, \theta ) := ( \varphi ( \theta) - \i \log (r) \varphi' ( \theta ) ) \chi (r)
\eeq
where $\chi (r)$ is a function that is $1$ on $[3/4,4/3]$ and $0$ outside of $[1/2, 2]$.  We may also assume that,
\beq
\chi(r) = \chi (r^{-1} ).
\eeq
If $\mu$ is any measure on the unit circle with Cauchy transform $f_\mu (z)$ we see that,
\begin{align} \label{eqn:HS}
\int \varphi ( \theta) \d \mu ( \theta) =& \frac{1}{ \pi}  \int_{\rr^2} \frac{ \del_{\bar{z}} \tilphi }{2 z } f_\mu (z) \d x \d y \notag\\
= & \frac{ 1}{ 4 \pi } \int ( \varphi ( \theta) - \i \log (r) \varphi' ( \theta ) ) \chi' (r) f_\mu (z) \d r \d \theta \notag\\
+ & \frac{1}{4 \pi  } \int \varphi'' ( \theta ) \log (r)  \chi (r) f_\mu (z) r^{-1} \d r \d \theta.
\end{align}
This is the Helffer-Sj{\"o}strand formula we require to establish our eigenvalue estimates and will be used in the next subsections. 

\subsection{Proof of Corollary \ref{cor:interval}} \label{sec:interval}

We follow Section 3.3 of \cite{huang2019rigidity} very closely, as the estimates of Theorem \ref{thm:bulk} and the formula \eqref{eqn:HS} combine in the exactly same manner there as here to prove Corollary \ref{cor:interval}.  Introduce,
\beq
\eta ( \theta ) = \inf_{ \eta > 1 } \{ (\eta-1) | \Re[ \tilf ( \eta \e^{ \i \theta } ) ] | \geq N^{\eps-1} \} .
\eeq
Note that from Lemma \ref{lem:calc-2}, 
\beq
\eta \to (\eta-1) | \Re[ \tilf ( \eta \e^{ \i \theta } ] |
\eeq
is an increasing function for $\eta >1$. 
We first assume $t < 2$. We consider $I = [ \theta_0, \pi]$ for some $ 0 < \theta_0 < \pi$. The case $- \pi < \theta_0 < 0$ requires only notational changes. Define,
\beq
\tileta := \inf_{ \eta : \log \eta \geq N^{1-\mfc} } \{ \eta : \max_{ \theta_0 \leq x \leq \theta_0 + \eta } \eta (x) \leq \eta \}.
\eeq
Let $\tiltheta$ be
\beq
\tiltheta := \mbox{argmax}_{ \theta_0 - \tileta \leq x \leq \theta_0 } \eta (x)
\eeq
so that
\beq
\eta ( \tiltheta) = \tileta .
\eeq
We let $\varphi$ be a function that is $1$ on $[\theta_0, \pi]$ and $0$ outside of $[ \theta_0 - \log \tileta \wedge \frac{1}{10}, \pi+N^{2\eps-1} ]$, with $| \varphi^{(k)} (x) | \leq C_k (N^{1-2\eps} )^k$ for $x$ near $\pi$ and $| \varphi^{(k)} (x) | \leq C_k (\log \tileta )^{-k}$ for $x$ near $\theta_0$.

By \eqref{eqn:HS} we have,
\begin{align}
\left| \frac{1}{N} \sum_i \varphi ( \lambda_i (t) ) - \int \varphi (\theta) \rho_t ( \theta) \d \theta \right| &\leq C \int ( | \varphi( \theta) | + | \varphi' ( \theta ) | ) | \chi' ( r) | S (z)| \d r \d \theta \notag\\
&+ \left| \int \varphi'' ( \theta) \log (r) \chi (r) S (z) r^{-1} \d r \d \theta \right|.
\end{align}
where
\beq
S (z) = f (z, t) - \tilf (z, t).
\eeq
On the event that the estimates of Theorem \ref{thm:bulk} hold we see that,
\beq
\int ( | \varphi( \theta) | + | \varphi' ( \theta ) | ) | \chi' ( r) | S (z)| \d r \d \theta  \leq C \frac{N^{\eps}}{N}.
\eeq
For the second term, note that the measure $r^{-1} \d r$ is invariant under the transformation $r \to r^{-1}$ so that,
\beq
\left| \int \varphi'' ( \theta) \log (r) \chi (r) S (z) r^{-1} \d r \d \theta \right|= 2 \left| \int_{0}^{2 \pi} \int_{r>1} \varphi''( \theta) \log (r) | \Re[S(z) ] | \chi(r) r^{-1} \d r \d \theta \right|
\eeq
Recall that $\varphi''(x) = 0$ unless $x \in [\pi, \pi+N^{2 \eps-1}]$ or $x \in [ \theta_0 -  \log \tileta \wedge \frac{1}{10}, \theta_0  ]$. Note that for $r \geq \eta_\mfc$ where $\log \eta_\mfa := N^{-\mfc}$ we have,
\beq \label{eqn:int-a1}
| \Re[ S (z)] | \leq C \frac{ N^{\eps}}{N \log |z|}
\eeq
using Lemma \ref{lem:calc-2} (see Remark 3.3 of \cite{huang2019rigidity} for a similar argument).  Therefore,
\begin{align}
& \int_{\theta \in [\pi, \pi+N^{2 \eps-1}] } \int_{r>1} | \varphi''( \theta) | \log (r) | \Re[S(z) ] | \chi(r) r^{-1} \d r \d \theta \notag\\
&\leq \int_{\theta \in [\pi, \pi+N^{2 \eps-1}] } \int_{r>\eta_\mfc} | \varphi''( \theta) | \log (r) | \Re[S(z) ] | \chi(r) r^{-1} \d r \d \theta \notag\\
&+ C \int_{\theta \in [\pi, \pi+N^{2 \eps-1}] } \int_{1< r< \eta_\mfc} | \varphi''( \theta) | \chi(r)  \d r \d \theta \notag\\
\leq & C \frac{N^{3 \eps}}{N} + \frac{C}{ N^\mfc} \leq C \frac{N^{3 \eps}}{N}.
\end{align}
We used \eqref{eqn:int-a1} in the first integral and the estimate $|(r-1) \Re[ S ( r \e^{ \i \theta} ) ] | \leq C$ for the second.  For the region where $\theta \in [ \theta_0 -  \log \tileta \wedge \frac{1}{10}, \theta_0  ] =: J$ we first bound,
\begin{align}
\int_{ \theta \in J} \int_{1 < r < \tileta} | \log(r) | | \Re[S(z) ] \varphi'' ( \theta) \chi (r) | \d r \d \theta & \leq \int_{ \theta \in J} \int_{1 < r < \eta_\mfc} | \log(r) | | \Re[S(z) ] \varphi'' ( \theta) \chi (r) | \d r \d \theta \notag\\
&+ \int_{ \theta \in J} \int_{\eta_\mfc < r < \tileta} | \log(r) | | \Re[S(z) ] \varphi'' ( \theta) \chi (r) | \d r \d \theta \notag\\
&\leq C \frac{ \eta_\mfc}{\tileta} + \frac{ C N^{\eps}}{N} \leq C N^{\eps-1}.
\end{align}
For the first integral we used $|(r-1) \Re[S] | \leq 2$ and that $\tileta \geq N^{1-\mfc}$. For the second region we used again \eqref{eqn:int-a1}. For the contribution of $r > \tileta$ we have, by partial integration
\begin{align}
\int_{ \theta \in J} \int_{ r> \tileta} \log(r) \varphi''(\theta) \chi (r) \Re[S] r^{-1} \d r \d \theta &= \int_{ \theta \in J} \theta' ( \theta) \log ( \tileta) \Im [ S]  \tileta^{-1} \d \theta \notag\\
&- \int_{ \theta \in J } \int_{r > \tileta} \varphi' ( \theta) \del_r ( \log(r) \chi (r) r^{-1} )\Im[S] \d r \d \theta.
\end{align}
By definition of $\tileta$, all $z$ appearing in the above integration lie in $\B_t$. Therefore, we may apply the estimate on $S$ of Theorem \ref{thm:bulk} and obtain that both of these integrals are bounded above by $C \log(N) N^{\eps-1}$. 

When these estimates hold we therefore conclude that,
\beq
\left| \left\{ i : \theta_i (t) \in I \right\} \right| \leq N \int \varphi ( \theta) \rho_t ( \theta ) \d \theta + C N^{3 \eps}.
\eeq
For $t < 2$, $\rho_t$ has no support for $| \theta - \pi | < \delta$ some $\delta >0$. Therefore,
\beq \label{eqn:int-a2}
N \int \varphi ( \theta) \rho_t ( \theta ) \d \theta \leq N \int_I \rho_t ( \theta) \d \theta + N \int_{\theta_0 - \log \tileta \wedge \frac{1}{10} }^{\theta_0} \rho_t ( \theta) \d \theta.
\eeq
There are two cases. If $\tileta = N^{-\mfc}$ then the second integral is bounded by $\tileta / \sqrt{t} = N^{-\mfc} t^{-1/2}$. Otherwise,
\beq
\int_{\theta_0 - \log \tileta \wedge \frac{1}{10} }^{\theta_0} \rho_t ( \theta) \d \theta \leq C (\tileta-1) | \Re[ \tileta \e^{ \i \tiltheta} ] | \leq C \frac{N^{ 3 \eps}}{N}.
\eeq
The lower bound for the number of $\theta_i \in I$ follows similarly.  For $t >2$, then the densities $\rho_t$ all satisfy that $\inf_{|\theta| < c} \rho_t (\theta) > c$ for some $c>0$. In this case, one uses intervals with one end-point at $\theta =0$. In \eqref{eqn:int-a2} there is then a second term on the RHS,
\beq
\int_{-N^{2\eps-1}}^0 \rho_t (\theta) \d \theta \leq C N^{2 \eps-1}.
\eeq
Everything else is identical. This proves the corollary. \qed

\subsection{Proof of Corollary \ref{cor:bulk-rig}} \label{sec:bulk-rig}

It follows from Proposition \ref{prop:com-evolution} that with overwhelming probability,
\beq
\sup_{0 \leq t \leq T} | \bar{\theta} (t) | \leq \frac{N^{\eps}}{N}
\eeq
for any $\eps >0$. Next, for each $i$, we may write,
\beq
\theta_i (t) = 2 \pi n_i (t) + \varphi_i (t)
\eeq
for $n_i(t)$ an integer and $- \pi \leq \varphi_i (t) \leq \pi$. For any fixed $t >0$ it follows from the fact that the eigenvalues never cross that the set $\{ n_1 (t), n_2 (t), \dots, n_N (t) \}$ contains at most two consecutive integers whose absolute values we will denote by $m (t), m(t) + 1$.  Let $N_1(t)$ be the number of eigenvalues s.t. $|n_i (t)| = m(t)$ and $N_2 (t)$ be the number of eigenvalues s.t. $| n_i (t) | = m(t) +1$. Let $\mu$ be a probability measure on $[- \pi, \pi]$ with a density. We have,
\beq
\int_0^\pi  \theta \d \mu ( \theta ) = \int_0^\pi \mu( [t, \pi] ) \d t,
\eeq
and since $\int_{-\pi}^{\pi} \theta \rho_t ( \theta) \d \theta = 0$ we have,
\beq
\left| \frac{1}{N} \sum_i \varphi_i (t) \right| \leq \frac{ N^{\eps}}{N}
\eeq
by Corollary \ref{cor:interval}. From this it follows that,
\beq
(1+m(t) ) N_2 (t) + m(t) N_1 (t) \leq N^{\eps}
\eeq
with overwhelming probability. In particular, $m(t) = 0$ and the number of eigenvalues s.t. $| \theta_i (t) | > \pi$ is at most $N^{\eps}$ with overwhelming probability. Using this, the remainder of Corollary \ref{cor:bulk-rig} follows from Corollary \ref{cor:interval} in a straightforward manner similar to the proof of Corollary 3.2 of \cite{huang2019rigidity}. \qed

\appendix

\section{Density of states calculations} \label{a:dos}

Recall the limiting spectral measure $\rho_t ( \theta)$ and its Cauchy transform,
\beq
\tilf (z, t) = \int_{0}^{2\pi} \frac{\e^{ \i \theta} + z}{ \e^{\i \theta} - z } \rho_t ( \theta) \d \theta.
\eeq
Note that,
\beq
\tilf (z, 0) = \frac{1+z}{1-z},
\eeq
and that $\tilf$ is constant along characteristics,
\beq
t \to z \exp \left[ \frac{t}{2} \tilf (z, 0) \right] = z \exp \left[ \frac{t}{2} \frac{1+z}{1-z} \right].
\eeq
Moreover, define the region,
\beq
\Gamma_t := \left\{ z \in \cc: \Re[z] > 0 \mbox{ and } \left| \frac{z-1}{z+1} \e^{\frac{t}{2} z } \right| <1 \right\}.
\eeq
Then, the function $\tilf (z, t)$ is a conformal map of the open unit disc into $\Gamma_t$ \cite{biane1997segal}. The domain $\Gamma_t$ has the following form. We let $x_+ (t)>1 $ be the smallest solution larger than $1$ of
\beq
\frac{ 1-x}{x+1} \e^{\frac{t}{2} x} = 1,
\eeq
and if $t >4$ we let, $x_-(t)$ be the largest solution less than $1$ of
\beq
\frac{ x-1}{x+1} \e^{\frac{t}{2} x} = 1.
\eeq
If $t \leq 4$ set $x_- (t) = 0$. For $x \in (x_- (t), x_+(t) )$ set
\beq
k_t (x) := \sqrt{ \frac{ (x+1)^2 - (x-1)^2 \e^{tx} }{ \e^{tx  }-1  }}.
\eeq
If $t < 4$ the boundary of $\Gamma_t$ is the union of the vertical line on the imaginary axis,
\beq
\{ z = x + \i y : x = 0, |y| \leq \sqrt{4t^{-1} -1 } \}
\eeq
and the curves
\beq
\{ z = x \pm \i k_t (x) : 0 < x \leq x_+ (t) \} .
\eeq
When $t\geq 4$ the boundary is the union of the curves,
\beq
\{ z = x \pm \i k_t (x) :  x_- (t) \leq x \leq x_+ (t) \}.
\eeq
Moreover, $\tilf(z, t)$ extends to a bijection of the closed unit disc to the closure of $\Gamma_t$. It is also a conformal map of the complement of the closed unit disc in $\cc$ to $- \Gamma_t$.

For the moment we consider the case $ t< 4$. We have that $\tilf(z, t)$ sends the arc
\beq
\{ z = \e^{ \i \theta } : \pi \geq | \theta | > \Theta_t \}
\eeq
to the boundary of $\Gamma_t$ that intersects the imaginary axis. Here,
\beq \label{eqn:theta-t-a-def-u1}
\Theta_t := \frac{1}{2} \sqrt{ (4-t)t} + \arccos \left( 1 - \frac{t}{2} \right) = \frac{1}{2} \sqrt{ (4-t)t} +2 \arcsin \left( \sqrt{\frac{t}{4}}\right) 
\eeq
Using the expansion,
\beq
\arcsin(1-x) = \frac{\pi}{2} - \sqrt{2} x^{1/2} - \frac{ \sqrt{2}}{12} x^{3/2} + \O ( x^{5/2} )
\eeq
one can check that the gap satisfies,
\beq
\Delta_t = 2 ( \pi - \Theta_t) = \frac{s^{3/2}}{3} + \O ( s^{5/2} )
\eeq
where $t = 4-s$. We require the following a-priori bound. First,
note that for any Cauchy transform of a measure on the unit circle and $ z= r \e^{ \i \varphi}$,
\begin{align}
\int_0^{2 \pi } \frac{\e^{ \i \theta } +z }{ \e^{\i \theta} - z} \rho ( \theta) \d \theta &= (1- |z|^2 ) \int_0^{2 \pi} \frac{1}{ | \e^{ \i \theta } - z |^2} \rho ( \theta) \d \theta + \i 2 r \int_0^{2 \pi} \frac{ \sin ( \varphi - \theta) }{ | \e^{ \i \theta } - z |^2} \rho ( \theta) \d \theta.
\end{align}
\bep \label{prop:a-apriori}
The following estimate holds for all $|z| <1$ and $t>0$,
\beq
|\tilf (z, t) | \leq 3 \left(1 + t^{-1/2} \right).
\eeq
\eep
\proof We have,
\begin{align}
\tilf (z, t) = \tilf ( z \e^{ -t/2 \tilf (z, t) }, 0) 
\end{align}
and so with $w = z \e^{ -t/2 \tilf (z, t) }$ and $|z| <1 $,
\begin{align}
\left| \tilf (z, t) \right|^2 \leq  (1+|w|)^2 \int_{0}^{2 \pi} \frac{1}{ | \e^{ \i \theta} - w|^2} \rho_0 ( \theta) \d \theta \leq  \frac{(1+|w|)^2}{ 1 - |w|^2} \Re[ \tilf (w, 0) ] = \frac{1+|w|}{ 1 - |w| } \Re[ \tilf (z, t) ]
\end{align}
Note since $\Re[ \tilf (z, t) ] >0$ we have $|w| \leq 1$. If $t\Re [ \tilf (z, t)] \geq 2$ then we see $1- |w| \geq 1 -\e^{-1}$ and so, e.g.,
\beq
|\tilf (z, t) | \leq 3.
\eeq
If $t \Re[ \tilf (z, t) ] \leq 2$ then the inequality $\e^{-x} \leq 1 - x + x^2/2$, valid for $x >0$ implies that $1- \e^{-x} \geq \frac{x}{2}$ for $0 <x < 1$. Therefore,
\beq
\frac{1}{1 - |w| } \leq \frac{4}{t \Re[ \tilf (z, t) ] },
\eeq
and so
\beq
| \tilf (z, t) |^2 \leq \frac{8}{t} .
\eeq
This yields the claim. \qed

The extension of $\tilf$ to the closed unit disc satisfies,
\beq \label{eqn:a-self-1}
\frac{ \tilf (z, t) -1 }{ \tilf (z, t) + 1} \e^{ \frac{t}{2} \tilf (z,t ) } = z .
\eeq
Let,
\beq
w_0 = \i \sqrt{ \frac{4}{t} -1 }, \qquad z_0 = \e^{ \i \Theta_t}.
\eeq
Note that this is the value of $\tilf (z_0, t) = w_0$.  Our goal is now to derive an approximate cubic equation via Taylor expansion for $\tilf (z, t)$ where $z$ are points on the boundary of the unit disc close to $z_0$. This will allow for conclusions about the behavior of $\rho_t$ near the edge later.  In order to facilitate this expansion,  we introduce the coordinates $E >0$ and $ q \in \cc$ via,
\beq
z = \e^{ \i ( \Theta_t - E) }, \qquad \tilf (z,t) = w_0 -\i q
\eeq
so that $\Im[q]  \geq 0$, as $\Re[ \tilf (z, t) ] > 0$ on the open unit disc, with strict inequality extending to points on the boundary inside the support of the spectral measure.

 With these coordinates,
\beq
\rho_t (E) = \frac{1}{ 2 \pi} \Im[q ((1^-)\e^{\i E } ) ]
\eeq
in the sense of boundary values. 
We now implement a  Taylor expansion of the self-consistent equation for $\tilf (z, t)$ to obtain the following.
\bep
Let $0 < t < 4$. Let $t =4 -s$.   Let $\tilf( z + z_0, t) = w_0 - \i q$. The following holds for $|q| \leq  t^{-1/2}/10  + \1_{ \{ t > 1 \} } 100$. First,
\begin{align}
z = \i q^2 \e^{ \i \Theta_t } \frac{ t^{3/2} s^{1/2}}{8} + \i q^3 \e^{ \i \Theta_t} \frac{(3-t) t^2}{24} + \O (t^{5/2} |q|^4 )
\end{align}
Second, if $z = \e^{ \i (\Theta_t - E) }$ then $\tilf (z, t) = w_0 - \i q$ satisfies,
\beq \label{eqn:cubic-edge}
0 = E - \i \frac{E^2}{2} + q^2 \frac{t^{3/2} s^{1/2}}{8} + q^3 \frac{ (3-t) t^2}{24} + \O ( t^{5/2} |q|^4 + |E|^3 )
\eeq
\eep
\remark Due to Proposition \ref{prop:a-apriori}, the assumption on $q$ holds for all $z$ for $1 < t < 4$. 

\proof 
Consider,
\beq
F(w) = \frac{w-1}{w+1} \e^{ \frac{t}{2} w }.
\eeq
The equation \eqref{eqn:a-self-1} can be written as,
\beq
F( \tilf (z, t) ) = z.
\eeq
Then, via Taylor expansion,
\begin{align}
F(w+w_0) - F(w_0) = -w^2 \e^{ \i \Theta_t} \i \frac{ t^{3/2} s^{1/2}}{8} + w^3 \e^{ \i \Theta_t} \frac{ (3-t) t^2}{24} + \O ( t^{5/2} |w|^4 ).
\end{align}
The restriction on the size of $|q|$ in the assumptions of the proposition allow us to obtain the explicit $t^{5/2}$ behavior of the error term for short times $t$. 
Note also that the appearance of $w$ in the denominator of $F$ causes no trouble in the above expansion as $\Re[ \tilf (z, t) ] \geq 0$.
The claim follows. \qed

Using the above expansion, we can derive the square-root behavior for times away from $0$ and $4$.

\bep
Let $\delta >0$ and assume $\delta < t < 4-\delta$. Then, for $|E| \leq c_\delta $ we have,
\beq
\rho_t (\Theta_t - E) = \frac{E^{1/2}}{ \pi} \sqrt{ \frac{2}{ t^{3/2} (4-t)^{1/2} } } \left( 1+ \O (E) \right).
\eeq
\eep
\proof From \eqref{eqn:cubic-edge} we see that for $|E| \leq c $ that $|q| \leq C E^{1/2}$. Writing $q = x + \i y$ one can first see that $y \asymp E^{1/2}$ and then consequently that $|x| \leq C E$. Therefore, taking real parts of \eqref{eqn:cubic-edge} we find,
\beq
y^2 \frac{ t^{3/2} s^{1/2}}{8} = E + \O ( |E|^2)
\eeq
and so the claim follows. \qed

At $t=4$, the following shows that we have an exact cusp. The proof will also serve as a useful warm-up for the near-cusp case.

\bep
Let $t=4$.  Then,
\beq
\rho_t (\pi - E ) = E^{1/3} \left( \frac{3}{2} \right)^{1/3} \frac{ \sqrt{3}}{4 \pi} \left( 1 + \O (|E|^{1/3} ) \right)
\eeq
\eep
\proof We have from \eqref{eqn:cubic-edge} that,
\beq
0 = E - \i \frac{E^2}{2} - \frac{2}{3} q^3 + \O ( |q|^4 + |E|^3 ).
\eeq
We see that $|q| \leq C E^{1/3}$. We rewrite this as,
\beq
q^3 = \zeta 
\eeq
where
\beq
\zeta =  \frac{3}{2} E ( 1 + \O ( |E|^{1/3} ) ).
\eeq
Note that $\Re[\zeta] >0$ for $0 < E < \delta$ some $\delta >0$. If $\zeta = r \e^{ \i \theta}$, then $|\theta | \leq C |E|^{1/3}$ and $r \asymp E$, and the three roots of the cubic are the functions,
\beq
r_+ = r^{1/3} \e^{2 \pi \i /3 + \i \theta/3}, \qquad r_0 = r^{1/3} \e^{ \i \theta/3}, \qquad r_- = r^{1/3}  \e^{-2 \pi \i /3 + \i \theta/3}.
\eeq
These are continuous and distinct for $E >0$. Since $q$ is a continuous function of $E$ we have for $|E| \leq \delta$ that $q$ equals one of these three functions. It cannot be $r_-$ since $\Im[r_-] <0$. The description of the boundary of $\Gamma_t$ shows that,
\beq
\lim_{E \dto 0 } \left| \frac{ \Im[q]}{\Re[q]} \right| = \sqrt{3}
\eeq
and so it cannot be $r_0$. This yields the claim. \qed

We will now consider the behavior of $\rho_t$ for times $t$ close to $4$. For this, we will follow relatively closely the arguments of Section 9 of the work \cite{qve} of Ajanki, Erd{\H{o}}s and Kr{\"u}ger. There, they showed how one may use Cardano's formula for roots of cubics to find expressions for spectral measures. As the coefficients of our equation are explicit, we are able to short-cut some of the arguments of \cite{qve} in adapting them to our setting.

Following \cite{qve}, we introduce the universal edge shape function,
\beq
\Psi_e ( \lambda )= \frac{ \sqrt{ (1+ \lambda) \lambda ) }}{ \left( 1 + 2 \lambda + 2 \sqrt{ (1+ \lambda) \lambda ) } \right)^{2/3} + \left( 1 + 2 \lambda - 2 \sqrt{ ( 1 + \lambda) \lambda } \right)^{2/3} + 1}
\eeq
As advertised, we characterize the spectral measure in the case of two nearby edges in terms of this shape function.
\bep \label{prop:cusp-shape-u1}
Let $3.5 < t < 4$. Then,
\beq
\rho_t (\Theta_t - E ) =  \Delta_t^{1/3} \left( \Delta_t^{1/6} \frac{ (t-3) t^{5/4} }{2 \sqrt{2}  s^{1/4} \pi} \right) \Psi_e \left( \frac{E}{ \Delta_t } \right) + \O \left( \min\{ E^{3/2}s^{-5/4}, E^{2/3} \} \right),
\eeq
or equivalently,
\beq
\rho_t (\Theta_t - E ) = s^{1/2} \left( \frac{ t^{5/4}(t-3)}{\pi 2 \sqrt{6} } \right) \Psi_e \left( \frac{3 E}{ s^{3/2}} \right) + \O \left( \min\{ E^{3/2}s^{-5/4}, E^{2/3} \} \right)
\eeq
\eep
\proof First, from \eqref{eqn:cubic-edge} we see that,
\beq \label{eqn:a-qbd-a1}
|q| \leq C \begin{cases} E^{1/2} s^{-1/4}, & E \leq s^{3/2} \\ E^{1/3} , & E \geq s^{3/2} \end{cases}
\eeq
for $|E| \leq \delta$, some $\delta >0$. Let,
\beq
a = \frac{ t -3}{24} t^2, \qquad b = \frac{ t^{3/2} s^{1/2}}{8}
\eeq
so that the cubic becomes
\beq
a q^3 - b q^2 - E = \O (|q|^4 + |E|^2 ). 
\eeq
Define,
\beq
3 a q = b( Q + 1) .
\eeq
This is our version of the ``normal coordinates'' introduced at the start of Section 9.2.2 of \cite{qve}. 
The equation for the cubic becomes,
\beq
Q^3 - 3 Q + (-2 - 27 a^2 b^{-3} E ) = b^{-3} \O ( |q|^4 + |E|^2 ).
\eeq
We rewrite this as,
\beq
Q^3 - 3 Q + 2  \zeta = 0.
\eeq
where 
\beq
\zeta = -1 - \frac{27}{2} \frac{a^2}{b^3} E (1 +  \mu (E) )
\eeq
and $\mu(E)$ is a function defined implicitly. Note that for $|E| \leq \delta$, the function $\mu(E)$ satisfies,
\beq
| \mu(E) | \leq C \begin{cases} E s^{-1}, & E \leq s^{3/2} \\ E^{1/3}, & E \geq s^{3/2} \end{cases} ,
\eeq
where we used \eqref{eqn:a-qbd-a1} to bound the $|q|^4$ term. 
We can take $E \leq \delta$ to ensure that $| \mu(E)| \leq \frac{1}{10}$.  So for $0 < E < \delta $ we have $\Re[\zeta] < -1$.  From Lemma 9.13 of \cite{qve} we have that the three roots of the cubic,
\beq
X^3 - 3 X + 2 \zeta =0
\eeq 
are $X_+ (\zeta)$, $X_- (\zeta)$ and $X_0 ( \zeta)$ where,
\beq
X_0 = - ( \Phi_+ + \Phi_- ) , \qquad X_\pm = \frac{1}{2} ( \Phi_+ + \Phi_- ) \pm \i \frac{ \sqrt{3}}{{2}} ( \Phi_+ - \Phi_- )
\eeq
and
\beq
\Phi_\pm ( \zeta ) = - \left( - \zeta \mp \sqrt{ \zeta^2-1} \right)^{1/3}
\eeq
for $\Re[ \zeta] \leq -1$. The three functions are continuous and distinct. It follows that $Q(E)$ is equal to one of these three root functions for all $|E| \leq \delta$. We must determine which one.  Here, we follow the ``Choice of $a_1$'' section of the proof of Lemma 9.14 of \cite{qve}.

First, note that for $E$ close to $0$ and so for $\zeta$ close to $-1$, we have that $X_0 = 2 + o(1)$. But by definition, for $E$ close to $0$, we see that $q = o (1)$ and so $Q  = -1 + o (1)$. Therefore, the required root cannot be $X_0$. 

One can has that $X_- (\zeta)$ is H{\"o}lder-$\frac{1}{2}$ on $\{ z \in \cc : \Re[z] < -1 \}$ and that,
\beq
\Im [ X_- (-1-x) ] \leq - c x^{1/2}
\eeq
for small $x >0$ (see Lemma 9.15 and (9.119) of \cite{qve}). Hence, for $\zeta = \zeta(E)$ as above, we have
\beq
\Im [ X_- ( \zeta) ] \leq - c E^{1/2} b^{-3/2}+ C |E \mu(E)|^{1/2} b^{-3/2}
\eeq
which is negative for sufficiently small $E$. This shows that,
\beq
Q = X_+ ( \zeta (E) )
\eeq
for $0 < E < \delta$, with $\delta >0$ as above.

In the remainder of the proof we follow some of the arguments in the proof of Proposition 9.8 at the end of Section 9.2.2 of \cite{qve} to replace the argument of $X_+$ by a real-valued quantity. From (9.111) of \cite{qve} we see that,
\begin{align}
\left| X_+ ( \zeta(E) ) - X_+ \left( -1 - \frac{27}{2} a^2 b^{-3} E \right) \right| \leq C \min\{ |E b^{-3}|^{1/2}, |E b^{-3} |^{1/3} \} | \mu (E) |.
\end{align}
Now,
\beq
\Im X_+ (-1 -2 \lambda) = 2 \sqrt{3} \Psi_e ( \lambda).
\eeq
This yields,
\beq
\rho_t ( \Theta_t -E ) = \frac{ b}{\pi \sqrt{3}} \Psi_e \left( \frac{27}{4} \frac{a^2}{b^3} E \right) + \O \left( \min\{ E^{3/2}s^{-5/4}, E^{2/3} \} \right).
\eeq
Note that,
\beq
\hat{\Delta}_t := \frac{ 4 b^3}{27 a^2} = \frac{ s^{3/2}}{3} + \O (s^{5/2} ) = \Delta_t + \O ( s^{5/2} ).
\eeq
Therefore, using (9.145) of \cite{qve} we see that,
\beq
\frac{ b}{\pi \sqrt{3}} \Psi_e \left( \frac{27}{4} \frac{a^2}{b^3} E \right) = \frac{ b}{\pi \sqrt{3}} \left( \frac{\Delta_t}{ \hat{\Delta}_t } \right)^{1/2} \Psi_e \left( \frac{E}{ \Delta_t} \right)+ \O \left( \min\{ E^{3/2}s^{-5/4}, E^{2/3} \} \right).
\eeq
We now have,
\beq
\frac{ b}{\pi \sqrt{3}} \left( \frac{\Delta_t}{ \hat{\Delta}_t } \right)^{1/2} = \Delta_t^{1/3} \left( \Delta_t^{1/6} \frac{ (t-3)t^{5/4} }{2 \sqrt{2}   s^{1/4} \pi} \right) \asymp \Delta_t^{1/3}
\eeq
This yields the first estimate.  The second is proven similarly. \qed

The asymptotics of the shape function are given in (9.63) of \cite{qve},
\beq \label{eqn:cusp-shape-u1}
\Delta_t^{1/3} \Psi_e \left( \frac{E}{ \Delta_t } \right) \asymp \begin{cases} \frac{E^{1/2}}{ \Delta_t^{1/6}}  , & E \leq \Delta_t \\ E^{1/3} , & E \geq \Delta_t \end{cases} .
\eeq

The short time regime is as follows.

\bep \label{prop:short-time-shape-u1}
There is a $\delta>0$ so that for all $t < \frac{1}{10}$ we have,
\beq
\rho_t ( \Theta_t - E) = \frac{ \sqrt{2}}{ \pi s^{1/4} t^{1/2}} \sqrt{ \frac{E}{ t^{1/2} } }\left( 1 + \O (E t^{-1/2} ) \right)
\eeq
for $|E| \leq t^{1/2} \delta$. 
\eep
\proof Let $\tilf (z, t) = w_0 - \i t^{-1/2} q$ and $z = \e^{ \i ( \Theta_t - t^{1/2} E)}$ so that from \eqref{eqn:cubic-edge} we have,
\beq
E - \i t^{1/2} \frac{E^2}{2} + q^2 \frac{s^{1/2}}{8} + q^3 \frac{(3-t)}{24} = \O \left( |q|^4 + t |E|^3 \right).
\eeq
We see that for $|E| \leq \delta$ we have $|q| \leq C |E|^{1/2}$. Letting $q = x + \i y$ we first see that $y \asymp E^{1/2}$ and then that $|x| \leq CE$. We find,
\beq
y^2 \frac{s^{1/2}}{8} = E + \O (|E|^2)
\eeq
and the claim follows. \qed

We consider now $ t> 4$. Let $w_0 = x_- (t) :=x$.  This is the value of $\tilf (z, t)$ at $ z = -1$. Note that it is purely real. The equation defining $w_0$ is
\beq
\e^{tx/2} = \frac{2}{ 1 - x} -1.
\eeq
Introducing $t = 4 + s$ we can expand this to find,
\beq \label{eqn:a-a1}
x^2 (4- t^3/24) + - x 2s(1+s/8) - s = \O (x^3).
\eeq
We see that $|x| \leq C s^{1/2}$, and then that,
\beq
x = s^{1/2} ( 4 - t^3/24)^{-1/2} ( 1 + \O ( s^{1/2} ) ) = \sqrt{3} s^{1/2}/2 + \O (s).
\eeq
Similar to the regime $ t< 4$ we derive a self-consistent equation for the Cauchy transform, but now in the regime $ t >4$. 
\bep  There is a $\delta >0$ so that for $ 4 < t < 4 + \delta$ the following holds.
Let $z = \e^{ \i ( \pi - E)}$  and $\tilf (z, t) = w_0  - \i q$. Then,
\begin{align} \label{eqn:a-cubic-2}
E + A q + B \i q^2 + C q^3 = \O ( |E|^2 + |q|^4 )
\end{align}
where
\beq
A = \frac{ 4 w_0^2 + s w_0^2 - s}{ 2 ( 1 + w_0) (1- w_0 )}
\eeq
and
\beq
B =  \left( w_0 (2t - t^2/8) - (t-4)^2/8+ w_0^2 (t+t^2/8) + w_0^3t^2/8 \right)  \frac{1}{ (1+w_0)^2 (w_0 -1 ) }
\eeq
and
\beq
C = \frac{1}{ (1+w_0)(w_0-1)}\left( \frac{t^3}{48} ( w_0^2-1) + \frac{t^2}{4} - \frac{t}{ 1+w_0} + \frac{2}{ (1+w_0)^2} \right)
\eeq
We have $A >0$ and $C<0$ and $B <0$, and $|A| \asymp s$ and $|B| \asymp s^{1/2}$ and $|C| \asymp 1$. 
\eep
\proof Letting $F(w) = \frac{w-1}{w+1} \e^{\frac{t}{2} w}$ as before, we have
\begin{align}
F(w+w_0) - F (w_0) &= w \frac{ \e^{t w_0/2}}{ 2 (1+w_0)^2} ( 4 w_0^2 + s w_0^2 - s) \notag\\
&+ w^2 \frac{ \e^{t w_0/2}}{ (1+w_0)^3 } \left( w_0 (2t - t^2/8) - (t-4)^2/8+ w_0^2 (t+t^2/8) + w_0^3t^2/8 \right) \notag\\
&+ w^3 \frac{ \e^{t w_0/2}}{ (1+w_0)^2 } \left( \frac{t^3}{48} ( w_0^2-1) + \frac{t^2}{4} - \frac{t}{ 1+w_0} + \frac{2}{ (1+w_0)^2} \right) + \O ( |w|^4).
\end{align}
This yields the equation. The asymptotics for the coefficients  are as follows. The asymptotics for $C$ are easy. For $A$, we just read off that
\beq
4 w_0^2 - s \asymp 4^3 w_0^2/24
\eeq
from the equation \eqref{eqn:a-a1}. For $B$, we have that $w_0 \asymp s^{1/2}$ which is much larger than the other terms contributing to $B$. 
 \qed
 
We now find an expression for $\rho_t$ near $\theta = \pi$ for $t$ just larger than $4$. In the proof below, we follow closely  Section 9.1 of \cite{qve}.  We define,
\beq
\Psi_m ( \lambda ) := \frac{ \sqrt{1+ \lambda^2}}{ ( \sqrt{ 1 + \lambda^2} + \lambda )^{2/3} + ( \sqrt{ 1 + \lambda^2} - \lambda )^{2/3} -1 } -1.
\eeq
 \bep
 There is a $\delta >0$ so that for $4 < t < 4 + \delta $ we have for $|E| \leq \delta $ that,
 \beq
 \rho_t (\pi - E) - \rho_t ( \pi ) = \frac{s^{1/2}}{ 4 \pi}(1+J )  \Psi_m \left(\frac{6E}{s^{3/2}} \right) + \O \left( \min \{ |E| s^{-1/2}, |E|^{2/3} \} \right)
 \eeq
 and
 \beq
 \rho_t ( \pi ) = \frac{\sqrt{3} s^{1/2}}{ 4 \pi } ( 1 + \O (s^{1/2} )),
 \eeq
 and above
 \beq
 J = 6 \gamma_1+ 3 \gamma_1^2 - 2 \gamma_2 = \O ( s^{1/2} )
 \eeq
 where $\gamma_i$ are introduced in the proof below. 
 \eep
\proof Dividing the equation \eqref{eqn:a-cubic-2} by $C$ we obtain a new cubic,
 \beq
 q^3 + \i b q^2 - a q +E/C = \O ( |E|^2 + |q|^4).
 \eeq
 We note that,
 \beq
 b = 3 \sqrt{3} s^{1/2}/2 + \O (s), \qquad a= 3s/2 + \O ( s^{3/2} ).
 \eeq
 Introducing $\alpha = \sqrt{3} s^{1/2} /2$ we have,
 \beq \label{eqn:a-cubic-3}
 q^3 + \i 3 \alpha (1 + \gamma_1 ) q^2 - 2 \alpha^2 (1+ \gamma_2 ) q + E/C = \O ( |E|^2 + |q|^4 )
 \eeq
 where $\gamma_i = \O ( s^{1/2} )$. Explicitly,
 \beq
 \gamma_1 = \frac{B}{3 C \alpha} -1, \qquad \gamma_2 = -\frac{A}{2 C\alpha^2} -1
 \eeq 
 The equation \eqref{eqn:a-cubic-3}  is now essentially in the same form as (9.34) of \cite{qve}. Working first with the coordinates,
 \beq
 Z = \alpha^{-1} q, \qquad \lambda = E/\alpha^3
 \eeq
 this equation becomes,
 \beq
 Z^3 + 3 \i (1+ \gamma_1 ) Z^2 - 2 (1+\gamma_2 ) Z + \lambda/C = \O ( \alpha^3 | \lambda|^2 + \alpha |Z|^4 ).
 \eeq
 For sufficiently small $\lambda$ we see that $|Z| \leq C | \lambda|$.  Then, comparing real and imaginary parts we see that with $Z= x + \i y$ that $|x| \asymp \lambda$ and then that $y \asymp \lambda^2$.  For $| \lambda| \geq c$ we then see that,
 \beq
 |Z| \leq C | \lambda|^{1/3}.
 \eeq
 Hence,
 \beq
 \Im[q] \leq C \begin{cases} \alpha^{-5} |E|^2, & |E| \leq \alpha^3 \\ |E|^{1/3}, & |E| \geq \alpha^3 \end{cases}
 \eeq
 and
 \beq
 |\Re[q] | \leq C \begin{cases} \alpha^{-2} |E| & |E| \leq \alpha^3 \\ |E|^{1/3} , & |E| \geq \alpha^3 \end{cases}.
 \eeq
 We now introduce normal coordinates similar to (9.35) of \cite{qve}. Letting now,
 \beq
 q = X Q + Y
 \eeq
 where $Y = -  \i \alpha (1+ \gamma_1 )$ and $\sqrt{3} X =\alpha (1 + 6 \gamma_1 +3 \gamma_1^2 - 2 \gamma_2 ) = \alpha (1+ \gamma_3) $ for some $\gamma_3 = \O ( s^{1/2})$ we obtain the cubic equation,
 \beq
 Q^3 + 3 Q + 2 \zeta = 0
 \eeq
 where
 \beq
 \zeta = 3 \sqrt{3} \frac{- E}{2 |C| \alpha^3(1+ \gamma_3)^3} + \gamma_8 + \alpha^{-3} \O ( |E|^2 + |q|^4)
 \eeq
 for some $\gamma_8 = \O (s^{1/2} )$, that is independent of $E$. 
The three roots of this cubic equation are given by the three functions $X_0 ( \zeta)$, and $X_\pm ( \zeta)$ where,
\beq
X_0 = - 2 \Phi_{o} \qquad X_\pm = \Phi_o \pm \i \sqrt{3} \Phi_e
\eeq
where $\Phi_o$ and $\Phi_e$ denote the odd and even parts of the function $\Phi = \left( \sqrt{1 + \zeta^2 } + \zeta \right)^{1/3}$. These functions coincide only at $\zeta = \pm \i$. Since $\zeta \neq \pm \i$ and the root functions are continuous, $Q $ equals one of the three branches for all $\zeta$. Now when $E=0$ we know that $q = 0$ and so $Q = \i \sqrt{3}$. By the argument near (9.49) of \cite{qve} we conclude that
\beq
Q( \zeta ) = X_+ ( \zeta)
\eeq
for $|E| \leq  \delta_1$ some $\delta_1 >0$ and for all sufficiently small $s>0$. We conclude similarly to the proof of Proposition 9.3 of \cite{qve} (see especially (9.52) and (9.57) appearing there) using 
\beq
\Psi_m ( \lambda) = \frac{ \Im [ X_+ ( \lambda) - X_+ (0) ] }{ \sqrt{3}}.
\eeq 
This completes the proof. \qed

The asymptotics of $\Psi_m$ are,
\beq
s^{1/2} \Psi_m \left( \frac{6E}{s^{3/2}} \right) \asymp \begin{cases} \frac{|E|^2}{s^{5/2}}, & |E| \leq s^{3/2} \\ |E|^{2/3} , & |E| \geq s^{3/2} \end{cases} .
\eeq

From the characterization of $\Gamma_t$ we see that the function $\rho_t ( \theta)$ is monotonic in $[0, \pi]$. Hence, we conclude the following.
\bel
Let $\delta >0$. Then there is a $c >0$ so that if either $t > 4 + \delta$ and $\theta \in [0, \pi]$ or $ t > \delta$ and $0 \leq \theta \leq \Theta_t - \delta$ (where we set $\Theta_t = \pi$ for $ t \geq 4$) we have,
\beq
\rho_t ( \theta)  \geq c.
\eeq
If $ t< \delta$ and $  0 \leq \theta < \Theta_t - \sqrt{t} \delta_1$ some $\delta_1 >0$ then there is a $c_1$ depending on $\delta$ and $\delta_1$ such that,
\beq
\rho_t ( \theta) \geq c_1 t^{-1/2}.
\eeq
\eel

\section{Auxilliary results} \label{a:calc}

\bel \label{lem:calc}
For $x >1$ we have,
\beq
 x \log(x) \leq \frac{ x^2-1}{2}
\eeq
\eel
\proof Let $f(x) = x \log (x)$. Then $f'(x) = \log(x) +1$ and $f''(x) = x^{-1}$. By Taylor's theorem,
\beq
f(x) -f(1) = (x-1) f'(1) + \frac{(x-1)^2}{2} f''(s)
\eeq
for some $1 < s < x$. Clearly $f''(s) \leq 1$ and the claim follows. \qed

\bel \label{lem:calc-2}
Let $f_\mu (z)$ be a Cauchy transform of a probability measure $\mu$ on the unit circle. Then the function,
\beq
r \to (r-1) | \Re[ f_\mu (r \e^{ \i \theta } ) ] |
\eeq
is increasing for $r>1$.
\eel
\proof It suffices to consider the case $\theta =0$ and $\mu$ is a point mass at $\e^{ \i \theta'}$. In this case we need to show that,
\beq
r \to (r-1) \frac{r^2-1}{1 + r^2 - 2 r x}
\eeq
is increasing where $x = \cos ( \theta')$. The derivative with respect to $r$ of the RHS is,
\beq \label{eqn:calc-a1}
\frac{1}{ (1 + r^2 - 2 r x )^2}\left\{ (3r^2-2r - 1) (1 + r^2 - 2 r x ) - (r^3 - r^2 - r + 1)(2r - 2 x) \right\}
\eeq
On the other hand,
\beq
\del_x ((3r^2-2r - 1) (1 + r^2 - 2 r x ) - (r^3 - r^2 - r + 1)(2r - 2 x) ) = - 4 r^3 + 2 r^2 +2 < 0
\eeq
for $r >1$. Therefore, the quantity in the parentheses of \eqref{eqn:calc-a1} is minimized when $x=1$, over the domain of possible values of $x \in [-1, 1]$.  But when $x=1$, this quantity is nothing more than,
\beq
 (1+r^2 - 2r)^2 \times \del_r \frac{ (r-1)(r^2-1)}{1+r^2- 2 r} .
\eeq
But
\beq
\frac{ (r-1)(r^2-1)}{1+r^2- 2 r}  = r+1
\eeq
which is obviously increasing.  Therefore, we conclude that the quantity in the parentheses of \eqref{eqn:calc-a1} is positive. This yields the claim. \qed

\noindent{\bf Proof of Lemma \ref{lem:char-path}.} We first observe,
\beq
\Im\left[ \frac{\e^{\i \tilde{\theta}} + \e^{\i \theta}}{\e^{\i \tilde{\theta}} - \e^{\i \theta}} +  \frac{e^{-\i \tilde{\theta}} + \e^{\i \theta}}{\e^{- \i \tilde{\theta}}- \e^{\i \theta}}\right] = \cot \left( \frac{\tilde{\theta} + \theta}{2} \right) + \cot \left( \frac{-\tilde{\theta} + \theta}{2} \right),
\eeq
so that
\beq \label{eqn:calc-b1}
\Im \left[  f ( \e^{ \i \theta} ) \right] = \int_0^{\pi} \left( \cot \left( \frac{\tilde{\theta} + \theta}{2} \right) + \cot \left( \frac{-\tilde{\theta} + \theta}{2} \right) \right) \rho ( \tilde{\theta} ) \d \tilde{\theta} ,
\eeq
by the assumed symmetry of $\rho ( \tilde{\theta} )$. 
If we consider $\pi \geq \theta > E$ then the integrand is non-zero only for  $\theta > \tilde{\theta}$ and $\tilde{\theta}< \pi$. Therefore the arguments of the cotangent functions are in $(0, \pi)$ where it is well-defined. Furthermore, since $\theta < \pi$ we have that $\frac{\theta - \tilde{\theta}}{2}< \pi - \frac{\tilde{\theta} + \theta}{2} $, so that
\beq
 \cot \left( \frac{\tilde{\theta} + \theta}{2} \right) + \cot \left( \frac{-\tilde{\theta} + \theta}{2} \right) = \cot \left( \frac{-\tilde{\theta} + \theta}{2} \right) -\cot \left(\pi - \frac{\tilde{\theta} + \theta}{2} \right) >0
\eeq
because $\cot(x)$ is strictly decreasing for $ x\in (0, \pi)$.  
  Thus, the integrand in \eqref{eqn:calc-b1} is strictly positive. Clearly it is $0$ for $\theta = \pi$.
  
  Now,
  \beq
 2 \frac{ \d}{ \d \theta} \left( \cot \left( \frac{\tilde{\theta} + \theta}{2} \right) + \cot \left( \frac{-\tilde{\theta} + \theta}{2} \right) \right)= -\csc^2\left( \frac{\theta + \tilde{\theta}}{2}\right) - \csc^2\left( \frac{-\tilde{\theta}+ \theta}{2}\right) < 0
  \eeq
  so that $\Im [ f ( \e^{ \i \theta} )]$ is strictly decreasing for $E < \theta < \pi$.  This proves the second part of the lemma.

  Now, we consider more general $z = r \e^{ \i \theta}$. First, 
\begin{equation}
    \Im \left[ \frac{\e^{\i \tilde{\theta}} + r \e^{\i \theta}}{\e^{\i \tilde{\theta}} -r \e^{\i \theta}}\right] = \frac{2r \sin(\theta - \tilde{\theta})}{1+ r^2 -2r \cos(\theta - \tilde{\theta})}.
\end{equation}
Calculating the derivative in $r$, we obtain
\begin{equation} 
    \partial_r\frac{2r \sin(\theta - \tilde{\theta})}{1+ r^2 -2r \cos(\theta - \tilde{\theta})} = \frac{2 \sin(\theta - \tilde{\theta})(1-r^2)}{(1 + r^2 -2 r \cos(\theta - \tilde{\theta}))^2}.
\end{equation}
Therefore,
\begin{align} \label{eqn:calc-b2}
\frac{ \del_r \Im [ f ( r \e^{ \i \theta } ) ]}{1-r^2} = \int_0^{\pi} \bigg\{ \frac{2 \sin(\theta - \tilde{\theta})}{(1 + r^2 -2 r \cos(\theta - \tilde{\theta}))^2} + \frac{2 \sin(\theta + \tilde{\theta})}{(1 + r^2 -2 r \cos(\theta + \tilde{\theta}))^2} \bigg\} \rho ( \tilde{\theta} ) \d \tilde{ \theta} 
\end{align}
Assume now that $E < \theta < \pi - \eps$. We claim that there is a $c_\eps >0$ so that the quantity on the RHS is strictly positive for all  $0 < r-1 < c_\eps$.  Assume first that $E < \frac{ \pi}{8}$. Then for $E < \theta < \frac{ \pi}{4}$ we see that every term in the parentheses on the RHS of \eqref{eqn:calc-b2} is positive. So we assume that $E + \eps < \theta < \pi - \eps$. For such $\theta$ we have,
\beq \label{eqn:calc-b3}
\lim_{r \to 1} \bigg\{ \frac{2 \sin(\theta - \tilde{\theta})}{(1 + r^2 -2 r \cos(\theta - \tilde{\theta}))^2} + \frac{2 \sin(\theta + \tilde{\theta})}{(1 + r^2 -2 r \cos(\theta + \tilde{\theta}))^2} \bigg\} = \frac{\cos(\frac{\theta - \tilde{\theta}}{2})}{4 \sin^3(\frac{\theta - \tilde{\theta}}{2})} +\frac{\cos(\frac{\theta + \tilde{\theta}}{2})}{4 \sin^3(\frac{\theta + \tilde{\theta}}{2})}
\eeq
Since we are assuming $E +\eps < \theta < \pi - \eps$ we have,
\beq \label{eqn:calc-b4}
\frac{\eps}{2} < \frac{ \theta - \tilde{\theta}}{2} < \frac{\pi}{2} - \frac{\eps}{2}, \qquad \frac{\eps}{2} < \frac{ \theta + \tilde{\theta}}{2} < \pi - \eps, \qquad \frac{ \theta - \tilde{\theta}}{2} < \pi - \frac{ \theta + \tilde{\theta}}{2} - \eps
\eeq 
Therefore, the cosine terms on the LHS of \eqref{eqn:calc-b3} are bounded away from $\pm 1$ and so the convergence is uniform over such $\tilde{\theta}$ and $\theta$. Moreover, the last inequality of \eqref{eqn:calc-b4} also implies that the quantity on the RHS of \eqref{eqn:calc-b3} is bounded below by a positive constant depending on $\eps >0$. Therefore, we see that there is a $c_\eps >0$ so that if $0 <r-1 < c_\eps$ the quantity on the LHS of \eqref{eqn:calc-b2} is positive. This completes the proof in the case that $E < \pi/8$.

We now assume that $\rho (\theta) \leq M$ and that $E > \frac{\pi}{8}$.  Fix now $\delta_1 >0$ so small that $\delta_1 (M+1) < \frac{1}{10}$.  Assume also $\delta_1 < \frac{ \eps}{10}$. For $E \leq \theta < \pi - \eps$ and $0 \leq \tilde{\theta} \leq E - \delta_1$, the estimates \eqref{eqn:calc-b4} all hold with $\eps$ replaced by $\min\{ \eps, \delta_1 \}$.  It follows from the above argument that there is a small $c_1 >0$ so that,
\beq
\int_0^{E-\delta_1} \bigg\{ \frac{2 \sin(\theta - \tilde{\theta})}{(1 + r^2 -2 r \cos(\theta - \tilde{\theta}))^2} + \frac{2 \sin(\theta + \tilde{\theta})}{(1 + r^2 -2 r \cos(\theta + \tilde{\theta}))^2} \bigg\} \rho ( \tilde{\theta} ) \d \tilde{ \theta} \geq 2 c_1 \int_{0}^{E-\delta_1} \rho (\tilde{\theta} ) \d \tilde{ \theta} \geq c_1
\eeq
by the choice of $\delta_1 >0$ and since $\rho$ is a probability measure (i.e., by the choice of $\delta_1$, the measure $\rho$ assigns at least mass $1/2$ to the interval $[-E+\delta_1, E-\delta_1]$), for all $0 < r-1 < c_1$. Now, $\delta_1 $ and $c_1$ are fixed for the remainder of the proof.

Fix now a $\delta_2 >0$ so small that $\delta_2 < \frac{\eps}{10}$ and $\delta_2 < \frac{ \pi}{16}$. We assume also $\delta_2 < \delta_1$.  Then for $E- \delta_2 < \tilde{\theta} \leq E$ and $E \leq \pi - \eps$ we see that
\beq
\frac{\pi}{8} \leq  \tilde{\theta} + \theta \leq 2 \pi - \eps, \qquad 0 \leq \theta -  \tilde{\theta} \leq \pi .
\eeq 
For $\tilde{\theta}$ and $\theta$ obeying the first of the above inequalities,
\beq
\sup_{1 < r < 5} \left| \frac{2 \sin(\theta + \tilde{\theta})}{(1 + r^2 -2 r \cos(\theta + \tilde{\theta}))^2} \right| \leq C_1
\eeq
for some $C_1>0$. So taking $\delta_2  < \frac{c_1}{10 (C_1+1)(M +1) }$ we see that,
\begin{align}
& \int_{E-\delta_2}^{E} \bigg\{ \frac{2 \sin(\theta - \tilde{\theta})}{(1 + r^2 -2 r \cos(\theta - \tilde{\theta}))^2} + \frac{2 \sin(\theta + \tilde{\theta})}{(1 + r^2 -2 r \cos(\theta + \tilde{\theta}))^2} \bigg\} \rho ( \tilde{\theta} ) \d \tilde{ \theta} \notag\\
 \geq & \int_{E-\delta_2}^{E}  \frac{2 \sin(\theta + \tilde{\theta})}{(1 + r^2 -2 r \cos(\theta + \tilde{\theta}))^2}  \rho ( \tilde{\theta} ) \d \tilde{ \theta} \geq - \frac{c_1}{2},
\end{align}
for all $1 < r < 5$ and $E \leq \theta \leq \pi - \eps$. On the other hand, for $E - \delta_1 \leq \tilde{\theta} \leq E - \delta_2$ and $E \leq \theta \leq E - \eps$, the inequalities \eqref{eqn:calc-b4} hold with $\eps$ replaced by $\min\{ \eps, \delta_2\}$.  It follows from the above arguments that there is a $c_2 >0$ so that for $1 < r < 1+c_2$ we have,
\beq
\int_{E-\delta_1}^{E-\delta_2} \bigg\{ \frac{2 \sin(\theta - \tilde{\theta})}{(1 + r^2 -2 r \cos(\theta - \tilde{\theta}))^2} + \frac{2 \sin(\theta + \tilde{\theta})}{(1 + r^2 -2 r \cos(\theta + \tilde{\theta}))^2} \bigg\} \rho ( \tilde{\theta} ) \d \tilde{ \theta} >0.
\eeq
We complete the proof by taking $c_\eps  = \min \{ c_1, c_2, 3\}$. \qed

\bel \label{lem:calc-3}
Let $\rho_t (E)$ be as above and $T <4$. There is an $\eps >0$ so that with $ z= r \e^{ \i ( \Theta_t + \kappa ) }$ with,
\beq
0 \leq r-1 < \eps, \qquad 0 < \kappa < \eps
\eeq
we have,
\beq
\Im [ \tilf ( z, t )  ] - \Im [ \tilf (\e^{ \i \Theta_t } , t)] \leq - d \sqrt{ \kappa}.
\eeq
for some $d>0$, uniformly in $0 < t < T$.
\eel
\proof By Lemma \ref{lem:char-path} it suffices to consider $r=1$, by taking $\eps$ much smaller than $\pi - \Theta_T$.  For $z = \e^{ \i ( \Theta_t + \kappa)}$ we have by direct calculation.
\beq
  \Im [ \tilf (\e^{ \i \Theta_t } , t)]  -\Im [ \tilf ( z, t )  ] =   \int_{0}^{2 \Theta_t} \frac{\rho_t(\Theta_t-\theta)}{2} \left\{  \cot \left( \frac{\theta}{2}\right)- \cot\left(\frac{\kappa+\theta}{2} \right)\right\} \d \theta.
\eeq
With our choice of $\eps >0$ we see that the quantity in the parentheses on the RHS is positive. Therefore,
\beq
  \Im [ \tilf (\e^{ \i \Theta_t } , t)]  -\Im [ \tilf ( z, t )  ] \geq  \int_{0}^{2 \Theta_t \wedge \kappa} \frac{\rho_t(\Theta_t-\theta)}{2} \left\{  \cot \left( \frac{\theta}{2}\right)- \cot\left(\frac{\kappa+\theta}{2} \right)\right\} \d \theta.
\eeq
For $0 < \theta < \kappa$ we have by the mean value theorem,
\beq
\cot \left( \frac{\theta}{2}\right)- \cot\left(\frac{\kappa+\theta}{2} \right) \geq \frac{c}{ \kappa}
\eeq
for some $c>0$.  First, for $t \geq \frac{1}{10}$ we may assume $\kappa < \Theta_t$ by decreasing $\eps >0$ and furthermore that $\rho_t ( \Theta_t-\theta ) \geq c \sqrt{ \theta}$. Therefore,
\beq
 \int_{0}^{2 \Theta_t \wedge \kappa} \frac{\rho_t(\Theta_t-\theta)}{2} \left\{  \cot \left( \frac{\theta}{2}\right)- \cot\left(\frac{\kappa+\theta}{2} \right)\right\} \d \theta \geq c \int_0^\kappa  \frac{ \sqrt{ \theta}}{\kappa} \d \theta \geq c \sqrt{\kappa} 
\eeq
as desired. For $t < \frac{1}{10}$ we have that $\Theta_t \asymp \sqrt{t}$ (directly from \eqref{eqn:theta-t-a-def-u1}) and from Proposition \ref{prop:short-time-shape-u1},
\beq
\rho(\Theta_t - \theta ) \geq c \sqrt{\theta}t^{-3/4}, \qquad \theta < c \sqrt{t},
\eeq
for some $c>0$. 
Therefore,
\begin{align}
\int_{0}^{2 \Theta_t \wedge \kappa} \frac{\rho_t(\Theta_t-\theta)}{2} \left\{  \cot \left( \frac{\theta}{2}\right)- \cot\left(\frac{\kappa+\theta}{2} \right)\right\} \d \theta \geq \frac{c}{t^{3/4}} \int_0^{\sqrt{t} \wedge \kappa} \frac{ \sqrt{ \theta}}{\kappa} \d \theta \geq \frac{ c ( \sqrt{t} \wedge \kappa )^{3/2}}{t^{3/4} \kappa}.
\end{align}
Finally,
\beq
\frac{  ( \sqrt{t} \wedge \kappa )^{3/2}}{t^{3/4} \kappa} \geq c \sqrt{ \kappa}
\eeq
for a small $c>0$. This yields the claim. \qed

\subsection{Cusp calculation}

In this section we consider general densities on the unit circle $\rho(x) \d x$ with support $[-\pi + \Delta/2, \pi - \Delta/2]$ where $0 < \Delta < \frac{1}{100}$. We assume that $\rho$ is symmetric, $\rho(x) = \rho (-x)$.  Let $\Theta = \pi - \Delta/2$ denote the edge. We moreover assume that there is a $C_1 >0$ so that
\beq
\frac{1}{C_1} \sqrt{x} \leq \Delta^{1/6} \rho (\Theta - x ) \leq C_1 \sqrt{x}, \qquad 0 < x < \Delta
\eeq
and 
\beq
\frac{1}{C_1} x^{1/3} \leq \rho (\Theta - x )  \leq C_1 x^{1/3}, \Delta < x < \Theta.
\eeq
\bep \label{prop:a-cusp}
Under the above assumptions there is a $\frac{1}{2} > c_1 >0$ and a $d>0$ depending only on $C_1>0$ so that the following holds.  For $ z= (1+\eta) \e^{ \i ( \Theta+ \kappa )}$ with $0 < \eta <  c_1 \Delta$ and $ 0 < \kappa < \Delta /2$,
\beq
\Im [ f (z) ] - \Im [ f ( \e^{ \i \Theta} ) ] \leq - d \Delta^{-1/6} \sqrt{ \kappa}.
\eeq
\eep
\proof By direct calculation with $z = r \e^{ \i \theta}$,
\begin{align}
- \Im[ f(z)] + \Im [f ( \e^{ \i \Theta} ) ] &= \int_{-\pi}^{\pi} \frac{ -2 r \sin ( \theta - x) }{ (1- r)^2+ 2 r \sin^2 \left( \frac{ x - \theta}{2} \right) } \rho (x) \d x \notag\\
&+ \int_{-\pi}^{\pi} \frac{ 2  \sin ( \Theta - x) }{+ 2  \sin^2 \left( \frac{ x - \Theta}{2} \right) } \rho (x) \d x
\end{align}
We will first find a lower bound for the contributions of the integrals for $ x \in [0, \pi]$, which equals,
\begin{align}
&\int_{0}^{\pi} \frac{ -2 r \sin ( \theta - x) }{ (1- r)^2+ 2 r \sin^2 \left( \frac{ x - \theta}{2} \right) } \rho (x) \d x + \int_{0}^{\pi} \frac{ 2  \sin ( \Theta - x) }{+ 2  \sin^2 \left( \frac{ x - \Theta}{2} \right) } \rho (x) \d x \notag\\
= & \int_{0}^{\Theta} \rho(\Theta - x) \left( \frac{ \sin(x) }{ 4 \sin^2 (x/2) } - \frac{ (1+\eta) \sin (x+\kappa)}{ \eta^2 + 4 (1+\eta) \sin^2 \left( \frac{x+\kappa}{2} \right) } \right) \d x
\end{align}
Consider first the case $ \eta < \kappa < \Delta/2$. Note that $\sin ( x + \kappa ) >0$ for $x$ in the integrand above, and that for any $y  \in \rr$,
\beq
\del_\eta \frac{ 1+ \eta }{ \eta^2 + (1+\eta) y } = \frac{ - 2 \eta - \eta^2}{ \left( \eta^2 + (1+\eta) y  \right)^2} < 0
\eeq
for $\eta >0$. Therefore, a lower bound is achieved when $\eta = 0$ and so,
\begin{align}
&\int_{0}^{\Theta} \rho(\Theta - x) \left( \frac{ \sin(x) }{ 4 \sin^2 (x/2) } - \frac{ (1+\eta) \sin (x+\kappa)}{ \eta^2 + 4 (1+\eta) \sin^2 \left( \frac{x+\kappa}{2} \right) } \right) \d x \notag\\
 \geq &  \frac{1}{4} \int_0^{\Theta} \rho(\Theta - x) \left( \frac{ \sin(x)}{  \sin^2 (x/2) } - \frac{ \sin ( x+\kappa) }{ \sin^2 \left( \frac{ x + \kappa}{2} \right) } \right) \d x \notag\\
 = & \frac{1}{2} \int_0^{\Theta} \rho ( \Theta - x) \left( \cot \left(\frac{x}{2} \right) - \cot \left( \frac{ x + \kappa}{2} \right) \right) \d x
\end{align}
For $\frac{\kappa}{2} < x < \kappa$ we have that,
\beq
\left( \cot \left(\frac{x}{2} \right) - \cot \left( \frac{ x + \kappa}{2} \right) \right) \geq \frac{c}{ \kappa}
\eeq
for a small $c >0$ and so,
\begin{align}
\int_0^{\Theta} \rho ( \Theta - x) \left( \cot \left(\frac{x}{2} \right) - \cot \left( \frac{ x + \kappa}{2} \right) \right) \d x \geq \int_{\kappa/2}^{\kappa}2 \Delta^{-1/6} \sqrt{x} \kappa^{-1} \d x \geq c \Delta^{-1/6} ( \kappa + \eta)^{1/2}
\end{align}
for some small $c>0$.  We consider now the case $ 0 < \kappa < \eta < \Delta/2$. Note that the function,
\beq
a \to \frac{ \sin (x+a)}{ \sin^2 \left( \frac{x+a}{2} \right) } = 2 \cot \left( \frac{x+a}{2} \right)
\eeq
is decreasing for $a \in (0, \Delta/2)$ and $x \in (0, \Theta)$. Therefore,
\begin{align}
 & \int_{0}^{\Theta} \rho(\Theta - x) \left( \frac{ \sin(x) }{ 4 \sin^2 (x/2) } - \frac{ (1+\eta) \sin (x+\kappa)}{ \eta^2 + 4 (1+\eta) \sin^2 \left( \frac{x+\kappa}{2} \right) } \right) \d x \notag\\
 \geq &  \int_0^{\Theta} \rho ( \Theta - x ) \sin(x+\kappa) \left( \frac{1}{ 4 \sin^2 \left( \frac{x+\kappa}{2} \right) } - \frac{ r }{ \eta^2 + 4 r \sin^2 \left( \frac{x+\kappa}{2} \right) } \right) \d x \notag\\
 = & \int_0^{\Theta} \rho ( \Theta - x ) \sin (x + \kappa )  \left( \frac{ \eta^2}{ 4 \sin^2 \left( \frac{x+\kappa}{2} \right) \left( \eta^2 + 4 r \sin^2 \left( \frac{x+\kappa}{2} \right) \right) } \right) \d x \notag\\
 \geq & c \Delta^{-1/6}  \int_\eta^{2 \eta } \sqrt{x} \frac{ \eta^2}{ (x+\kappa)( \eta^2 + (x+\kappa)^2 } \d x \geq c  \Delta^{-1/6} \eta^{1/2} \geq c \Delta^{-1/6} ( \eta + \kappa)^{1/2}.
\end{align}
We have therefore proven that,
\beq
\int_0^{\pi} \left( \Im \left[ \frac{ \e^{ \i  x } +z}{ \e^{ \i x } - z } \right] - \Im \left[ \frac{ \e^{ \i  x } + \e^{ \i \Theta}}{ \e^{ \i x } - \e^{ \i \Theta} } \right] \right) \rho (x) \d x \geq c_1 \Delta^{-1/6} ( \kappa + \eta)^{1/2}
\eeq
for some small $c_1 >0$. On the other hand, using that,
\beq
\left| \frac{\e^{ \i x } +z}{ \e^{ \i x} - z} -  \frac{\e^{ \i x } +w}{ \e^{ \i x} - w}\right| \leq C |z-w|\left( \frac{1}{ | \e^{ \i x } - z |^2} + \frac{1}{ | \e^{ \i x } - w |^2} \right)
\eeq
for $z, w$ in bounded regions of the complex plane, one easily sees that for $|z| \geq 1$ in the upper-half plane,
\begin{align}
\left| \int_0^{-\pi} \left( \Im \left[ \frac{ \e^{ \i  x } +z}{ \e^{ \i x } - z } \right] - \Im \left[ \frac{ \e^{ \i  x } + \e^{ \i \Theta}}{ \e^{ \i x } - \e^{ \i \Theta} } \right] \right) \rho (x) \d x \right| \leq C |z - \Theta| \int_0^\Theta \frac{ \rho (\Theta- x) }{ x^2 + \Delta^2} \d x.
\end{align}
A routine calculation shows that the integral is $\O ( \Delta^{-2/3})$. It follows that there is a small $c_2 >0$ so that for $|z - \Theta| \leq c_2 \Delta$ we have,
\beq
\Im [ f ( \e^{ \i \Theta} ) ] - \Im [ f (z) ] \geq \Delta^{-1/6} \frac{c_1}{2} \sqrt{ \kappa + \eta }.
\eeq
From the proof of Lemma \ref{lem:char-path} we see that, since $\rho (x)$ is symmetric, that
\beq
\theta \to \Im [ f ( \e^{ \i \theta } ) ]
\eeq
is decreasing for $ \Theta < \theta < \pi$. Therefore, there is a $c_3 >0$ so that for $c_2 \Delta /2 \leq \kappa \leq \pi- \Theta$ we have,
\beq
\Im [ f ( \e^{ \i \Theta } ) ] - \Im [ f( \e^{ \i ( \Theta + \kappa ) } ] \geq c_3 \Delta^{1/3}.
\eeq
On the other hand, for $c_2 \Delta /2 \leq \kappa \leq \pi - \Theta$ we see from the above argument that,
\beq
\left| \Im [ f ( \e^{ \i ( \Theta + \kappa ) } ) ] - \Im [ f ( (1+\eta) \e^{ \i ( \Theta + \kappa ) } ) ] \right| \leq C \eta \Delta^{-2/3}
\eeq
It follows that there is a $c_4 >0$ so that for $\eta \leq c_4 \Delta$ and $c_2 \Delta /2 \leq \kappa \leq \pi - \Theta$ that,
\beq
 \Im [ f ( \e^{ \i \Theta } ) ]  -\Im [ f ( (1+\eta) \e^{ \i ( \Theta + \kappa ) } ) ] \geq c_4 \Delta^{1/3}.
\eeq
This completes the proof. \qed

\subsection{Extended spectral domain}

The following proposition is used in the Proof of Theorem \ref{thm:cusp}. It uses only as inputs Theorems \ref{thm:bulk} and \ref{thm:edge}. 

\bep \label{prop:a-extended}
Let $c < t < 4-c$ for some $c>0$. Let $\eps, \delta >0$. In the domain
\beq
\D_1 := \{ z  = (1+\eta) \e^{ \i \theta} : N^{\delta-1} < \eta < \delta^{-1} \}
\eeq
the estimate
\beq \label{eqn:ext-ll-1}
| \tilf (z, t) - f (z, t) | \leq \frac{N^{\eps}}{N \eta}
\eeq
holds with overwhelming probability. In the domain,
\beq
\D_2 := \{ z = (1+\eta) \e^{ \i \theta} : N^{\delta-1} < \eta < \delta^{-1} , \Theta_t + N^{-2/3+\delta} < |\theta | \leq \pi \}
\eeq
the estimate
\beq \label{eqn:ext-ll-2}
| \tilf (z, t) - f (z, t) | \leq \frac{N^{\eps}}{N ( \eta + \kappa ) }
\eeq
holds with overwhelming probability.
\eep
\proof Let $\lambda_i = \e^{ \i \theta_i}$ be the eigenvalues and $\gamma_i$ be the quantiles of $\rho_t (\theta)$. We know that
\beq \label{eqn:a-rig}
| \theta_i - \gamma_i | \leq \frac{N^{\eps/10}}{N^{2/3}  \min \{ i^{1/3}, (N+1-i )^{1/3} \} }
\eeq
with overwhelming probability. We have,
\begin{align}
\left| f(z, t) - \tilf (z, t) \right| \leq 2 |z| \sum_{i=1}^N \int_{\gamma_i}^{\gamma_{i+1}} \left| \frac{ \e^{ \i x} - \e^{ \i \theta_i } }{ ( \e^{ \i x} -z)( \e^{ \i \theta_i } - z) } \right| \rho_t ( x) \d x .
\end{align}
Note the identity,
\beq
| \e^{ \i \varphi } - (1+\eta) \e^{ \i \theta } |^2 = \eta^2 + 4(1+\eta) \sin^2 \left( \frac{ \varphi - \theta}{2} \right),
\eeq
which is useful for estimating from below the quantities in the denominator appearing above. 
WLOG assume $z$ is in the upper half-plane. For the first estimate we break the sum into the terms near $z$ and far from $z$. Let $z = (1+\eta) \e^{ \i \theta}$ and let $i_0$ be the index such that $|\theta_{i_0} - \theta|$ is minimized.  Then,
\begin{align}
  \sum_{ |i - i_0 | < N^{\eps/2} } \int_{\gamma_i}^{\gamma_{i+1}} \left| \frac{ \e^{ \i x} - \e^{ \i \theta_i } }{ ( \e^{ \i x} -z)( \e^{ \i \theta_i } - z) } \right| \rho_t ( x) \d x 
 & \leq  \sum_{ |i - i_0 | < N^{\eps/2} } \int_{\gamma_i}^{\gamma_{i+1}} \left( \frac{1}{ | \e^{ \i x } -z|} + \frac{1}{ | \e^{ \i \theta_i } - z | } \right) \rho_t ( x) \d x \notag\\
& \leq \frac{CN^{\eps/2}}{N \eta}.
\end{align}
For $|i -i_0 | > N^{\eps/2}$ on the event that \eqref{eqn:a-rig} holds we have always that 
\beq
2 | \theta - \theta_i | \geq  | \theta - x|, \qquad  x \in [ \gamma_i, \gamma_{i+1} ].
\eeq
Moreover, one can check that \eqref{eqn:a-rig} implies,
\beq
| \e^{ \i x } - \e^{ \i \theta_i } | \leq  C \frac{N^{\eps/10}}{N \rho_t (x) }, \qquad x \in [ \gamma_i, \gamma_{i+1} ]
\eeq
Therefore,
\begin{align}
 & \sum_{ |i - i_0 | > N^{\eps/2} } \int_{\gamma_i}^{\gamma_{i+1}} \left| \frac{ \e^{ \i x} - \e^{ \i \theta_i } }{ ( \e^{ \i x} -z)( \e^{ \i \theta_i } - z) } \right| \rho_t ( x) \d x  \notag\\
 \leq & \frac{N^{\eps/10}}{N} \sum_{ |i - i_0 | > N^{\eps/2} } \int_{\gamma_i}^{\gamma_{i+1}} \frac{\rho_t(x)}{ \rho_t(x) ( \eta^2 + ( \gamma_i - \theta)^2)} \d x  \leq C \frac{ N^{\eps/10}}{N \eta}
\end{align}
as desired. The second set of estimates is proven similarly. Instead of splitting into two sets of sums, we can directly find for a lower bound for the denominator,
\beq
|( \e^{ \i x } - z)( \e^{ \i \theta_i } -z ) | \geq c( (\eta + \kappa)^2 + ( x - \theta)^2 ), \qquad x \in [ \gamma_i, \gamma_{i+1} ].
\eeq
This yields the claim. \qed

\subsection{Proof of Lemma \ref{lem:edge-sqr}} \label{a:edge-sqr-proof}

By direct calculation with $z= r \e^{ \i \theta}$ we have,
\beq \label{eqn:edge-a2}
- \Re\left[  \tilf (z, t) \right] = (r^2-1) \int_{-\pi}^\pi \frac{1}{ (r-1)^2 + 4 r \sin^2 \left( \frac{x-\theta}{2} \right) } \rho_t (x) \d x.
\eeq
We start with the case that $ \delta < t < 4- \delta$. 
We first assume that $0 < \eta < \eps $ and $0  < \kappa < \eps$ for some small $\eps >0$. Then as long as $\eps >0$ is sufficiently small we have that there are constants $c_1, c_2 >0$ so that the integral is bounded below by,
\beq
\int_{-\pi}^\pi \frac{1}{ (r-1)^2 + 4 r \sin^2 \left( \frac{x-\theta}{2} \right) } \rho_t (x) \d x \geq c_2 \int_{0}^{c_1} \frac{1}{ \eta^2 + (\kappa + x)^2} \sqrt{ x} \d x .
\eeq
We may assume that $\eps < c_1$. For the integral on the RHS, consider first the case that $\kappa \geq \eta$. Than the integral is bounded below by,
\beq
\int_{0}^{c_1} \frac{1}{ \eta^2 + (\kappa + x)^2} \sqrt{ x} \d x \geq \int_{0}^{\kappa} \frac{1}{ \eta^2 + (\kappa + x)^2} \sqrt{ x} \d x \geq c \kappa^{-2} \int_0^{ \kappa} \sqrt{x} \d x \geq c \kappa^{-1/2}.
\eeq
If $\eta \geq \kappa$, then we proceed similarly, by considering the integral over $[0, \eta]$ and bounding the integrand below by $c \sqrt{x} \eta^{-2}$.  This completes the lower bound in the case that $\kappa, \eta < \eps$. For the upper bound in this regime, we again have for a small $c_1$ that, for all $\eps >0$ sufficiently small,
\beq
\int_{-\pi}^\pi \frac{1}{ (r-1)^2 + 4 r \sin^2 \left( \frac{x-\theta}{2} \right) } \rho_t (x) \d x \leq C \int_{0}^{c_1} \frac{1}{ \eta^2 + (\kappa + x)^2} \sqrt{ x} \d x  + C_1,
\eeq
where we used the fact that for $- \Theta_t < x < \Theta_t - c_1$ we have for $\Theta_t \leq \theta \leq \pi$ that
\beq
c_3 < \frac{ | \theta - x|}{ 2} < \pi - c_3
\eeq
for some $c_3 >0$. For the upper bound, first 
\beq
\int_0^{ \kappa + \eta} \frac{ \sqrt{x}}{ \eta^2 + ( \kappa + x)^2} \d x \leq \frac{C}{ ( \kappa + \eta)^2} \int_0^{ \kappa + \eta} \sqrt{x} \d x \leq C (\kappa + \eta)^{-1/2}.
\eeq
On the other hand,
\beq
\int_{\kappa+\eta}^1  \frac{ \sqrt{x}}{ \eta^2 + ( \kappa + x)^2} \d x \leq \int_{\kappa+\eta}^1 x^{-3/2} \d x \leq C (\kappa+\eta)^{-1/2}.
\eeq
This completes the proof of the upper bound when $\kappa + \eta < 2 \eps$ as $C_1 \leq C_1 / \sqrt{\eps}$. In the regime where either $\kappa > \eps$ or $\eta > \eps$ we see that,
\beq
c \leq | z - \e^{ \i x } | \leq C
\eeq
for all $|x| \leq \Theta_t$ by the representation used for the denominator in \eqref{eqn:edge-a2}. This completes the proof of the first estimate.  The estimate for $t < \frac{1}{2}$ follows directly from \eqref{eqn:edge-a2}. \qed

\subsection{Proof of Lemma \ref{lem:cusp-no-eig}} \label{a:cusp-no-eig-proof}

For a contradiction assume that there is an eigenvalue $\lambda = \e^{ \i \theta}$ with  $\Theta_t + N^{-2/3+5\eps} \Delta_t < | \theta | \leq \pi $. As in the proof of Lemma \ref{lem:edge-no-eig} we see that at $ z= (1+\eta ) \e^{ \i \theta}$ we have
\beq
| \Re[ f (z, t) ] | \geq \frac{c}{N \eta}.
\eeq
Take $\eta = N^{-2/3+\eps} \Delta_t^{1/9}$. Since $\kappa \geq N^{-2/3+5 \eps} \Delta_t^{1/9}$ we see that,
\beq
\frac{N^{\eps}}{N \sqrt{ \eta} \sqrt{ \eta + \kappa } } \leq N^{-\eps} \frac{1}{N \eta} .
\eeq
Note that the assumption $\Delta_t \geq N^{-3/4+9 \eps}$ implies that
\beq
\eta = N^{-2/3+\eps} \Delta_t^{1/9} \leq N^{-7\eps} \Delta_t,
\eeq
allowing us to apply the estimate from Lemma \ref{lem:ref-cusp}. Therefore, 
\beq
\Re[ \tilf (z) ] \leq C \frac{ \eta}{ \sqrt{ \kappa + \eta}} \Delta_t^{-1/6}  \leq C N^{-\eps/2} \frac{1}{N \eta} 
\eeq
by our choice of $\kappa$ and $\eta$. This completes the proof. \qed

\section{General $\beta$}

In this section we discuss extensions of our result to general $\beta >0$. In this case, we consider the $N$ particles all on $\rr$,
\beq \label{eqn:general-theta}
\d \theta_i (t) = \sqrt{ \frac{2}{N \beta} } \d W_i (t)  + \frac{1}{2 N} \sum_{ j \neq i } \cot  \left( \frac{ \theta_i(t) - \theta_j (t) }{2} \right) \d t
\eeq
where $\theta_i (0) = 0$ for all $i$ and $\{ W_i \}_i$ are a family of independent standard Brownian motions. We leave aside discussion of the existence of solutions to the above system for now and proceed at a formal level. Setting, $\lambda_i (t) = \e^{ \i \theta_i (t)}$, we find
\beq
\d \lambda_i  = \i \lambda_i \sqrt{ \frac{2}{N \beta} } \d W_i  - \frac{1}{N} \sum_{j \neq i } \frac{ \lambda_i \lambda_j}{ \lambda_i - \lambda_j} \d t - \frac{1}{2} \lambda_i \d t  + \frac{ \lambda_i( \beta-2)}{2 \beta N} \d t
\eeq
recovering the dynamics \eqref{eqn:unit-dbm} when $\beta =2 $. Letting now,
\beq
f (z, t) := \frac{1}{N} \sum_{i=1}^N \frac{ \lambda_i +z}{ \lambda_i - z }
\eeq
as before we find,
\beq \label{eqn:general-f}
\d f(z, t) = - \frac{ z f(z, t)}{2} \del_z f (z, t) + \d M_t^\beta (z) + \frac{ \beta -2}{2 \beta N} \left( z \del_z f (z, t) + z^2 \del_z^2 f (z, t) \right).
\eeq
where 
\beq
\d M^\beta_t (z) = -\sqrt{ \frac{2}{N \beta}} \sum_{i=1}^N \frac{2 \i z \lambda_i (t)}{ ( \lambda_i (t) - z)^2} \d W_i (t)
\eeq
In order to extend our main results, Theorems \ref{thm:bulk}, \ref{thm:edge} and \ref{thm:cusp} to the case of general $\beta$ we must examine how the equation \eqref{eqn:general-f} changes when $\beta \neq 2$. First, the martingale $M^\beta_t (z)$ can be handled the exactly the same as $\beta=2$ as the quadratic variations are the same up to constant factors.

What is at first non-obvious is whether or not the last term on the RHS \eqref{eqn:general-f} drastically changes the behavior of $f$ or not. Recall that our proof strategy is to introduce a stopping time $\tau$ which is the first time an inequality of the form $| f (z_t, t) - \tilf (z_t, t)| \leq N^{\eps/2} \Phi (z_t)$ is violated along some collection of characteristics, where $\Phi(z_t)$ is some control parameter. In the proof of Theorem \ref{thm:bulk} we take $\Phi (z) = 1 / N \eta$ (where we recall the notation $\eta = |z| -1$) and in the other cases $\Phi (z) = B(z)$ with $B(z)$ defined in \eqref{eqn:Bz-def}.

In order for our proofs to go through it suffices to check that for $t < \tau$ along any characteristic $z_t$ we consider that,
\beq
\frac{1}{N} \int_0^t \left( | \del_z f(z_s, s)| + | \del_z^2 f (z_s, s) | \right) \d s \leq N^{\eps/3} \Phi (z_t) \ll N^{\eps/2} \Phi (z_t). 
\eeq
In all three cases (bulk, edge and cusp) we have that along the characteristics we consider we have $N^{\eps} \Phi (z_t) \leq | \Re[ \tilf (z_t, t) ] |$ (see Lemma \ref{lem:edge-ref} and Lemma \ref{lem:cusp-small-error} for the edge and cusp cases; the bulk case is obvious, see e.g., \eqref{eqn:ref-char-bd}) and so we may estimate the integrand on the LHS by
\beq \label{eqn:general-u1}
\frac{1}{N}\left( | \del_z f(z_t, t)| + | \del_z^2 f (z_t, t) |  \right) \leq C \frac{ | \Re[ \tilf (z_t, t) ] |}{ \eta_t \mathrm{dist} (z, \boldsymbol{\lambda} ) }
\eeq
where $\boldsymbol{\lambda}$ denotes the set of eigenvalues. For the bulk case, Theorem \ref{thm:bulk}, we can use simply $ \mathrm{dist} (z, \boldsymbol{ \lambda} ) \geq \eta_t$ and find that the integral of the RHS of \eqref{eqn:general-u1} is bounded above by $(N \eta_t)^{-1}$, which is sufficient.

For the edge regime, Theorem \ref{thm:edge}, one can  modify appropriately the part of the proof of Proposition \ref{prop:edge-bdg} that begins with \eqref{eqn:general-u2} and ends with \eqref{eqn:general-u3} to find an estimate of $C N^{\eps/4} \log(N) (N ( \kappa_t + \eta_t ))^{-1}$ for the integral of the RHS of \eqref{eqn:general-u1} over the time interval $[0, t]$. Similarly, in the cusp regime one can modify the part of the proof of Lemma \ref{lem:cusp-qv} that starts with \eqref{eqn:general-u4} and ends with \eqref{eqn:general-u5} to arrive at the same estimate as in the cusp case.

The conclusion is then that formally the arguments of Sections \ref{sec:bulk}-\ref{sec:cusp} extend to the more general \eqref{eqn:general-f} with only minor changes.  However, all of this rests on the definition of notion of solution of \eqref{eqn:general-theta} which we now briefly discuss.

In the work \cite{cepa2001brownian}, C{\'e}pa and L{\'e}pingle have developed a solution theory for the equation \eqref{eqn:general-theta} and various generalizations. In particular they show that for all choices of $\beta >0$ and all choices of initial data, the equation \eqref{eqn:general-theta} has a unique strong solution. The caveat is that when $\beta <1$, the particles may collide. In such cases, the work \cite{cepa2001brownian} shows that the local time of collisions is sufficiently small so that, for example, the cotangent on the RHS \eqref{eqn:general-theta} is a.s. locally integrable. 

Alternatively, at least in the case $\beta \geq 1$, one can completely avoid discussion of collisions as follows. First, note that if the particles are distinct then one can construct a unique strong solution in an elementary way following, e.g., Chapter 4.3 of \cite{AGZ}, and substituting in the Lyaponuv function used in the proof of Theorem 3.1 of \cite{cepa2001brownian} to rule out collisions. 

Suppose now that we consider \eqref{eqn:general-theta} with initial data $\theta^\eps_i = i \eps$ where eventually we take $\eps \to 0$. As long as, say, $\eps \leq \e^{-N}$, it should not be much trouble to recover all of the results of this paper for the $\theta^\eps$ process, uniformly in $\e^{-N} \geq \eps>0$.

Convergence of the sample paths $\theta^\eps_i (t)$ to some limiting set of sample paths may then be proven using the maximum principle argument of Lemma 2.3 of \cite{huang2019rigidity} which establishes $\sup_{i, 0 < t < T} | \theta^\eps_i(t) - \theta^\delta_i(t)| \leq C_T \sup_i | \theta^\eps_i(0) - \theta^\delta_i(0)|  $. The argument of Lemma 2.3 of \cite{huang2019rigidity} goes through for the process \eqref{eqn:general-theta} after usage of the partial fraction expansion,
\beq
\pi \cot ( \pi  x) = \lim_{n \to \infty} \sum_{|i| \leq n } \frac{1}{ x + i } .
\eeq


\bibliography{mybib}{}
\bibliographystyle{abbrv}

\end{document}